\documentclass[12pt,leqno]{article} 
\usepackage{bezier,amssymb,amsmath}
\input{rk.sty}
\newcommand{\LD}{\Lambda^{\!\mb{\scr\rm D}}}
\newcommand{\sLD}{\Lambda^{\!\mb{\sscr\rm D}}}

\newcommand{\sLDE}{\Lambda^{\!\mb{\rm\sscr D,\,E}}}

\DeclareMathOperator{\Zyk}{Z}

\begin{document}
\title{Some additive galois cohomology rings}
\author{M.\ K\"unzer, H.\ Weber}
\maketitle

\begin{small}
\begin{quote}
\begin{center}{\bf Abstract}\end{center}\vspace*{2mm}
Let $p\geq 3$ be a prime. We consider the cyclotomic extension $\Z_{(p)}[\zeta_{p^2}]\, |\,\Z_{(p)}$, with galois group 
$G = (\Z/p^2)^\ast$. Since this extension is wildly ramified, the $\Z_{(p)} G$-module $\Z_{(p)}[\zeta_{p^2}]$ is not 
projective. We calculate its cohomology ring $\HH^\ast(G,\Z_{(p)}[\zeta_{p^2}];\Z_{(p)})$, carrying the cup product 
induced by the ring structure of $\Z_{(p)}[\zeta_{p^2}]$. Formulated in a somewhat greater generality, our results 
also apply to certain Lubin-Tate extensions.
\end{quote}
\end{small}

\renewcommand{\thefootnote}{\fnsymbol{footnote}}
\renewcommand{\thefootnote}{\arabic{footnote}}

\begin{footnotesize}
\renewcommand{\baselinestretch}{0.7}
\parskip0.0ex
\tableofcontents
\parskip1.2ex
\renewcommand{\baselinestretch}{1.0}
\end{footnotesize}

\setcounter{section}{-1}
\vspace*{-3mm}

\section{Introduction}
\vspace*{-2mm}

\subsection{Results}
\vspace*{-1mm}

\subsubsection{A cohomology ring}
\vspace*{-1mm}

Suppose given a purely ramified extension $T|S$ of discrete valuation rings of residue characteristic $p\geq 3$, with maximal ideals generated by 
$s\in S$ and $t\in T$, respectively. Let $K = \fracfield S$, $L = \fracfield T$, and assume $L|K$ to be galois with $\Gal(L|K) = C_p = \spi{\sa}$.
In particular, $T|S$ is wildly ramified, and so, as {\sc Speiser} remarked {\bf\cite[\rm\S 6]{Sp16}}, $T$ is not projective as a module 
over $SC_p$. We are interested in its cohomology ring $\HH^\ast(C_p,T;S) \aufgl{def} \Ext^\ast_{SC_p}(S,T)$, carrying the 
cup product induced by multiplication on $T$.

We make the assumption that $\val_s(p)\geq b - \ul{b}$, where $b \aufgl{def} -1 + \val_t(t^\sa - t)$ and $b = p\ul{b} + \ol{b}$ 
with $\ol{b}\in [0,p-1]$. This assumption is void if $K$ is of characteristic $p$. As a result, we obtain
\[
\ba{rl}
       & \HH^\ast(C_p,T;S) \\
\defgl & \Ext^\ast_{SC_p}(S,T) \\
  =    & \chi_0^{(0)} S \\
  \ds  & \left(\chi_0^{(1)}(S/s^{\ul{b}+1}S) \ds \cdots \ds \chi_{\ol{b} - 1}^{(1)}(S/s^{\ul{b}+1}S)\right) \ds 
         \left(\chi_{\ol{b}+1}^{(1)}(S/s^{\ul{b}}S) \ds \cdots \ds \chi_{p-1}^{(1)}(S/s^{\ul{b}}S)\right) \\
  \ds  & \chi_0^{(2)}(S/s^{b - \ul{b}}S) \\
  \ds  & \left(\chi_0^{(3)}(S/s^{\ul{b}+1}S) \ds \cdots \ds \chi_{\ol{b} - 1}^{(3)}(S/s^{\ul{b}+1}S)\right) \ds 
         \left(\chi_{\ol{b}+1}^{(3)}(S/s^{\ul{b}}S) \ds \cdots \ds \chi_{p-1}^{(3)}(S/s^{\ul{b}}S)\right) \\
  \ds  & \chi_0^{(4)}(S/s^{b - \ul{b}}S) \\
  \ds  & \cdots\; ,\\
\end{array}
\]
with $S$-linear generators $\chi_i^{(k)}\in \HH^k(C_p,T;S)$; where the multiplication $\chi_0^{(2)}\cup -$ acts as a degree shift by $2$, and where 
multiplication of odd degree elements is given by
\[
\chi_j^{(2l+1)}\cup \chi_k^{(2m+1)} \= \dell_{\,\ol{j+k},\ol{2b}}\, s^{b + \ul{j+k-2b}}\, (\ol{b-j})^{-1} \chi_0^{(2l + 2m + 2)}
\]
for $l,m\geq 0$ and $j,k\in [0,p-1]\ohne\{\ol{b}\}$ (\ref{ThAG11}), where $\dell$ denotes the Kronecker delta.

\subsubsection{A reduction isomorphism}

Let $U|T$ be a further purely ramified extension of discrete valuation rings, and denote $E = \fracfield U$. Suppose $E|K$ to be finite galois, write  
$H = \Gal(E|K)$, and suppose that $[E:L]$ is coprime to $p$. So we are given $K\auf{p}{\tm} L\auf{p'}{\tm} E$, containing $S\tm T\tm U$, 
respectively. We have the reduction isomorphism
\[
\HH^\ast(H,U;S)\;\iso\; \HH^\ast(C_p,T;S)
\]
of graded $S$-algebras (\ref{PropND7}, \ref{PropCup4}), thus enabling us to disregard the upper $p'$-part.

\subsubsection{A cyclotomic example}

The results apply in particular to the following situation.

\begin{tabular}{p{9cm}p{5cm}}
$
\ba{rclcrcl}
S      & = & \Z_{(p)}[\pi_{n-1}]   & \;\; & s      & = & \pi_{n-1} \\
T      & = & \Z_{(p)}[\pi_n]       &      & t      & = & \pi_n \\
U      & = & \Z_{(p)}[\zeta_{p^n}] &      & u      & = & \zeta_{p^n} - 1 \\
H      & = & C_{p-1}\ti C_p        &      &        &   & \\
b      & = & (p^{n-1} - 1)/(p-1)   &      &        &   & \\
\ul{b} & = & (p^{n-2} - 1)/(p-1)   &      & \ol{b} & = & 1 \; ,\\
\ea
$ & 
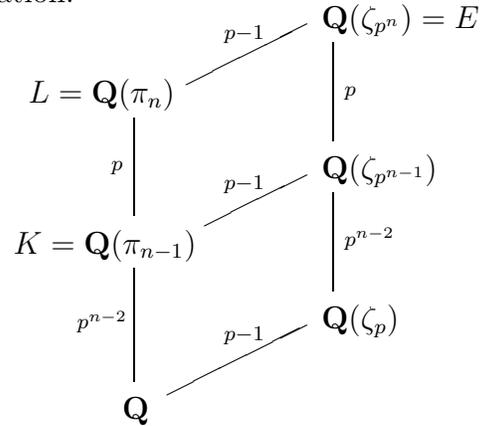
\begin{picture}(500,300)
\put(  10, 200){$L = \Q(\pi_n)$}
\put( 150, 180){\line(0,-1){130}}
\put( 120, 110){$\scm p$}
\put( -10,   0){$K = \Q(\pi_{n-1})$}
\put( 150, -20){\line(0,-1){150}}
\put(  75,-100){$\scm p^{n-2}$}
\put( 135,-220){$\Q$}

\put( 220, 225){\line(2,1){160}}
\put( 270, 285){$\scm p-1$}
\put( 244,  37){\line(2,1){136}}
\put( 270,  85){$\scm p-1$}
\put( 194,-188){\line(2,1){186}}
\put( 270,-115){$\scm p-1$}

\put( 400, 300){$\Q(\zeta_{p^n}) = E$}
\put( 415, 280){\line(0,-1){130}}
\put( 430, 210){$\scm p$}
\put( 400, 100){$\Q(\zeta_{p^{n-1}})$}
\put( 415,  80){\line(0,-1){130}}
\put( 430,  10){$\scm p^{n-2}$}
\put( 400,-100){$\Q(\zeta_p)$}
\end{picture}
\ru{-20}
\end{tabular}

where $n\geq 2$ and 
\[
\pi_n \; \aufgl{def}\; \prod_{j\in [1,\, p-1]} (\zeta_{p^n}^{j^{p^{n-1}}} - 1)\= \Nrm_{E|L}(\zeta_{p^n} - 1)\; ,
\]
cf.\ (\ref{ExCyc2}). In the case $n = 2$, we note that $\pi_1 = p$.

\subsection{Method}

We reinterpret $\HH^\ast(C_p,T;S) \aufgl{def} \Ext^\ast_{SC_p}(S,T) \iso \Ext^\ast_{T\wr C_p}(T,T)$ by adjunction, invoking the 
{\it twisted group ring} $T\wr C_p\om SC_p$ which carries the multiplication $(\rh y)(\ta z) = (\rh\ta)(y^\ta z)$, where $\rh,\ta\in C_p$ and 
$y,\, z\,\in\, T$. In this way, we have gained freedom in our choice of a projective resolution --- when resolving $T$ over $T\wr C_p\,$, we are not 
bound to take a projective resolution of $S$ over $SC_p$ and to tensor it with $T\wr C_p$ over $SC_p\,$. The cup product on 
$\Ext^\ast_{SC_p}(S,T)$ corresponds to the Yoneda product on $\Ext^\ast_{T\wr C_p}(T,T)$ (\ref{PropCup4}).

Still, we need an interpretation of $T\wr C_p$ that facilitates the Yoneda product calculation.
Rationally, there is the Wedderburn isomorphism $L\wr C_p\;\lraisoa{\omega}\;\End_K L = K^{p\ti p}$, sending 
$\rho y\lramapsa{\omega} (z\lramaps z^\rho y)$. Restricting to the locally integral situation and using the resulting Wedderburn {\it embedding} 
$\;T\wr C_p\;\hraa{\omega}\;\End_S T = S^{p\ti p}$, we get a workable isomorphic copy 
$(T\wr C_p)\omega\tm S^{p\ti p}$ of $T\wr C_p$. 

Namely, a matrix in $S^{p\ti p}$ is contained in $(T\wr C_p)\omega$ if and only if it satisfies a set of 
congruences between its entries that can be deduced from the single congruence
\[
t^\sa\;\;\con_{t^{1+b}}\;\; t\;\; ;
\]
see (\ref{ThFT16}, \ref{CorFT14}). With this tie description of $(T\wr C_p)\omega$ inside $S^{p\ti p}$ at our disposal, it is easy to 
calculate the cohomology ring. 

In particular, we believe that it is easier to use this tie description to resolve projectively than to work with the classical projective 
resolution of $T$ over $T\wr C_p\,$. Using the latter, it seems that at some point the operation of $\sa$ on $T$ as 
an element of $\End_S T$, i.e.\ the matrix $(\sa)\omega$, enters the picture; for instance, when solving equations occurring in the resolution 
of cocycles (cf.\ \ref{RemAG4}). Using our tie description, we circumvent this problem by choosing an $S$-linear basis of $(T\wr G)\omega$ 
{\it without} specifying the coefficients of $(\sa)\omega$ therein. So it is no longer necessary to determine the matrix $(\sa)\omega$,
which, indeed, we have not been able to control.

\vspace*{-2mm}
\subsection{Known results and some historical remarks}

\subsubsection{Galois module structure, wild case}

Let $L|K$ be a finite galois extension with $G = \Gal(L|K)$, and let $S\tm K$ be a Dedekind domain with field of fractions
$K$ and integral closure $T$ in $L$. Since $T$ is not isomorphic to $SG$ as a module over $SG$ as soon as a prime ideal of $S$ wildly ramifies in $T$, {\sc Leopoldt} 
split the galois module structure problem in two parts. If $K = \Q$, $S = \Z$ and $G$ is abelian, he determined, firstly, generators for the
{\it associated order}
\[
\Al_{L|K} \= \{\xi\in KG \; :\; T\xi\tm T\}\; ; 
\]
secondly, he showed that $T\iso\Al_{L|K}$ as a module over $\Al_{L|K}$ by construction of an isomorphism {\bf\cite[\rm Satz 6]{Le59}}. 

We recall some of the recent extensions of this result. 
\begin{itemize}
\item[(1)] If $K = \Q(\zeta_n)$, $K\tm L\tm \Q(\zeta_{mn})$ and $S = \Z[\zeta_n]$, then {\sc Byott} and {\sc Lettl} gave a 
description of $\Al_{L|K}$ and showed that $T\iso\Al_{L|K}$ {\bf\cite{BL96}}.
\item[(2)] {\sc Aiba} showed that the analogue of (1) holds in the Carlitz-Hayes function field case if $L$ equals the analogue of $\Q(\zeta_{mn})$ 
and if the analogue of $m$ divides $n$. On the other hand, if this divisibility condition fails, then $T\not\iso\Al_{L|K}$ 
{\bf\cite[\rm th.\ 4]{Ai98}}.
\item[(3)] If $L|K$ is an extension of finite extensions of $\Q_p$ with $L|\Q_p$ abelian, and $S$ the valuation ring of $K$, {\sc Lettl} gave a 
description of $\Al_{L|K}$ and showed that $T\iso\Al_{L|K}$ {\bf\cite{Le98}}.
\item[(4)] As {\sc Byott} observed, similar phenomena as in (2) occur for Lubin-Tate extensions over $\Q_p$ {\bf\cite[\rm th.\ 5.1]{By97}}.
\item[(5)] If $L|K$ is an extension of finite extensions of $\Q_p$, and if $G = C_{p^n}$, {\sc Elder} calculated $T$ as a $\Z_p G$-module 
{\bf\cite{El02}}.
\end{itemize}

It is quite possible that such a description of $T\iso\Al_{L|K}$ as a module over $SG$ might be used to calculate cohomology, and also to calculate 
the cup product on $\HH^\ast(G,T;S)$. Nonetheless, generally speaking, Yoneda products are somewhat easier to calculate than cup products.

In the case of an extension $L|K$ of number fields, with $S$ the ring of algebraic integers in $K$ and $T|S$ at most tamely ramified, 
{\sc Fr\"ohlich} conjectured a connection between the class of $T$ in $\mb{Cl}(\Z G)$ and the Artin root numbers, which has been proven by 
{\sc Taylor}; see e.g.\ {\bf\cite[\rm th.\ 3]{Ta85}}. The extensions of this result from the tame to the wild case start by replacing $T$ by an 
object better suited for class group considerations, thus putting the emphasis on the Artin root numbers; see e.g.\ {\bf\cite[\rm sec.\ 4.2]{CCFT92}}.

\vspace*{-2mm}
\subsubsection{Wedderburn embedding}

Suppose given a Dedekind domain $R$ and an $R$-order $\Omega$ which is rationally semisimple, that is, for which 
$\b \Omega \aufgl{def} (\fracfield R)\ts_R \Omega$ is semisimple. 

The Wedderburn embedding method consists of restricting the Wedderburn isomorphism $\omega$ from $\b \Omega$ to $\Omega$, and to describe the image 
$\Omega\omega$ inside a product of rationally simple maximal orders via congruences of matrix entries, called {\it ties.}

For the first time, this method surfaced in the proof of the {\sc Brauer-Nesbitt} theorem, in which the assumption of the existence of a certain 
non-maximal overorder of a quasiblock is led to a contradiction {\bf\cite[\rm eq.\ (36)]{BN41}}. Here, the {\it quasi}blocks of $\Omega$ are the
$R$-orders $\Omega\eps$ for {\it rational} primitive central idempotents $\eps\in\b \Omega$.

Around the same time, {\sc Higman} calculated the ties for $\Z (C_p\sdp C_{p-1})$ with  
$C_p\sdp C_{p-1} = \spi{a,b \,:\, a^p,\, b^{p-1}, a^b = a^r}$, 
where $r$ is a generator of $\Fu{p}^\ast$ {\bf\cite{Hi40}}. In particular, he calculated the hereditary 
quasiblock of size $(p-1)\ti (p-1)$, which is isomorphic to the twisted group ring $\Z[\zeta_p]\wr C_{p-1}$. See {\bf\cite[\rm sec.\ 6]{Sa83}}.

In the commutative case, descriptions of suborders via congruences have also been given by {\sc Leopoldt} 
{\bf\cite[\rm p.\ 125, p.\ 134, p.\ 140]{Le59}}. 

{\sc Milnor} gave a pullback description of $\Z C_{p^n}$, the iteration of which yields the ties describing this ring 
{\bf\cite[\rm p.\ 601 f.]{Ba68}}. 

The first systematic study to describe Wedderburn embeddings via ties has been undertaken
by {\sc Plesken} {\bf\cite{Pl80}}, {\bf\cite{Pl83}}, following hints of {\sc Zassenhaus.} This has been particularly successful
when the endomorphism rings of the indecomposable projective modules of the quasiblocks under consideration are isomorphic to the discrete
valuation ring $R$.

We list some subsequent work related to the twisted group ring $T\wr G$, where $T|S$ is a finite galois extension of discrete valuation rings with
$G = \Gal(T|S)$.
\begin{itemize}
\item[(1)] {\sc Auslander} and {\sc Rim} showed that $T|S$ is tamely ramified if and only if $T\wr G$ is hereditary {\bf\cite[\rm prop.\ 3.5]{AR63}}, 
which in turn corresponds to an upper triangular shape of its single quasiblock.
\item[(2)] {\sc Benz} and {\sc Zassenhaus} showed that in the wildly ramified case the radical idealisator process, started with $T\wr G$, ends 
up with a hereditary overorder that is unique with respect to certain conditions {\bf\cite[\rm Satz (a)]{BZ85}}.
\item[(3)] {\sc Cliff} and {\sc Weiss} showed that the radical idealisator process of (2) ends up with the unique minimal hereditary overorder of 
$T\wr G$, they determined its shape, and they calculated the number of steps of this process in terms of the different and the ramification index of 
$T|S$ {\bf\cite{CW86}}. As an example, they figured out the ties for $T\wr G$ in the case of a wildly ramified quadratic extension 
{\bf\cite[\rm p.\ 98, 97]{CW85}}.
\item[(4)] {\sc Wingen} calculated the ties of the associated order for $G$ cyclic {\bf\cite[\rm p.\ 82]{Wi93}}.
\item[(5)] {\sc Weber} calculated the ties of $T\wr G$ in the case $T = \Z_{(p)}[\zeta_{p^n}]$, where $n\ge 2$, and $G = C_p$ {\bf\cite{We03}}.
\end{itemize}

\vspace*{-3mm}
\subsection{Acknowledgements}
We would like to thank {\sc G.\ Nebe} for the method of reduction by blockwise decomposition.
We would like to thank {\sc B.\ Keller} for the argument that shows that Yoneda product and cup product coincide.
We would like to thank several colleagues for comments and suggestions.

The first author would like to thank {\sc P.\ Littelmann} for kind hospitality in Strasbourg, 
where the first part of this work was done, and where he received support from the EU TMR-network `Algebraic Lie Representations', grant no.\ 
ERB FMRX-CT97-0100. Another part of the work was financed by a DFG research grant.

\subsection{Leitfaden}
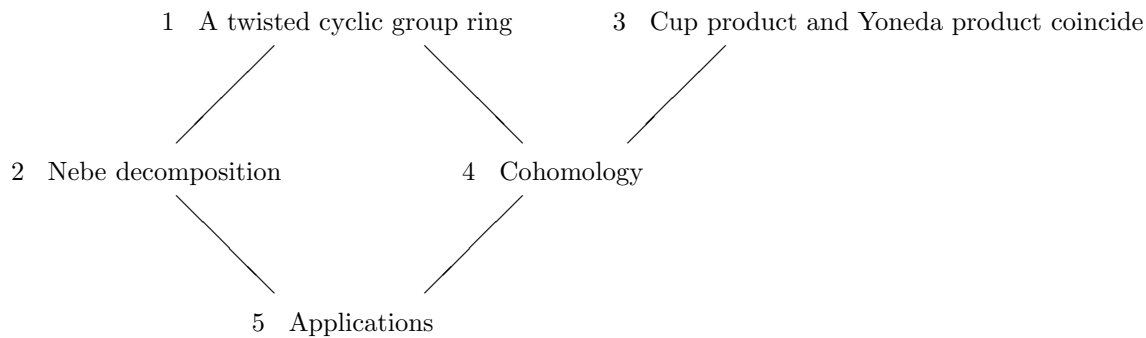
\begin{footnotesize}
\begin{center}
\begin{picture}(1500,450)
\put( 200, 400){\ref{SecTies}\hspace*{2mm} A twisted cyclic group ring}
\put( 800, 400){\ref{SecCupYoneda}\hspace*{2mm} Cup product and Yoneda product coincide}
\put( 350, 380){\line(-1,-1){130}}
\put( 550, 380){\line(1,-1){130}}
\put( 950, 380){\line(-1,-1){130}}
\put(   0, 200){\ref{SecNebe}\hspace*{2mm} Nebe decomposition}
\put( 600, 200){\ref{SecAddGal}\hspace*{2mm} Cohomology}
\put( 220, 180){\line(1,-1){130}}
\put( 680, 180){\line(-1,-1){130}}
\put( 320,   0){\ref{SecAppl}\hspace*{2mm} Applications}
\end{picture}
\end{center}
\end{footnotesize}
\vspace*{-3mm}

\subsection{Notations and conventions}

\begin{footnotesize}
\begin{itemize}
\item[(i)] Composition of maps is written on the right, $\lraa{a}\lraa{b} = \lraa{ab}$. Exception is made for `standard' maps, such as traces, 
characters, derivations \dots
\item[(ii)] For $a,b\in\Z$, we denote by $[a,b] := \{c\in\Z\;|\; a\leq c\leq b\}$ the integral interval.
\item[(iii)] Given elements $x, y$ of some set $X$, we let $\dell_{x,y} = 1$ in case $x = y$ and $\dell_{x,y} = 0$ in case $x\neq y$.
\item[(iv)] Given $a,b\in\Z$, the binomial coefficient $\smatze{a}{b}$ is defined to be zero unless $0\leq b\leq a$.
\item[(v)] If $R$ is a discrete valuation ring with maximal ideal generated by $r$, we write $\val_r(x)$ for the valuation of $x\in R\ohne\{ 0\}$ 
at $r$, i.e.\ $x/r^{\val_r(x)}$ is a unit in $R$. Moreover, we let $\val_r(0) = +\infty$.
\item[(vi)] If $R$ is a discrete valuation ring and $M\hraa{f} N$ an injective $R$-linear map between $R$-modules $M$ and $N$ with cokernel of 
finite length in the sense of Jordan-H\"older, we refer to this length $\ell_R(N/M)$ as the $R$-linear colength of $M$ in $N$.
\item[(vii)] Let $n\geq 1$, let $A$ be a ring. The ring of $n\ti n$-matrices over $A$ is denoted by $A^{n\ti n}$. 
\item[(viii)] Given a ring $A$, by an $A$-module we mean a finitely generated right $A$-module, unless specified otherwise.
\item[(ix)] Given a ring $A$, we denote by $K^-(A)$ the homotopy category of complexes of $A$-modules bounded to the right.
\item[(x)] Given a category $\Cl$, and objects $X$, $Y$ in $\Cl$, we denote the set of morphisms from $X$ to $Y$ by $\liu{\Cl}{(X,Y)}$. If 
$\Cl = \modr\! A$ for a ring $A$, we abbreviate $\liu{A}{(X,Y)} := \liu{\modr\!\! A}{(X,Y)}$.
\end{itemize}

\end{footnotesize}
\section{A twisted cyclic group ring}
\label{SecTies}

\bq
 We give a description of a certain twisted group ring of a cyclic group as a subring of a matrix ring. We proceed in descending generality, ending
 up with a complete description in the case of a cyclic group of prime order as the galois group of a purely ramified extension of discrete valuation 
 rings.  For an attempt in the next larger case $C_{p^2}$, see appendix \ref{AppExp}.

 For the sake of illustration, a continuing example is included, indicated by (cont.).
\eq

\subsection{Subrings defined by derivations}
\label{SubsecDer}

Let $A$ be a ring, let $I\tm A$ be an ideal. Let $k\geq 1$, let $x = (x_j)_{j\in [1,k]}$ be a tuple of elements of $A$ and let the 
corresponding inner derivations be denoted by
\[
\barcl
A & \lrah{D_{x_j}} & A \\
y & \lramaps       & y x_j - x_j y\; . \\
\ea
\]
We note that if $x_i x_j = x_j x_i$, then $D_{x_i}\circ D_{x_j} = D_{x_j}\circ D_{x_i}$.

\begin{Lemma}
\label{LemFT0_5}
Let $h = (h_j)_{j\in [1,k]}$ be a tuple of positive integers (the {\rm height}). Let 
$l = (l_j)_{j\in [1,k]}$ be a tuple of positive integers (the {\rm length}). Then the abelian subgroup
\[
A(x,h,l)_I := \Big\{ y\in A\; :\; D_{x_1}^{i_1}\circ\cdots\circ D_{x_k}^{i_k}(y)\in I^{i_1 l_1 + \cdots + i_k l_k} 
\;\mb{\rm\ for $(i_j)_{j\in [1,k]}\in \prodd{j\in [1,k]} [0,h_j]$} \Big\} \;\tm\; A
\]
is a subring of $A$. 

\rm
{\it Proof.} Given $y$ and $z$ in $A(x,h,l)_I$, we need to show that the product $yz$ is contained in $A(x,h,l)_I$. So suppose given 
$(i_j)_{j\in [1,k]}\in \prod_{j\in [1,k]} [0,h_j]$. The term
\[
D_{x_1}^{i_1}\circ\cdots\circ D_{x_k}^{i_k}(yz)
\]
equals a sum over terms of type
\[
\Big(D_{x_1}^{i'_1}\circ\cdots\circ D_{x_k}^{i'_k}(y)\Big) \cdot \Big(D_{x_1}^{i''_1}\circ\cdots\circ D_{x_k}^{i''_k}(z)\Big) 
\]
where $(i'_j)_j, (i''_j)_j \in \prod_{j\in [1,k]} [0,h_j]$ such that $(i'_j)_j + (i''_j)_j = (i_j)_j$; each such summand is contained in 
$I^{i_1 l_1 + \cdots + i_k l_k}$.\qed
\end{Lemma}

\subsection{A Wedderburn embedding}

Let $T|S$ be a finite and purely ramified extension of discrete valuation rings of residue characteristic $p\geq 3$. Let $L|K$ be the corresponding 
extension of fields of fractions, assumed to be galois with galois group $G$ of order $g$. Let $t$ generate the maximal ideal of $T$, let 
$s := (-1)^{g+1}\Nrm_{L|K}(t)$ generate the maximal ideal of $S$. Let 
\[
\Dfk_{T|S}^{-1} \; :=\; \{ x\in L : \Tr_{L|K}(xT)\tm S\}\;\tm\; L
\] 
define the different ideal $\Dfk_{T|S}\tm T$, and let $\Delta_{T|S} := \Nrm_{L|K}(\Dfk_{T|S})\tm S$ denote the \mb{discriminant} ideal of $T|S$. Write the 
minimal polynomial of $t$ over $K$ as
\[
\mu_{t,K}(X) \; =\; X^g + \left(\sum_{j\in [1,g-1]} e_j X^j\right) - s\;\in S[X]\; .
\]
We recall that $T = S[t]$, that $Tt^g = Ts$, that $S/Ss \lraiso T/Tt$, that $\mu_{t,K}(X)\con_s X^g$ and that
$\Dfk_{T|S} = (\mu'_{t,K}(t))\tm T$ {\bf\cite[\rm III.\S 6, cor.\ 2]{Se62}}.

By $\;T\wr G = \{ \sum_{\rho\in G} \rho y_\rho\; :\; y_\rho\in T\}\;$ we denote the {\it twisted group ring} carrying the multiplication induced by 
$(\rho y)(\tau z) = (\rho\tau)(y^\tau z)$, where $\rho,\,\tau\,\in\, G$ and $y,\,z\,\in\, T$. Let $\Xi$ denote the image of the Wedderburn embedding 
(of $S$-algebras)
\[ 
\barcl
T\wr G & \lraa{\omega} & \End_S T \; =:\; \Gamma    \\
y      & \lramaps      & (\dot y : x\lramaps xy) \\
\rho   & \lramaps      & (\dot\rho : x\lramaps x^\rho) \\
\ea
\]
We consider the subring
\[
\Lambda \; :=\; \Big\{ f\in\Gamma \; :\; (T t^i) f \tm T t^i  \mb{\rm\ for all } i\geq 0\Big\} \;\tm\; \Gamma
\]
which contains $\Xi$, i.e.\ 
\[
T\wr G \;\lraisoa{\omega}\; \Xi \; \tm \; \Lambda\; \tm \Gamma\; ,
\]
by a slight abuse of the notation $\omega$. Since $\dim_L L\ts_T\Lambda = g$, we have  
$\ell_T(\Lambda/\dot t\Lambda) = \ell_T(\Lambda/\Lambda\dot t) = g$, and so we obtain 
\[
\dot t\Lambda \; =\; \Lambda \dot t \; =\; \{ f\in\Gamma \; :\; (Tt^i) f \tm Tt^{i+1} \mb{ for all } i\geq 0\}\; ,
\]
which is thus an ideal of $\Lambda$. Moreover, note that we may write
\[
\Lambda \; =\; \{ f\in\Gamma \; :\; (T t^i) f \tm T t^i  \mb{\rm\ for all } i\in [0,g-1]\}\; .
\]

\begin{Remark}
\label{RemFT3}\rm
Using matrices with respect to the $S$-linear basis $(t^0,t^1,\dots,t^{g-1})$ of $T$ to represent elements of $\Gamma = \End_S T$, the subring 
$\Lambda$ of the full matrix ring $S^{g\ti g}$ is given by the set of matrices with strictly lower triangular entries contained in $Ss$. The ideal 
$\dot t\Lambda$ is given by the set of matrices with lower triangular entries contained in $Ss$, including the main diagonal. In this interpretation, 
we shall refer to matrix positions using the coordinates $[0,g-1]\ti [0,g-1]$.
\end{Remark}

\begin{Remark}[\rm cf.\ {\bf\cite[\rm prop.\ 3.5]{AR63}}]
\label{RemFT3_5}
We have $g\not\con_p 0$ if and only if
\[
T\wr G \;\lraisoa{\omega}\; \Xi \= \Lambda\; .
\]
In particular, in this case $T$ is projective as a module over $T\wr G$ {\rm (cf.\ {\bf\cite[\rm prop.\ 2.4]{K01}}).}

\rm
{\it Proof.} By {\bf\cite[\rm IV.\S 2, cor.\ 1, cor.\ 3]{Se62}}, $g\not\con_p 0$ is equivalent to $\val_t(t^\rho - t) = 1$ for all $\rho\in G\ohne\{ 1\}$.

Now by {\bf\cite[\rm Cor 2.17]{K01}}, the $S$-linear colength of $\Xi$ in $\Gamma$ is given by $\ell_S(\Gamma/\Xi) = g\cdot\val_s(\Delta_{T|S})/2$. 
On the other hand, $\ell_S(\Gamma/\Lambda) = g(g-1)/2$. Therefore, the inclusion $\Xi\tm\Lambda$ is an equality if and only if 
$\val_s(\Delta_{T|S}) = g - 1$. But 
\[
\val_s(\Delta_{T|S}) \= g^{-1} \sumd{\substack{\rho,\,\tau\,\in\, G,\;\\ \rho\;\neq\;\tau}} \val_t(t^\rho - t^\tau)
\= \sumd{\rho\,\in\, G \ohne\{ 1\}} \val_t(t^\rho - t)\; . 
\]
\qedfl{-100}
\end{Remark}

\bq
 \begin{Example}[cont.]
 \label{ExHelp1}\rm
 Let $S = \Z_{(3)}$, $s = 3$, $t = \pi_2 := (\zeta_9 - 1)(\zeta_9^{-1} - 1)$, and $T = \Z_{(3)}[\pi_2]$ (cf.\  
 sec.\ \ref{SubsecCNF}). Let $\sa:\zeta_9\lramaps\zeta_9^4$, restricted from $\Q(\zeta_9)$ to $T$. We have $G = \{1,\sa,\sa^2\}\iso C_3$. We shall 
 use the matrix interpretation explained in (\ref{RemFT3}).

 The Wedderburn embedding sends $t$ to $\dot t = \rsmateckdd{0}{1}{\; 0}{0}{0}{1}{3}{-9}{6}$, the last row resulting from\linebreak 
 $t^3 = 3 - 9 t + 6 t^2$.
 Furthermore, it sends $\sa$ to $\dot\sa = \smateckdd{1}{0}{0}{6}{-5}{1}{24}{-21}{4}$, the second row resulting from $t^\sa = 6 - 5 t + t^2$.
 In particular, $t^\sa - t = t^2 - 6t + 6$ has valuation $2$ at $t$, so $\Xi\tme\Lambda$ is a proper subring, i.e.\ the Wedderburn embedding does not 
 induce an isomorphism of $T\wr G$ with $\Lambda$.

 We have $\Lambda = \smateckdd{S}{S}{S}{3S}{S}{S}{3S}{3S}{S}$, containing $\dot\sa$, and 
 $\dot t\Lambda = \Lambda\dot t = \smateckdd{3S}{S}{S}{3S}{3S}{S}{3S}{3S}{3S}$.
 \end{Example}
\eq

\subsection{The intermediate ring $\LD$}
\label{SubsecInter}

We assume that
\[
b \; :=\; -1 + \min_{\rh\,\in\, G}\val_t(t^\rh - t) \; \geq\; 1\; .
\]
For example, if $G = \spi{\sa} \iso C_g$, then $b = -1 + \val_t(t^\sa - t)$ since $\val_t(t^{\sa^i} - t) \geq \val_t(t^\sa - t)$ for $i\in [0,g-1]$.

Let
\[
\LD \; :=\; \Lambda\left((\dot t),(g-1),(1+b)\right)_{\dot t\Lambda} \;\tm\; \Lambda,
\]
where $k = 1$ in the notation of (\ref{LemFT0_5}). Explicitly, we have
\[
\LD \=
\left\{ u\in\Lambda \; :\;\sumd{j\in [0,i]} (-1)^j \smatze{i}{j} \dot t^j u \dot t^{i-j} 
\con_{\dot t^{(1+b)i}\Lambda} 0 \mbox{ for all } i\in [0,g-1] \right\}\; .
\]

\begin{Lemma}
\label{LemFT4}
We have
\[
\Xi\;\tm\; \LD \;\;\; (\;\tm\;\Lambda\;\tm\;\Gamma\;)\;.
\]

\rm
{\it Proof.} Since $\LD$ is a subring of $\Lambda$, it suffices to show that $\dot t \in \LD$, which follows from  $D_{\dot t}^i(\dot t) = 0$ for $i\geq 1$, and 
that $\dot\rh\in\LD$ for $\rh\in G$, i.e.\ that $D_{\dot t}^i(\dot\rh)$ lies in $\dot t^{i(b+1)}\Lambda$ for $\rh\in G$ and $i\in [0,g-1]$. But
\[
\ba{rll}
0 & \con_{t^{(1+b)i}} & (t^\rh - t)^i \\
  & =                 & \sumd{j\in [0,i]} (-1)^{i-j} \smatze{i}{j} (t^j)^\rh t^{i-j}\\
\ea
\]
implies
\[
\ba{rll}
0 & \con_{\dot t^{(1+b)i}\Lambda} & \sumd{j\in [0,i]} (-1)^{i-j} \smatze{i}{j} (\dot t^j)^{\dot\rh} \dot t^{i-j} \\
  & =                             & \dot\rh^{-1}\left(\sumd{j\in [0,i]} (-1)^{i-j} \smatze{i}{j} \dot t^j \dot\rh \dot t^{i-j}\right) \\
  & =                             & \dot\rh^{-1}\, D_{\dot t}^i(\dot\rh)\; . \\
\ea
\]
\qedfl{-50}
\end{Lemma}

\subsection{A description of $\LD$ via ties}
\label{SubsecApprox}

For $i\in\Z$, we write 
\[
i =: g\ul{i} + \ol{i}\Icm\mb{with $\ol{i}\in [0,g-1]$.}
\]
We note that $\ul{-i} = -1 - \ul{i-1}$, and that $\ol{-i} = g-1 - \ol{i-1}$.

For $i\in [0,g-1]$ and $j\geq 0$, we consider the element $\eps_{i,j}\in\Lambda\tm\Gamma$ defined on the $S$-linear basis $(t^0,\dots,t^{g-1})$ of 
$T$ by
\[
\fbox{$\;\;\ru{-2}
(t^l)\eps_{i,j} := \dell_{i,l}t^{\ol{l + j}} s^{\ul{l + j}}
\;\;$}
\]
for $l\in [0,g-1]$. In the interpretation of (\ref{RemFT3}), $\eps_{i,j}$ is the matrix with a single non-zero entry $s^{\ul{i+j}}$ at position 
$(i,\ol{i+j})$. The tuple $(\eps_{i,j})_{i,j\in [0,g-1]}$ forms an $S$-linear basis of $\Lambda$. For instance,
\[
\dot t\; =\; 
\left(\sum_{i\in [0,g-1]} \eps_{i,1}\right) - s^{-1}\left(\sum_{j\in [1,g-2]} e_j\eps_{g-1,j+1} \right) - e_{g-1}\eps_{g-1,0}\;\in\;\Lambda\; . 
\]

\bq
 The elements $\dot\rh$ seem to be harder to describe in this manner (cf.\ \ref{ExHelp1}). In case $G = C_p$, we shall circumvent this problem by 
 giving a $S$-linear basis of $\Xi$ {\it without} giving the transition matrix from the `initial' $S$-linear basis 
 $(\dot t^i\dot\sa^j)_{i,j\in [0, g-1]}$ of $\Xi$ to it. To do so, we use a description of $\Xi$ by congruences, or {\it ties,} between the 
 coefficients with respect to the $(\eps_{i,j})_{i,j}$-basis of $\Lambda$.
\eq

\begin{Lemma}
\label{LemFT11}
For $i,i'\in [0,g-1]$ and $j,j'\geq 0$, we obtain
\[
\fbox{$\;\;
\eps_{i,j}\eps_{i',j'} = \dell_{i',\ol{i+j}}\eps_{i,j+j'}\; .
\;\;$}
\]

\rm
{\it Proof.} This follows by evaluation on $t^l$ for $l\in [0,g-1]$.\qed
\end{Lemma}

Let
\[
\ddot t \; :=\; \sum_{i\in [0,g-1]} \eps_{i,1} \;\in\;\Lambda\; . 
\]

\begin{Lemma}
\label{LemFT12}
For $j\geq 0$, we obtain
\[
\ddot t^{\, j} \= \sumd{i\in [0,g-1]} \eps_{i,j}\; .
\]

\rm
{\it Proof.} This follows using induction on $j$ together with (\ref{LemFT11}).\qed
\end{Lemma}

\begin{Lemma}
\label{LemFT8}
We have $\;\ddot t\Lambda \= \Lambda \ddot t \= \dot t\Lambda \= \Lambda\dot t\;$. An $S$-linear basis of $\ddot t^{\, k}\Lambda$ for $k\geq 0$ is 
given by 
\[
(s^{-\ul{j-k}}\eps_{i,j})_{i,j\in [0,g-1]}\; .
\]

\rm
{\it Proof.} To see this, we may use the matrix interpretation of (\ref{RemFT3}). \qed
\end{Lemma}

\begin{Assumption}
\label{AssFT5}\rm
Assume that $\val_s(e_j)\geq 1 + b - \ul{b+j}$ for $j\in [1,g-1]$. 
\end{Assumption}

For example, if $b = 1$, then the assumption (\ref{AssFT5}) reads $e_j\in Ss^2$ for $j\in [1,g-2]$ and $e_{g-1}\in Ss$. 

\begin{Lemma}
\label{LemFT9}
Assuming {\rm (\ref{AssFT5}),} we get $\;\ddot t - \dot t\;\in\;\ddot t^{\, 1+b(g-1)}\Lambda\;$.

\rm
{\it Proof.} We have $\ddot t - \dot t = e_{g-1}\eps_{g-1,0} + s^{-1}\sum_{j\in [1,g-2]} e_j\eps_{g-1,j+1}$. Thus, by (\ref{LemFT8}), we need to prove the 
inequalities 
\[
\left\{
\ba{rclcll}
\val_s(e_j)     & \geq & -\; \ul{(j+1) - (1 + b(g-1))} + 1 & = & 1 + b - \ul{j+b}     & \mb{for $j\in [1,g-2]$}  \\
\val_s(e_{g-1}) & \geq & -\; \ul{0 - (1 + b(g-1))}         & = & 1 + b - \ul{(g-1)+b}\;\; , \\ 
\ea
\right.
\]
which in turn have been assumed in (\ref{AssFT5}).\qed
\end{Lemma}

\bq
 Using (\ref{LemFT9}), we are in position to substitute $\dot t$ by $\ddot t$ in the construction of the subring $\LD$.
\eq

\begin{Lemma}
\label{LemFT10}
Assume {\rm (\ref{AssFT5}).} Given $\ga\geq 0$, we obtain
\[
\barcl
\dot t^\ga\LD 
& \aufgl{\rm (1)} & \dot t^\ga\Lambda\!\left((\dot t),(g-1),(1+b)\right)_{\dot t\Lambda} \\
& \aufgl{\rm (2)} & \left\{ u\in\Lambda  \; :\; \sumd{h\in [0,l]} (-1)^h \smatze{l}{h} \dot t^h u \dot t^{l-h}
                                                     \con_{\dot t^{(1+b)l + \ga}\Lambda} 0 \mbox{\rm\ for all } l\in [0,g-1] \right\}\vspace*{2mm} \\
& \aufgl{\rm (3)} & \left\{ u\in\Lambda  \; :\; \sumd{h\in [0,l]} (-1)^h \smatze{l}{h} \ddot t^{\, h} u \ddot t^{\, l-h}
                                                                 \con_{\ddot t^{(1+b)l + \ga}\Lambda} 0 \mbox{\rm\ for all } l\in [0,g-1] \right\} \\
& \aufgl{\rm (4)} & \ddot t^\ga\Lambda\!\left((\ddot t),(g-1),(1+b)\right)_{\ddot t\Lambda}\; . \\
\ea
\]

\rm
{\it Proof.} Equation (1) follows by definition of $\LD$. Let us prove (3). First, we remark that for $i\geq 0$, (\ref{LemFT9}, \ref{LemFT8}) give
\[
\dot t^i \;\con_{\ddot t^{i+b(g-1)}\Lambda}\; \dot t^{i-1}\ddot t \;\con_{\ddot t^{i+b(g-1)}\Lambda}\; \cdots 
\;\con_{\ddot t^{i+b(g-1)}\Lambda}\; \ddot t^{\, i}\; .
\]
Now, both sets are subsets of $\dot t^\ga\Lambda = \ddot t^\ga\Lambda$, as we see by putting $l = 0$.
So suppose given $u\in\dot t^\ga\Lambda = \ddot t^\ga\Lambda$. Using (\ref{LemFT9}, \ref{LemFT8}) we obtain 
\[
\dot t^h u \dot t^{l-h} \con_{\ddot t^{l + b(g-1) + \ga}\Lambda} \ddot t^h u \dot t^{l-h} 
\con_{\ddot t^{l + b(g-1) + \ga}\Lambda} \ddot t^h u \ddot t^{l-h}\; .
\]
for $0\leq h\leq l\leq g-1$. Since $l + b(g-1) + \ga \geq (1+b)l + \ga$, this shows (3). Let us prove (4). Again, both sides are contained in 
$\ddot t^\ga\Lambda$. So if $v\in\Lambda$, then $\ddot t^\ga v$ is in the right hand side of (4) if and only if 
$v\in \Lambda\left((\ddot t),(g-1),(1+b)\right)_{\ddot t\Lambda}$, for multiplication with $\ddot t$ is injective. This in turn is the case if
and only if $\ddot t^\ga v$ is contained in the left hand side of (4), again by injectivity of the multiplication with $\ddot t$. This proves (4).
The argument for (2) is analogous.\qed
\end{Lemma}

\begin{Lemma}
\label{LemFT12_5}
For $x,y\in\Z$ with $y$ coprime to $g$, we obtain
\[
\sumd{i\in [0,g-1]} \ul{x + iy} \; =\; g^{-1}\sumd{i\in [0,g-1]} \left((x + iy) - (\ol{x+iy})\right) \; =\; x + (g-1)(y-1)/2 \; .
\]
\end{Lemma}

\bq
 Now we shall give a description of $\LD$ via {\it ties;} that is, we give certain congruences between the coefficients with respect to the 
 $S$-linear 
 $\eps$-basis of $\Lambda$ that are necessary, and, taken together, also sufficient for an element of $\Lambda$ to lie in $\LD$. Necessity will 
 follow from (\ref{LemFT4}), whereas sufficiency needs a comparison of colengths.
\eq

\begin{Lemma}
\label{PropFT13}
Assume {\rm (\ref{AssFT5}).} Given $\ga\geq 0$, we get
\[
\ba{rclcl}
\dot t^\ga\LD 
& =   & \Bigg\{\;\sumd{i,j\in [0,g-1]} a_{i,j}\eps_{i,j} & : & a_{i,j}\in S\;\;, \mb{\rm\ and  for all $i,j,l\in [0,g-1]$, we have } \\
&     &               &  & \val_s\!\left(\sumd{h\in [0,l]} (-1)^h \smatze{l}{h} a_{\ol{i+h},j}\right) \;\geq\; 1 + \ul{bl - j - 1 + \ga} \;\Bigg\} \\
& \tm & \Lambda \; .  &  & \\
\ea
\]

\rm
{\it Proof.} The $l$th condition of (\ref{LemFT10}) ($l\in [0,g-1]$) for an element $u = \sum_{i,j\in [0,g-1]} a_{i,j}\eps_{i,j}\in\Lambda$ ($a_{i,j}\in S$) to 
lie in $\dot t^\ga \LD$ reads
\[
\ba{rll}
0
& \;\: \con_{\ddot t^{(1+b)l+\ga}\Lambda} & 
                         \sum_{h\in [0,l]}\;\sum_{i,j\in [0,g-1]} (-1)^h \smatze{l}{h} \ddot t^{\, h} a_{i,j}\eps_{i,j} \ddot t^{\, l-h} \\
& \auf{\mb{\scr (\ref{LemFT12})}}{=}  & 
                           \sum_{h\in [0,l]}\;\sum_{i,j,i',i''\in [0,g-1]} (-1)^h \smatze{l}{h} a_{i,j} \eps_{i',h} \eps_{i,j} \eps_{i'',l-h} \\
& \auf{\mb{\scr (\ref{LemFT11})}}{=}  & 
             \sum_{h\in [0,l]}\;\sum_{i,j,i',i''\in [0,g-1]} (-1)^h \smatze{l}{h} a_{i,j} \dell_{i,\ol{i'+h}}\dell_{i'',\ol{i'+h+j}} \eps_{i',l+j} \\
& \;\: =  & \sum_{j,i'\in [0,g-1]}\left(\sum_{h\in [0,l]} (-1)^h \smatze{l}{h} a_{\ol{i'+h},j}\right) \eps_{i',l+j}\; . \\
\ea
\]
For $a\in S$ and $k\geq 0$, the element $a\eps_{i,j}$ is contained in $\ddot t^{\, k}\Lambda$ if and only if $\val_s(a)\geq -\ul{j - k}$ 
(cf.\ \ref{LemFT8}). Thus we obtain the set of conditions
$$
\val_s\!\left(\sum_{h\in [0,l]} (-1)^h \smatze{l}{h} a_{\ol{i+h},j}\right) \;\geq\; -(\ul{-bl + j - \ga})\; =\; 1 + \ul{bl-j-1+\ga} 
\leqno (\ast_{i,j,l})
$$
for $i,j,l\in [0,g-1]$, as claimed.\qed
\end{Lemma}

\subsection{The cyclic case of order $p$}
\label{SubsecComp}

\bq
 As we have seen, in general we have an inclusion $T\wr G\lraiso\Xi\tm\sLD$, and we dispose of a workable description of $\sLD$. Now we shall
 consider a case in which this inclusion will turn out to be an equality.
\eq

Consider the case $g = p$, i.e.\ assume $G = \spi{\sa} = C_p$ to be of order equal to the residue characteristic $p$ of $S$. Note 
that then $b = -1 + \val_t(t^\sa - t)$, and recall that we have stipulated $b\geq 1$.

\begin{Remark}
\label{RemFT13_1}
We have 
\[
\val_t(\Dfk_{T|S}) \= \val_s(\Delta_{T|S}) \= p^{-1}\val_t\!\left(\prod_{\substack{\rh,\,\ta\,\in\, G\; ,\\ \rh\;\neq\;\ta}} (t^\rh - t^\ta) \right) 
\= \val_t\!\left(\prod_{\rh\in G\ohne\{1\}} (t^\rh - t) \right) \= (p-1)(1+b)
\]
(cf.\ {\bf \cite[\rm V.\S 3, lem.\ 3]{Se62}}).
\end{Remark}

\begin{Remark}
\label{RemFT13_2}
In the present case $G = C_p$, assumption {\rm (\ref{AssFT5})} is fulfilled.

\rm
{\it Proof.} Since $\Dfk_{T|S} = (\mu'_{t,K}(t))$, (\ref{RemFT13_1}) implies in particular that 
\[
t^{(p-1)(1+b)} \; |\; \mu'_{L|K}(t) \= p t^{p-1} + \sum_{j\in [1,p-1]} j e_j t^{j-1}\; , 
\]
i.e.\ that $\val_t(e_j) \geq (p-1)(1+b) - (j - 1)$ for $j\in [1,p-1]$, for the valuations of the summands are pairwise different. But since 
$e_j\in S$, this implies 
\[
\val_s(e_j)\;\geq\; \ul{(p-1)(1+b) - (j - 1) + (p - 1)} \= 1 + b - \ul{b + j}
\]
for $j\in [1,p-1]$, and thus, assumption (\ref{AssFT5}) is fulfilled.\qed
\end{Remark}

\begin{Lemma}
\label{LemFT13_5}
Suppose given $\alpha_0,\dots,\alpha_{p-1}\geq 0$. Consider the $S$-linear submodule
\[
\barcl
M & :=  & 
\left\{ (x_h)_{h\in [0,p-1]}\in S^p 
\; :\; \mb{\rm for all $l\in [0,p-1]$, we have $\val_s\!\left(\sum_{h\in [0,l]} (-1)^h\smatze{l}{h} x_h\right) \geq \alpha_l$}\right\} \vspace*{1mm}\\
  & \tm & S^p\; . \\
\ea
\]
An $S$-linear basis of $M$ is given by the tuple $\left((s^{\alpha_l}\smatze{i}{l})_{i\in [0,p-1]}\right)_{l\in [0,p-1]}$.

\rm
{\it Proof.} Let $\w M$ be the $S$-linear submodule of $S^p$ spanned by $\left((s^{\alpha_{l'}}\smatze{i}{l'})_{i\in [0,p-1]}\right)_{l'\in [0,p-1]}$. Since
\[
\sum_{h\in [0,l]} (-1)^h \smatze{l}{h}\left(s^{\alpha_{l'}}\smatze{h}{l'}\right)\; =\; (-1)^{l'} \dell_{l,l'}s^{\alpha_{l'}}
\]
for $l,l'\in [0,p-1]$, we have $\w M\tm M$. The $S$-linear colength of $M$ in $S^p$ is given by $\sum_{l\in [0,p-1]} \alpha_l$, and so is the 
colength of $\w M$. Hence $\w M = M$.\qed
\end{Lemma}

\begin{Proposition}
\label{CorFT14}
Suppose given $\ga\geq 0$ such that $\;\val_s(p)\;\geq\; b - \ul{b-\ga}\;$. Then
\[
\fbox{$\;\;
\ba{rclcl}
\dot t^\ga\LD 
& =   & \Bigg\{\;\sumd{i,j\in [0,p-1]} a_{i,j}\eps_{i,j} & : & a_{i,j}\in S\;\;, \mb{\rm\ and  for all $j,l\in [0,p-1]$, we have } \\
&     &               &  & \val_s\!\left(\sumd{h\in [0,l]} (-1)^h \smatze{l}{h} a_{h,j}\right) \;\geq\; 1 + \ul{bl - j - 1 + \ga} \;\Bigg\} \\
& \tm & \Lambda \; .  &  & \\
\ea
\;\;$}
\]
An $S$-linear basis of $t^\ga\LD$ is given by
\[
\fbox{$\;\;
(\mu^{(\ga)}_{l,j})_{l,\, j\,\in\, [0,p-1]} 
\; :=\; \left(s^{1 + \ul{bl - j - 1 + \ga}} \sumd{i\in [0,p-1]} \smatze{i}{l} \eps_{i,j}\right)_{\! l,\, j\,\in\, [0,p-1]} \;\;$} \;\; .
\]
The according basis of $\LD$ will also be written $(\mu_{l,j})_{l,j} := (\mu^{(0)}_{l,j})_{l,j}$. The $S$-linear colength of $\LD$ in $\Lambda$ 
is given by $bp(p-1)/2$.

\rm
{\it Proof.} Comparing with (\ref{PropFT13}), we have to show the redundancy of the conditions $(\ast_{i,j,l})$ for $i\in [1,p-1]$, $j,l\in [0,p-1]$. First of 
all, we have $(-1)^h\smatze{p-1}{h} \con_{s^{1+\ul{b(p-1) - j - 1 + \ga}}} 1$ for $h\in [0,p-1]$, since this holds modulo $p$. Hence the validity 
of $(\ast_{i,j,p-1})$ is independent of $i\in [0,p-1]$. By downwards induction on $l\in [0,p-2]$, forming the difference, equivalence of 
$(\ast_{i,j,l})$ and $(\ast_{i+1,j,l})$ ensues from $(\ast_{i,j,l+1})$, where $i\in [0,p-2]$.

Now (\ref{LemFT13_5}), applied to a fixed $j\in [0,p-1]$, yields the $S$-linear basis as claimed. The colength of $\LD$ in $\Lambda$ is given by
\[
\barcl
\sum_{l\in [0,p-1]}\;\sum_{j\in [0,p-1]} (1 + \ul{bl-j-1}) 
& \aufgl{(\ref{LemFT12_5})} & p^2 + \sum_{l\in [0,p-1]} (bl - p) \\
& = & bp(p-1)/2. \\ 
\ea
\]
\qedfl{-50}
\end{Proposition}

\begin{Remark}
\label{RemFT15}\rm
In (\ref{CorFT14}), we can as well fix any $m\in [0,p-1]$ and impose $(\ast_{m,j,l})$ on elements of $\Lambda$ for $j,l\in [0,p-1]$ to describe 
$\dot t^\ga\LD$ inside. Accordingly, we obtain an $S$-linear basis 
$(s^{1 + \ul{bl - j - 1 + \ga}} \sum_{i\in [0,p-1]} \smatze{i}{l} \eps_{\ol{i+m},j})_{l,j}$ of $\dot t^\ga\LD$.
\end{Remark}

\begin{Theorem}
\label{ThFT16}
Recall that we assume $b = -1 + \val_t(t^\sa - t) \geq 1$, where $C_p = \spi{\sa}$. Suppose in addition that $\val_s(p)\geq b - \ul{b}$.

The Wedderburn embedding factors over the $S$-algebra isomorphism
\[
\fbox{$\;\;
T\wr C_p\;\lraisoa{\omega}\;\LD\;\tm\;\Lambda\;\tm\;\Gamma
\;\;$}\;\; .
\]

\rm
{\it Proof.} By (\ref{LemFT4}), it suffices to show that the $S$-linear colengths of $(T\wr G)\omega = \Xi$ and of $\LD$ in $\Lambda$ are equal. The colength of 
$\Lambda$ in $\Gamma$ equals $p(p-1)/2$, so that, by (\ref{CorFT14}), the colength of $\LD$ in $\Gamma$ equals $(1+b)p(p-1)/2$. On the other hand,
by {\bf\cite[\rm cor.\ 2.17]{K01}}, the colength of $\Xi$ in $\Gamma$ equals $p\val_s(\Delta_{T|S})/2 \aufgl{(\ref{RemFT13_1})} p(1+b)(p-1)/2$.\qed
\end{Theorem}

\bq
 \begin{Example}[cont.]
 \label{ExHelp2}\rm
 We have $\dot t = \rsmateckdd{0}{1}{\; 0}{0}{0}{1}{3}{-9}{6}$ and $\ddot t = \smateckdd{0}{1}{0}{0}{0}{1}{3}{0}{0}$.

 We have $g = p = 3$ and $b = 1$. Since $\mu_{t,\Q}(X) = X^3 - 6 X^2 + 9 X - 3$, assumption (\ref{AssFT5}) is satisfied, in accordance with 
 (\ref{RemFT13_2}). Consequently, we have $\dot t - \ddot t = \rsmateckdd{0}{0}{\; 0}{0}{0}{0}{0}{-9}{6} \in 3\Lambda$ (cf.\ \ref{LemFT9}).

 If $\ga \leq 1$, then $1 = \val_s(p)\geq b - \ul{b-\ga} = 1$, and so (\ref{CorFT14}) can be applied to give
 \[
 \hspace*{-5mm}
 \ba{rcll}
 \Z_{(3)}[\pi_2]\wr C_3 
 & =   & T\wr C_3 \\
 & \lraisoa{\omega} &  \LD \\
 & =   & \left\{\ru{-2}\right.
         \rsmateckdd{a_{0,0}}{a_{0,1}}{a_{0,2}}{3 a_{1,2}}{a_{1,0}}{a_{1,1}}{3 a_{2,1}}{3 a_{2,2}}{a_{2,0}} \; :\; 
           & a_{0,0}\in 3^0 S,\; a_{0,0} - a_{1,0}\in 3^1 S,\; a_{0,0} - 2 a_{1,0} + a_{2,0}\in 3^1 S; \\
 &     &   & a_{0,1}\in 3^0 S,\; a_{0,1} - a_{1,1}\in 3^0 S,\; a_{0,1} - 2 a_{1,1} + a_{2,1}\in 3^1 S; \\
 &     &   & a_{0,2}\in 3^0 S,\; a_{0,2} - a_{1,2}\in 3^0 S,\; a_{0,2} - 2 a_{1,2} + a_{2,2}\in 3^0 S;  \left.\ru{-2}\right\} \\
 & =   & \left\{\ru{-2}\right.
         \rsmateckdd{a_{0,0}}{a_{0,1}}{a_{0,2}}{3 a_{1,2}}{a_{1,0}}{a_{1,1}}{3 a_{2,1}}{3 a_{2,2}}{a_{2,0}} \; :\; 
           & a_{0,0}\con_3 a_{1,0}\con_3 a_{2,0} ;\;\; a_{0,1} + a_{1,1} + a_{2,1} \con_3 0 \left.\ru{-2}\right\} \vspace*{1mm}\\
 & \tm & S^{3\ti 3}\; , \\
 \ea
 \]
 where $S = \Z_{(3)}$. In particular, we have $\dot\sa = \smateckdd{1}{0}{0}{6}{-5}{1}{24}{-21}{4}\in\LD$, as expected. Similarly,
 \[
 \ba{rcll}
 \dot t\LD 
 & =   & \left\{\ru{-2}\right.
         \rsmateckdd{a_{0,0}}{a_{0,1}}{a_{0,2}}{3 a_{1,2}}{a_{1,0}}{a_{1,1}}{3 a_{2,1}}{3 a_{2,2}}{a_{2,0}} \; :\; 
           & a_{0,0}\in 3^1 S,\; a_{0,0} - a_{1,0}\in 3^1 S,\; a_{0,0} - 2 a_{1,0} + a_{2,0}\in 3^1 S; \\
 &     &   & a_{0,1}\in 3^0 S,\; a_{0,1} - a_{1,1}\in 3^1 S,\; a_{0,1} - 2 a_{1,1} + a_{2,1}\in 3^1 S; \\
 &     &   & a_{0,2}\in 3^0 S,\; a_{0,2} - a_{1,2}\in 3^0 S,\; a_{0,2} - 2 a_{1,2} + a_{2,2}\in 3^1 S;  \left.\ru{-2}\right\} \\
 & =   & \left\{\ru{-2}\right.
         \rsmateckdd{3a_{0,0}}{a_{0,1}}{a_{0,2}}{3 a_{1,2}}{3a_{1,0}}{a_{1,1}}{3 a_{2,1}}{3 a_{2,2}}{3a_{2,0}} \; :\; 
           & a_{0,1}\con_3 a_{1,1}\con_3 a_{2,1} ; \;\;\;\; a_{0,2} + a_{1,2} + a_{2,2} \con_3 0 \left.\ru{-2}\right\} \vspace*{1mm}\\
 & \tm & S^{3\ti 3}\; . \\
 \ea
 \]
 \end{Example}

 \begin{Example}
 \label{ExFT17}\rm
 Suppose that $S$ contains a primitive $p$-th root of unity $\zeta_p$ (so in particular, $\ch K = 0$). Let $T := S[\sqrt[p]{s}]$ and 
 $t := \sqrt[p]{s}$, i.e.\ let $\mu_{t,K}(X) = X^p - s$. Then $L|K$ is galois with galois group $C_p$, and we have
 \[
 b \= -1 + \val_t(t\zeta_p - t) \= p\val_s(1 - \zeta_p) = \frac{p}{p-1}\val_s(p)\; ,
 \]
 i.e.\ $\val_s(p) = b - \ul{b}$, so that (\ref{ThFT16}) may be applied. 
 \end{Example}
\eq

\begin{Corollary}
\label{CorFT18}
Suppose given two discrete valuation rings $T$ and $T'$ over $S$ with $T|S$ and $T'|S$ both purely ramified and galois with galois group $C_p$. 
Suppose that 
\[
\val_s(\Delta_{T|S}) \= \val_s(\Delta_{T'|S}) \;\leq\; p\val_s(p) + p - 1\; . 
\]
Then
\[
T\wr C_p\; \iso\; T'\wr C_p\; 
\]
as $S$-algebras.

\rm
{\it Proof.} Since $\val_s(\Delta_{T|S}) = (p - 1)(b + 1)$ by (\ref{RemFT13_1}), similarly for $T'$, and since $\LD$ depends only on $S$ and $b$ 
(cf.\ \ref{CorFT14}), we may conclude
\[
T\wr C_p\;\auf{(\ref{ThFT16})}{\iso}\;\LD\;\auf{(\ref{ThFT16})}{\iso}\; T'\wr C_p\; .
\]
\qedfl{-80}
\end{Corollary}

\bq
 \begin{Example}[cont.]
 \label{ExHelp2_5}\rm
 Let $S = \Z_{(3)}$, $T = \Z_{(3)}[\pi_2]$ and let $T' = S[t']$ with $\mu_{t',K}(X) = X^3 + 3X^2 - 18X + 48$ {\bf\cite{JJ01}}. Then 
 $\Delta_{T|S} = \Delta_{T'|S} = Ss^4$, whence $T\wr C_3\iso T'\wr C_3$ by (\ref{CorFT18}).
 The fields of fractions of $T$ and $T'$, however, are not isomorphic.
 \end{Example}
\eq
\section{Nebe decomposition}
\label{SecNebe}

\bq
 The purpose of this section is, translated to the basic example $\Z_{(p^2)}[\zeta_{p^2}]|\Z_{(p)}$, the reduction of the cohomology calculation
 from the galois group $C_p\ti C_{p-1}$ to its quotient $C_p$, being the galois group of a subextension.
\eq

\subsection{A block decomposition}
\label{SubsecBD}

Let $T|S$ and $U|T$ be finite purely ramified extensions of discrete valuation rings with maximal ideals generated by $s\in S$, $t\in T$ and 
$u\in U$, respectively. Let $K = \fracfield S$, $L = \fracfield T$ and $E = \fracfield U$, and suppose that $E|K$ is galois with galois group $H$, 
and that $L|K$ is galois with galois group $G$. In particular, if
\[
1\;\lra\; N\; \lra\; H\;\lra\; G\;\lra\; 1 
\]
is short exact, then $E|L$ is galois with galois group $N$. Denote $h := |H|$, $g := |G|$ and $n := |N|$, so $h = gn$. 

The situation can be depicted as follows.
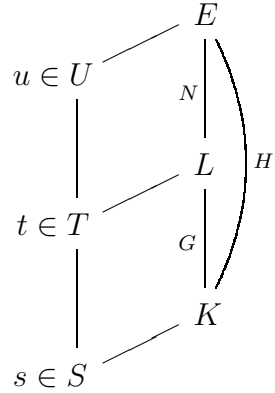
\begin{center}
\begin{picture}(300,520)
\put(  30,   0){$s\in S$}
\put( 115,  50){\line(0,1){130}}
\put(  35, 200){$t\in T$}
\put( 115, 250){\line(0,1){130}}
\put(  30, 400){$u\in U$}

\put( 150,  30){\line(2,1){100}}
\put( 150, 230){\line(2,1){100}}
\put( 150, 430){\line(2,1){100}}

\put( 270,  80){$K$}
\put( 285, 130){\line(0,1){130}}
\put( 250, 180){$\scm G$}
\put( 270, 280){$L$}
\put( 285, 330){\line(0,1){130}}
\put( 250, 380){$\scm N$}
\put( 270, 480){$E$}
\bezier{200}(300,130)(380,295)(300,460)
\put( 350, 290){$\scm H$}
\end{picture}
\end{center}

Suppose given an $S$-linear submodule $V$ of $U$ spanned by a $T$-linear basis of $U$, i.e.\ such that $T\ts_S V\lraisoa{\phi} U$, 
$x\ts y\lramaps xy$. In general, $V$ is not a subring of $U$. We write 
\[
\ba{lcl}
\Gamma     & := & \End_S U \\
\Gamma'    & := & \End_S T \\
\Gamma''   & := & \End_T U \\
\Gamma''_0 & := & \End_S V\; . \\
\ea
\]
Here we change our notation slightly in that the ring $\End_S T$, which has previously been denoted by $\Gamma$, is now denoted by $\Gamma'$. 
Similarly, we denote the Wedderburn embeddings by
\[
\barcl
U\wr H & \hraa{\omega}   & \Gamma \\
T\wr G & \hraa{\omega'}  & \Gamma' \\
U\wr N & \hraa{\omega''} & \Gamma''\; , \\
\ea
\]
again switching notation from what has been denoted by $\omega$ to $\omega'$. Accordingly, we denote the images by
\[
\ba{lcl}
\Xi   & := & (U\wr H)\omega \\
\Xi'  & := & (T\wr G)\omega' \\
\Xi'' & := & (U\wr N)\omega'' \\
\ea
\]
Usage of the notation $\Gamma'$, $\omega'$, $\Xi'$ pertains only to the present section~\ref{SecNebe}, in which we consider the passage from 
$U\wr H$ to $T\wr G$.

Note that we have an isomorphism of $S$-algebras
\[
\ba{rllclcl}
\Gamma' & \ts_S & \Gamma''_0 & \lraiso  & \End_S (T\ts_S V)                & \lraiso  & \Gamma \\
\al     & \ts   & \be_0      & \lramaps & (x\ts y \lramaps x\al\ts y\be_0) &          &        \\
        &       &            &          & \eta                             & \lramaps & \phi^{-1}\eta\phi\; , \\
\ea
\]
denoted by $\Gamma'\ts_S \Gamma''_0\lraisoa{\theta}\Gamma$. Moreover, we have an isomorphism of $T$-algebras
\[
\ba{rllcrllcl}
T & \ts_S & \Gamma''_0 & \lraiso  & \End_T (T & \ts_S & V)    & \lraiso  & \Gamma'' \\
1 & \ts   & \be_0      & \lramaps &         1 & \ts   & \be_0 & \lramaps & \phi^{-1} (1\ts\be_0) \phi \; ,\\
\ea
\]
which shall be denoted by $\;T\ts_S \Gamma''_0 \;\lraisoa{\psi}\; \Gamma''\;$. Thus, we have an $S$-linear isomorphism
\[
\Gamma'\ts_T\Gamma''\;\;\llaisoa{1\ts\psi}\;\; 
\Gamma'\ts_T T\ts_S\Gamma''_0 \;\;\llaiso\;\;\Gamma'\ts_S \Gamma''_0\;\;\lraisoa{\theta}\;\; \Gamma \;\;,
\]
denoted by 
\[
\Gamma'\ts_T\Gamma''\;\;\lraisoa{\tht}\;\; \Gamma\;\; .
\]

\begin{Lemma}
\label{LemIns}
If $(x_1,\dots,x_g)$ is an $S$-linear basis of $T$, and $(y_1,\dots,y_n)$ is a $T$-linear basis of $U$ lying in $V$, and given $\al\in\Gamma'$ and 
$\be\in\Gamma''$, then $(\al\ts\be)\tht$ sends $x_i y_j$ to $(x_i\al)(y_j\be)$, where $i\in [1,g]$, $j\in [1,n]$. 

\rm
{\it Proof.} By $T$-bilinearity we may assume that $\be = (1\ts\be_0)\psi$ with $\be_0\in\Gamma''_0$, in which case the assertion follows by construction.\qed
\end{Lemma}

We use $\tht$ for a transport of the $S$-algebra-structure from $\Gamma$ to $\Gamma'\ts_T\Gamma''$; i.e.\ given 
$\ga_1\, ,\ga_2\,\in\,\Gamma'\ts_T\Gamma''$, we let 
\[
\ga_1\cdot\ga_2 \; := \; ((\ga_1\tht)\cdot(\ga_2\tht))\tht^{-1}\; ,
\]
so that $\tht$ becomes an isomorphism of $S$-algebras.

Now we choose $V$ to be the $S$-linear span of $(u^0,\dots,u^{n-1})$ in $U$. Let 
\[
\Lambda'' \; :=\; \{\phi''\in\Gamma''\; :\; u^k\phi''\in Uu^k \;\mb{ for $k\geq 0$ } \} \;\tm\;\Gamma'' \; ,
\]
which is a sub-$T$-algebra.

\begin{Lemma}
\label{LemND1}
The $S$-linear submodule $\;\Xi'\ts_T\Lambda''\;\tm\; \Gamma'\ts_T\Gamma''\;$ is a sub-$S$-algebra.

\rm
{\it Proof.} Given $k,l\in [0,n-1]$, we let $e''_{k,l}\in\Gamma''$ be defined by $u^j e''_{k,l} = \dell_{j,k} u^l$ for $j\in [0,n-1]$. Then
$$
\Xi'\ts_T\Lambda'' 
\= \left\{\sum_{k,l\in [0,n-1]} \xi'_{k,l} \ts e''_{k,l} \; :\; 
   \mb{$\xi'_{k,l}\in\Xi'$ for $k\leq l$, $\xi'_{k,l}\in \dot t\Xi'$ for $k > l$} \right\}\; .
\leqno (\ast)
$$
Since we can write $e''_{k,l} = (1\ts e''_{k,l;0})\psi$ with $e''_{k,l;0}\in\Gamma''_0$, i.e.\ since $e''_{k,l}$ restricts to an $S$-linear 
endomorphism of $V$, multiplication of elements written in the form as in $(\ast)$ is given by 
\[
(\xi' \ts e''_{k,l})(\w\xi'\ts e''_{\w k,\w l}) \= (\xi'\w\xi') \ts (\dell_{l,\w k} e''_{k,\w l})\; ,
\]
where $\w\xi,\,\w\xi'\,\in\,\Xi'$ and $k,\,l,\,\w k,\,\w l\,\in\, [0,n-1]$. Now if $k > \w l$ and $l = \w k$, then $k > l$ or $\w k > \w l$, hence
$\xi'\w\xi'\in\dot t\Xi'$. Therefore, the $S$-linear submodule $\Xi'\ts_T\Lambda''$ is closed under multiplication.\qed
\end{Lemma}

\begin{Lemma}
\label{LemND2}
We have $\;\Xi\tht^{-1}\;\tm\;\Xi'\ts_T\Lambda''\;$.

\rm
{\it Proof.} Since $\Xi'\ts_T\Lambda''$ is a sub-$S$-algebra of $\Gamma'\ts_T\Gamma''$ (\ref{LemND1}), it suffices to show that 
$u\omega\tht^{-1}\in \Xi'\ts_T\Lambda''$ and that $\rh\omega\tht^{-1}\in \Xi'\ts_T\Lambda''$ for $\rh\in H$.

We have $(1\ts u\omega'')\tht = u\omega$ by (\ref{LemIns}), hence $u\omega\tht^{-1}\in \Xi'\ts_T\Lambda''$.

Suppose given $\rh\in H$. Let $(x_1,\dots,x_g)$ be an $S$-linear basis of $T$. The element $x_i u^j$ is sent to $x_i^\rh (u^\rh)^j$ by $\rh\omega$, 
where $i\in [1,g]$, $j\in [0,n-1]$ (\ref{LemIns}). The element of $\Gamma'$ determined by $x_i\lramaps x_i^\rh$ is just $\rh\omega'$, which is in 
$\Xi'$. The element of $\Gamma''$ determined by $u^j\lramaps (u^\rh)^j$ is contained in $\Lambda''$. The tensor product of these elements is thus
contained in $\Xi'\ts_T\Lambda''$, as was to be shown.\qed
\end{Lemma}

\bq
 \begin{Remark}
 \label{RemIns2}\rm
 In general, the element $u^j\lramaps (u^\rh)^j$, where $j\in [0,n-1]$, is {\it not} contained in $\Xi''$ since $\rh$ need not be in $N$. 
 \end{Remark}
\eq

By (\ref{LemND1}, \ref{LemND2}), we have a commutative diagram of $S$-algebras
\begin{center}
\begin{picture}(350,250)
\put( -20, 200){$\Gamma'\ts_T\Gamma''$}
\put( 150, 210){\vector(1,0){130}}
\put( 200, 225){$\scm\tht$}
\put( 200, 190){$\scm\sim$}
\put( 300, 200){$\Gamma$}
\put(  50,  60){\vector(0,1){120}}
\put(  40,  60){\oval(20,20)[b]}
\put( 315,  60){\vector(0,1){120}}
\put( 305,  60){\oval(20,20)[b]}
\put( -22,   0){$\Xi'\ts_T\Lambda''$}
\put( 270,  10){\vector(-1,0){120}}
\put( 270,  20){\oval(20,20)[r]}
\put( 200,  25){$\scm\tht^{-1}$}
\put( 300,   0){$\Xi\; .$}
\end{picture}
\end{center}

\begin{Proposition}
\label{PropND3}
If $n\not\con_p 0$, then 
\[
\Xi\tht^{-1} \= \Xi'\ts_T\Lambda''\; .
\]
In other words, the Wedderburn embedding $U\wr H\hraa{\omega}\Gamma$ factors as
\[
U\wr H\;\;\lraisoa{\omega\tht^{-1}}\;\; \Xi'\ts_T\Lambda'' \;\hraa\; \Gamma'\ts_T\Gamma'' \;\lraisoa{\tht}\; \Gamma\; .
\]

\rm
The tensor product $\Xi'\ts_T\Lambda''$ is called the {\it Nebe decomposition} of $U\wr H$.

{\it Proof.} First, we remark that if $n\not\con_p 0$, then $\Xi'' = \Lambda''$ by (\ref{RemFT3_5}).

We need to show that the $S$-linear colengths coincide. By {\bf\cite[\rm cor.\ 2.17]{K01}}, the colength of $\Xi\tm\Gamma$ equals
\[ 
h\val_s(\Delta_{U|S})/2 \= h\val_u(\Dfk_{U|S})/2 \= h\val_u(\Dfk_{U|T}\Dfk_{T|S})/2 \= h(n-1)/2 + hn\val_s(\Delta_{T|S})/2\; . 
\]
On the other hand, by {\bf\cite[\rm cor.\ 2.17]{K01}}, the colength of the composite embedding \linebreak
$\Xi'\ts_T\Lambda''\tm\Gamma'\ts_T\Lambda''\tm\Gamma'\ts_T\Gamma''$ equals   
\[
(g\val_s(\Delta_{T|S})/2)\cdot n^2 + g\cdot n(n-1)/2 \= hn\val_s(\Delta_{T|S})/2 + h(n-1)/2\; .
\]
\qedfl{-90}
\end{Proposition}

\bq
 \begin{Example}[cont.]
 \label{ExHelp3}\rm
 Let $S = \Z_{(3)}$, $s = 3$, $T = \Z_{(3)}[\pi_2]$, $t = \pi_2$, $U = \Z_{(3)}[\zeta_{3^2}]$, $u = \zeta_{3^2} - 1$, so that 
 $H = (\Z/3^2)^\ast\iso C_3\ti C_2$, $G = C_3$ and $N = C_2$. We use the \mb{$S$-linear} basis $(u^0,u^1)$ of $V$ and the $S$-linear basis 
 $(t^0,t^1,t^2)$ of $T$, and thus the $S$-linear basis $(u^0 t^0, u^0 t^1, u^0 t^2, u^1 t^0, u^1 t^1, u^1 t^2)$ of $U$.
 By (\ref{PropND3}) and using $(\ast)$, the Wedderburn descriptions of $\Xi'$ and $\dot t\Xi'$ obtained in (\ref{ExHelp2}) can be inserted as 
 blocks into
 \[
 \ba{l}
 U\wr H \;\;\lraisoa{\omega}\;\; \Xi \;\=\; \vs\\
 \left\{
 \left[
 \etwasenger
 \ba{rrr|rrr}
 \scm a_{0,0;0,0}  &\scm   a_{0,0;0,1} &\scm   a_{0,0;0,2} &\scm   a_{0,1;0,0} &\scm   a_{0,1;0,1} &\scm   a_{0,1;0,2} \\
 \scm 3a_{0,0;1,2} &\scm   a_{0,0;1,0} &\scm   a_{0,0;1,1} &\scm  3a_{0,1;1,2} &\scm   a_{0,1;1,0} &\scm   a_{0,1;1,1} \\
 \scm 3a_{0,0;2,1} &\scm  3a_{0,0;2,2} &\scm   a_{0,0;2,0} &\scm  3a_{0,1;2,1} &\scm  3a_{0,1;2,2} &\scm   a_{0,1;2,0} \\\hline
 \scm 3a_{1,1;0,0} &\scm   a_{1,1;0,1} &\scm   a_{1,1;0,2} &\scm   a_{1,0;0,0} &\scm   a_{1,0;0,1} &\scm   a_{1,0;0,2} \\
 \scm 3a_{1,1;1,2} &\scm  3a_{1,1;1,0} &\scm   a_{1,1;1,1} &\scm  3a_{1,0;1,2} &\scm   a_{1,0;1,0} &\scm   a_{1,0;1,1} \\
 \scm 3a_{1,1;2,1} &\scm  3a_{1,1;2,2} &\scm  3a_{1,1;2,0} &\scm  3a_{1,0;2,1} &\scm  3a_{1,0;2,2} &\scm   a_{1,0;2,0} \\
 \ea
 \weiter
 \right]
 \;
 \text{\huge :}
 \;
 \ba{l}
 \left\{
 \ba{l}
 a_{0,0;0,0}\con_3 a_{0,0;1,0}\con_3 a_{0,0;2,0},\\ 
 a_{0,0;0,1} + a_{0,0;1,1} + a_{0,0;2,1}\con_3 0; \\
 \ea 
 \right.\vspace*{1mm} \\
 \left\{
 \ba{l}
 a_{0,1;0,0}\con_3 a_{0,1;1,0}\con_3 a_{0,1;2,0}, \\
 a_{0,1;0,1} + a_{0,1;1,1} + a_{0,1;2,1}\con_3 0; \\
 \ea 
 \right.\vspace*{1mm} \\
 \left\{
 \ba{l}
 a_{1,0;0,0}\con_3 a_{1,0;1,0}\con_3 a_{1,0;2,0}, \\
 a_{1,0;0,1} + a_{1,0;1,1} + a_{1,0;2,1}\con_3 0; \\
 \ea 
 \right.\vspace*{1mm} \\
 \left\{
 \ba{l}
 a_{1,1;0,1}\con_3 a_{1,1;1,1}\con_3 a_{1,1;2,1}, \\
 a_{1,1;0,2} + a_{1,1;1,2} + a_{1,1;2,2} \con_3 0 \\
 \ea 
 \right. \vspace*{1mm}\\
 \ea
 \right\}\vs\\
 \tm\; (S^{3\ti 3})^{2\ti 2}\; . \\
 \ea
 \]
 \end{Example}
\eq

\subsection{A reduction isomorphism}
\label{SubsecRed}
We maintain the notation of section \ref{SubsecBD}, having chosen $V = S\spi{u^0,\dots,u^{n-1}}\tm U$. 

\begin{Definition}
\label{DefND4}\rm
Given a $\Xi'$-module $Y$, the $S$-module $Y\ts_T U$ decomposes into a direct sum 
\[
Y\ts_T U = \Ds_{j\in [0,n-1]} Y\ts u^j,
\]
so that we may endow it with the structure of a $\Xi'\ts_T\Lambda''$-module by means of
\[
(y\ts u^j) (\xi'\ts e''_{k,l}) := y \xi'\ts \dell_{j,k} u^l,
\]
where $y\in Y$, $\xi'\in\Xi'$, $j,k,l\in [0,n-1]$, $\xi'\in\dot t\Xi$ if $k > l$, and where $e''_{k,l}$ is as in $(\ast)$ above. 

More naturally explained, we use the isomorphism
$Y\ts_T U \iso Y\ts_S V$ to transport the structure of an $\Xi'\ts_S\Gamma''_0$-module from $Y\ts_S V$ to $Y\ts_T U$. 
Restricting the $S$-algebra isomorphism $\Gamma'\ts_S\Gamma''_0 \lraiso \Gamma'\ts_T\Gamma''$ to 
$\Xi'\ts_S\Gamma''_0 \lraiso \Xi'\ts_T\Gamma''$, we obtain on $Y\ts_T U$ the structure of a $\Xi'\ts_T\Gamma''$-module, which we restrict to  
$\Xi'\ts_T\Lambda''$. This also shows how the $\Xi'\ts_T\Lambda''$-module structure on $Y\ts_T U$ depends on the choice of $V$.

We obtain an exact functor
\[
\begin{array}{rcl}
\modr\Xi'            & \mra{F := (-\ts_T U)} & \modr(\Xi'\ts_T\Lambda'') \\
(Y\lraa{\alpha} Y')  & \mramaps              & (Y\ts_T U \;\lraa{\alpha\ts 1}\; Y'\ts_T U)\; . \\
\end{array}
\]
\end{Definition}

\begin{Lemma}
\label{LemND5}
The functor $F$ maps projective $\Xi'$-modules to projective $\Xi'\ts_T\Lambda''$-modules. Moreover, $F$ is full and faithful. 

\rm
{\it Proof.} To show that $F\Xi'$ is projective over $\Xi'\ts_T\Lambda''$, we remark that 
\[
\begin{array}{rcl}
F\Xi' \= \Xi'\ts_T U & \lraiso  & (1\ts e''_{0,0})(\Xi'\ts_T\Lambda'') \\
         \xi'\ts u^j & \lramaps & \xi'\ts e''_{0,j}\; , \\
\end{array}
\]
where $j\in [0,n-1]$, is an isomorphism of $\Xi'\ts_T\Lambda''$-modules.

We shall prove that $F$ is full and faithful. Given $Y_1,\,Y_2\,\in\,\modr\Xi'$, we claim that 
\[
\liu{\Xi'}(Y_1,Y_2)\;\lraa{F}\;\liu{\Xi'\ts_T\Lambda''}(FY_1,FY_2)
\]
is an isomorphism. Using a two-step free resolution of $Y_1$ and exactness of $F$, we may assume that $Y_1 = \Xi'$. Likewise, using projectivity of
$F\Xi'$ and a two-step free resolution of $Y_2$, we may assume $Y_2 = \Xi'$. But in this case, we have an isomorphism
\[
\begin{array}{rclclcl}
\Xi' & \lraiso  & \liu{\Xi'}(\Xi',\Xi') & \lraisoa{F} & \liu{\Xi'\ts_T \Lambda''}(F\Xi',F\Xi') & \lraiso  & 
                                                                                              (1\ts e''_{0,0})(\Xi'\ts_T\Lambda'')(1\ts e''_{0,0}) \\
\xi' & \lramaps & \xi'\cdot(-)          & \lra        & (\xi'\ts 1)\cdot (-)                   & \lramaps & \xi'\ts e''_{0,0}. \\
\end{array}
\]
\qedfl{-40}
\end{Lemma}

\begin{Lemma}
\label{LemND6}
There is an isomorphism
\[
\Ext^\ast_{\Xi'}(T,T)\;\mraisoa{\Ext^\ast(F)}\;\Ext^\ast_{\Xi'\ts_T\Lambda''}(U,U)
\]
of graded $S$-algebras with respect to the Yoneda product, where $U$ is a $\Xi'\ts_T\Lambda''$-module via $U \iso T\ts_T U = FT$.

\rm
{\it Proof.} Since $F$ preserves projectivity (\ref{LemND5}), we may apply $F$ to a projective resolution of $T$ to obtain a projective resolution of $U$, and 
use this resolution to calculate $\Ext^\ast_{U\wr H}(U,U)$. Application of $F$ to morphisms of complexes modulo homotopy now yields the morphism of 
graded $S$-algebras $\Ext^\ast(F)$. But since $F$ is full and faithful (\ref{LemND5}), this map is an isomorphism.\qed
\end{Lemma}

\begin{Proposition}
\label{PropND7}
If $n\not\con_p 0$, then 
\[
\Ext^\ast_{T\wr G}(T,T)\;\iso\;\Ext^\ast_{U\wr H}(U,U)
\]
as graded $S$-algebras.

\rm
{\it Proof.} By (\ref{PropND3}), we have $\Xi'\ts_T\Lambda''\iso\Xi\iso U\wr H$. So (\ref{LemND6}) gives the assertion, provided the module structures of 
$\Xi'\ts_T\Lambda''$ and of $U\wr H$ on $U$ coincide. Suppose given $\xi'\ts e''_{k,l}\in \Xi'\ts_T\Lambda''$, $\xi'\in\Xi'$, $k,\, l\,\in\, [0,n-1]$.
Let $(x_1,\dots,x_g)$ be an $S$-linear basis of $T$. 

In the interpretation $U\iso T\ts_T U = FT$, we have 
$(x_i\ts u^j)(\xi'\ts e''_{k,l}) = (x_i\xi')\ts (\dell_{j,k} u^l) = 1\ts (x_i\xi')(\dell_{j,k} u^l)$ by definition of $F$.

On the other hand, by (\ref{LemIns}), the element $\xi'\ts e''_{k,l}$ is mapped via $\tht$ to the element in $\Gamma$ that maps 
$x_i u^j$ to $(x_i\xi')(\dell_{j,k} u^l)$. So after identifying $T\ts_T U\iso U$, the operations coincide.\qed
\end{Proposition}

\bq
\begin{Example}[cont.]
\label{ExHelp4}\rm
We have
\[
\Ext^\ast_{\Z_{(3)}[\zeta_{3^2}]\wr (C_3\ti C_2)}(\Z_{(3)}[\zeta_{3^2}],\Z_{(3)}[\zeta_{3^2}]) 
\;\iso\; \Ext^\ast_{\Z_{(3)}[\pi_2]\wr C_3}(\Z_{(3)}[\pi_2],\Z_{(3)}[\pi_2])
\]
as graded $\Z_{(3)}$-algebras.
\end{Example}
\eq
\section{Cup product and Yoneda product coincide}
\label{SecCupYoneda}

In this section, we let $T$ be a commutative ring, we let $G$ be a group acting on $T$ via a group morphism $G\lra\Aut_\mb{\scr ring} T$, and define
$S := \Fix_G T$ to be the fixed ring of this operation. We assume that $T$ is a free module over $S$.

We consider the twisted group ring $T\wr G$ with respect to this operation, which is an $S$-algebra.

We write $\ts := \ts_S$. Let 
\[
P := (\cdots\lraa{d_3} SG^{\ts 3} \lraa{d_2} SG^{\ts 2}\lraa{d_1} \ub{SG^{\ts 1}}_{\mb{\scr degree 0}} \lra 0 \lra \cdots)
\]
denote the bar resolution, which is a projective resolution of the trivial module $S$ over $SG$. The differential is given by
\[
\begin{array}{rcl}
SG^{\ts i+1} & \lraa{d_i} & SG^{\ts i} \\
g_{[0,i]}    & \lramaps   & \sumd{l\in [0,i]} (-1)^{i-l} g_{[0,i]\ohne\{ l\}}\\
\ea
\]
where $i\geq 1$, and where we write shorthand
\[
g_{A} \; := \; g_{a_1} \ts g_{a_2} \ts \cdots \ts g_{a_k},
\]
for $A\tm\Z$, where $k := \# A < \infty$, where $a_i\in A$ and $g_{a_i}\in G$ for $i\in [1,k]$, and where $a_1 < a_2 < \dots < a_k$. Let 
$SG\lraa{d_0} S$ denote the augmentation map, sending each element of $G$ to $1\in S$. We have a quasiisomorphism $P\lraa{d_0} S$, where $S$ is 
regarded as a complex concentrated in degree $0$ and where $d_0$ also denotes the morphism of complexes given by $d_0$ in degree $0$.

Let
\[
\begin{array}{rcl}
\Modr SG            & \lraa{L} & \Modr(T\wr G) \\
M                   & \lramaps & LM \; :=\; M\ts_{SG} T\wr G \\
\Modr SG            & \llaa{R} & \Modr(T\wr G) \\
N|_{SG}  \; =:\; RN & \llamaps & N \; . \\
\ea
\]
Then $L$ is left adjoint to $R$, with adjunction morphisms denoted by
\[
\begin{array}{rcll}
LRN    & \lraa{\eta N} & N  & \Icm\mb{counit} \\
n\ts x & \lramaps      & nx & \\
RLM    & \llaa{\eps M} & M  & \Icm\mb{unit}\\
m\ts 1 & \llamaps      & m  & \\
\ea
\] 
We note that $T\iso LS$, by means of $x\lramaps 1\ts x$, $x^g\llamaps 1\ts xg$, $x\in T$, $g\in G$, which we use as identification. 
Since $T\wr G$ is free as a module over $SG$, both $R$ and $L$ are exact. To complexes, $R$ and $L$ are applied entrywise. 
So for instance, $RLS$ is the complex of $SG$-modules having entry $T$ in degree $0$ and entry $0$ elsewhere.

We use the description
\[
\Ext^i_{T\wr G}(T,T) \= \liu{\KK^-(T\wr G)}{(LP,LP[i])}\; , 
\]
where $[i]$ denotes the shift of complexes by $i$ steps to the left, i.e.\ $(X[i])_j := X_{j-i}$, where $i\geq 0$, $j\in\Z$. Suppose given 
$u \in \liu{\KK^-(T\wr G)}{(LP,LP[i])}$,
$v \in \liu{\KK^-(T\wr G)}{(LP,LP[j])}$, 
for some degrees $i,j\geq 0$. The {\it Yoneda product} $u\cdot v$ is given as the composite
\[
(LP\lraa{u\cdot v} LP[i+j]) \; := \; (LP \lraa{u} LP[i] \lraa{v[i]} LP[i+j])\;\in\; \liu{\KK^-(T\wr G)}{(LP,LP[i+j])}\; . 
\]
It turns $\Ext^\ast_{T\wr G}(T,T)$ into a graded $S$-algebra.

Moreover, we use the description of the cohomology of $G$ with coefficients in $T$ over the ground ring $S$
\[
\HH^i(G,T;S) \; :=\; \Ext^i_{SG}(S,T) \= \liu{\KK^-(SG)}{(P,RLS[i])} \= \HH^i(CP)\; , 
\]
where $i\geq 0$, and where $C:\Modr SG\lra \Modr S:X\lramaps\liu{SG}{(X,T)}$ is applied entrywise to complexes. For a morphism $f$,
we denote $f^\ast := Cf$. Let $\Zyk^i(CP) := \Kern (CP_i\lraa{d_{i+1}^\ast} CP_{i+1})$ denote the $S$-module of $i$-cocycles.

We will also make use of the alternative interpretation
\[
\barcl
\HH^i(G,T;S) \= \liu{\KK^-(SG)}{(P,RLS[i])} & \llaisoa{\lm} & \liu{\KK^-(SG)}{(P,RLP[i])} \\
                          u (RL d_0[i]) & \llamaps      & u\; . \\
\ea
\]
Given $a\in CP_i = \liu{SG}{(SG^{\ts i+1},T)}$, $b\in CP_j = \liu{SG}{(SG^{\ts j+1},T)}$, for some degrees $i,j\geq 0$, their {\it cup product} 
$a\cup b\in  CP_{i+j} = \liu{SG}{(SG^{\ts i+j+1},T)}$ is defined by
\[
(g_{[0,i+j]})(a\cup b) \; := \; (g_{[0,i]})a \cdot (g_{[i,i+j]})b\; ,
\]
where $g_l\in G$ for $l\in [0,i+j]$. Because of the Leibniz rule 
\[
(a\cup b)d_{i+j+1}^\ast \= (-1)^j (a d_{i+1}^\ast\cup b) + (a\cup b d_{j+1}^\ast)\; ,
\]
the cup product restricts to $\Zyk^i(CP)\ti \Zyk^j(CP)\lraa{\cup} \Zyk^{i+j}(CP)$ and induces a map 
\[
\HH^i(G,T;S)\ti \HH^j(G,T;S) \;\lraa{\cup}\; \HH^{i+j}(G,T;S)\; .
\]
Given cocycles $a\in\Zyk^i(CP)$, $b\in\Zyk^j(CP)$, where $i,j\geq 0$, we let $c_{a,b}\in CP_{i+j-1}$ be defined by
\[
h_{[0,i+j-1]} c_{a,b} \; :=\; \sumd{m\in [0,j-1]} (-1)^{m(i+j-1)}(h_{[m,m+i]})a\cdot (h_{[m+i,i+j-1]}\ts h_{[0,m]})b\; ,
\]
where $h_l\in G$ for $l\in [0,i + j - 1]$. We obtain
\[
(g_{[0,i+j]}) (c_{a,b} d_{i+j}^\ast) \= g_{[0,i]} a\cdot g_{[i,i+j]} b - (-1)^{ij} g_{[0,j]} b\cdot g_{[j,i+j]} a\; ,
\]
whence $a\cup b = (-1)^{ij}\, b\cup a$ as elements of $\HH^{i+j}(G,T;S)$. Thus the cup product turns $\HH^\ast(G,T;S)$ into a graded commutative $S$-algebra.

\bq
The following proposition we owe to {\sc B.\ Keller.}
\eq

\begin{Proposition}
\label{PropCup4}
The isomorphism $\lm^{-1}\ka$ of graded $S$-modules, given by
\[
\barcl
\HH^\ast(G,T;S) \;\llaisoa{\lm}\; \liu{\KK^-(SG)}{(P,RLP[\ast])} & \lraisoa{\ka}  & \liu{\KK^-(T\wr G)}{(LP,LP[\ast])} \= \Ext^\ast_{T\wr G} (T,T) \\
                                                               u & \lramaps          & (Lu)(\eta LP[\ast]) \\
                                                   (\eps P)(Ru') & \llamaps          & u' \; , \\
\ea
\]
is in fact an isomorphism of graded $S$-algebras, with respect to the cup product on \linebreak 
$\HH^\ast(G,T;S)$ and with respect to the Yoneda product on $\Ext^\ast_{T\wr G} (T,T)$. 

In particular, $\Ext^\ast_{T\wr G} (T,T)$ is a graded commutative $S$-algebra. 

Note that $\HH^0(G,T;S) \iso \Hom_{T\wr G} (T,T) \iso S$ as $S$-algebras. 

\rm
{\it Proof.} Given $u\in\liu{\KK^-(SG)}{(P,RLP[i])}$, $v\in\liu{\KK^-(SG)}{(P,RLP[j])}$, where $i,j\geq 0$, we calculate
\[
\barcl
((u\ka)\cdot (v\ka))\ka^{-1}
& = & \Big((Lu)(\eta LP[i])(Lv[i])(\eta LP[i+j])\Big)\ka^{-1} \vspace*{1mm}\\
& = & (\eps P)(RLu)(R\eta LP[i])(RLv[i])(R\eta LP[i+j]) \vspace*{1mm}\\
& = & u (RLv[i])(R\eta LP[i+j])\; . \\
\ea
\]
Now suppose given $a\in\Zyk^i(CP)$, $b\in\Zyk^j(CP)$. Letting $P_{i+l}\lraa{\w a_{i+l}} RLP_l$ be defined by
\[
(g_{[0,i+l]})\w a_{i+l} \; :=\;  (-1)^{il} g_{[0,l]}\ts (g_{[l,i+l]}a)\; ,
\]
where $l\geq 0$ and $g_k\in G$ for $k\in [0,i+l]$, we obtain a morphism of complexes $P\lraa{\w a} RLP[i]$ that is mapped by $\lm$ to 
$P\lraa{a} RLS[i]$. Hence
\[
\ba{rl}
  & (g_{[0,i+j]})\Big(((a\lm^{-1}\ka)\cdot (b\lm^{-1}\ka))\ka^{-1}\lm\Big) \vspace*{1mm}\\
= & (g_{[0,i+j]})\Big(((\w a\ka)\cdot (\w b\ka))\ka^{-1}\lm\Big) \vspace*{1mm}\\
= & (g_{[0,i+j]}) \w a_{i+j} (RL\w b_j) (R\eta L(SG)) (RL d_0[i+j]) \vspace*{1mm}\\
= & (-1)^{ij}(g_{[0,j]}\ts (g_{[j,i+j]} a))(RL\w b_j) (R\eta L(SG)) (RL d_0[i+j]) \vspace*{1mm}\\
= & (-1)^{ij}(g_0\ts (g_{[0,j]} b)\ts (g_{[j,i+j]} a)) (R\eta L(SG)) (RL d_0[i+j]) \vspace*{1mm}\\
= & (-1)^{ij}(g_0\ts (g_{[0,j]} b)\cdot (g_{[j,i+j]} a)) (RL d_0[i+j]) \vspace*{1mm}\\
= & (-1)^{ij} (g_{[0,j]} b)\cdot (g_{[j,i+j]} a) \vspace*{1mm}\\
= & (-1)^{ij} (g_{[0,i+j]}) (b\cup a)\; , \\
\ea
\]
i.e.\
\[
(a\lm^{-1}\ka)\cdot (b\lm^{-1}\ka) \= (a\cup b)\lm^{-1}\ka\; .
\]
\qedfl{-60}
\end{Proposition}

\bq
 \begin{Remark}
 \label{RemAIP}\rm
 The graded commutativity of $\Ext^\ast_{T\wr G} (T,T)$ can also be obtained by {\bf\cite[\rm 2.1]{Su04}}.
 \end{Remark}
\eq
\section{Cohomology}
\label{SecAddGal}

\setcounter{subsection}{-1}
\subsection{A classical approach}

\bq
 \begin{Example}
 \label{ExLEO}\rm
 We shall calculate the cohomology directly in an example, still disregarding the cup product, however. 

 Let $\pi_2 = \prod_{j\in [1,p-1]} (\zeta_{p^2}^{j^p} - 1)$, cf.\ section \ref{SubsecCNF} below. Let
 \[
 \ba{rclrcl}
 S & = & \Z_{(p)} \; ,       & K & = & \Q\; , \\
 T & = & \Z_{(p)}[\pi_2]\; , & L & = & \Q(\pi_2)\; . \\
 \ea
 \]
 Then $G = C_p = \spi{\sa}$, where $\sa$ is the restriction to $\Q(\pi_2)$ to the automorphism $\zeta_{p^2} \lramaps \zeta_{p^2}^{1+p}$ of $\Q(\zeta_{p^2})$. 
 Let $e := \frac{1}{p}\sum_{j\in [0,p-1]} \sa^j \in \Q G$. We have a split short exact sequence of $SG$-modules
 \[
 \begin{array}{rcl}
 0 \; \lra \; M\;\lra\; T & \lra     & S\;\lra\; 0 \\
                        x & \lramaps & xe \\ 
                        y & \llamaps & y\; , \\
 \end{array}
 \]
 which is welldefined since $\Tr_{L|K}(x)$ is divisible by $p$ for all $x\in T$ because $T|S$ is wildly ramified.

 Now $M$ is an $SG(1-e)$-lattice of rank $\rk_S M = p-1$. Since $SG(1-e)$ is isomorphic to $S[\zeta_p]$ via $\sigma(1-e)\lramapsiso\zeta_p$ as an $S$-algebra, we 
 have $M\iso SG(1-e)$. Using the $2$-periodic projective resolution 
 \[
 \cdots\;\lraa{b}\; SG\;\lraa{a}\; SG\;\lraa{b}\; SG\;\lra\; 0
 \]
 of $S$, where $a:1\lramaps pe$ and $b: 1\lramaps \sigma - 1$, we obtain
 \[
 \ba{rcl}
 \HH^j(G,T;S)
 &\iso & \HH^j(G,SGe\ds SG(1-e);S) \vspace*{2mm}\\
 &\iso & \HH^j(G,SGe;S)\ds \HH^j(G,SG(1-e);S) \vspace*{2mm}\\
 &\iso & \left\{
         \ba{ll}
         0        &\mb{$j$ odd} \\
         \Fu{p}   &\mb{$j\geq 2$ even} \\
         S        &\mb{$j = 0$} \\
         \ea
         \right\}\ds 
         \left\{
         \ba{ll}
         \Fu{p}   &\mb{$j$ odd} \\
         0        &\mb{$j$ even} \\
         \ea
         \right\}\vspace*{2mm}\\
 &\iso & \left\{
         \ba{ll}
         \Fu{p}     &\mb{$j\geq 1$} \\
         S          &\mb{$j = 0$} \\
         \ea
         \right\} 
 \ea
 \]
 as an $S$-module. This will be confirmed by (\ref{PropAG8}, \ref{ExCyc2}) below, since $\HH^j(G,T;S) \iso \Ext^j_{T\wr G}(T,T)$ by adjunction.
 \end{Example}

 \begin{Remark}
 \label{RemAG4}\rm
 We attempt to explain why in the case of a cyclic galois group, and in presence of a classical periodic resolution, the description in
 (\ref{ThFT16}) is actually useful  to calculate products in cohomology.

 For the purpose of this remark, let $G = \spi{\sigma\,|\, \sigma^g = 1}$ be a cyclic group of order $g\geq 1$ acting on the discrete valuation 
 ring $T$, let $S = \mb{\rm Fix}_G T$ be the fixed ring in $T$ under $G$, and let $L = \fracfield T$, $K = \fracfield S$. Let 
 \[
 \ba{rclcrcl}
 T\wr G & \lraa{\al_0}  & T\wr G & : & 1 & \lramaps & \left(\sum_{i\in \Z/g}\sigma^i\right) \\
 T\wr G & \lraa{\be_0}  & T\wr G & : & 1 & \lramaps & (\sa - 1) \\
 T\wr G & \lraa{\eps_0} & T      & : & 1 & \lramaps & 1 \\
 \ea
 \]
 be $T\wr G$-linear maps. We obtain a periodic projective resolution of period $2$
 \[
 P_0 := \left(\cdots\;\lra\; T\wr G \;\lraa{\al_0}\; T\wr G \;\lraa{\be_0}\; T\wr G \;\lraa{\al_0}\; T\wr G 
              \;\lraa{\be_0}\hspace*{-1mm}\ub{T\wr G}_\mb{\scr degree $0$}\hspace*{-1mm}\lra\; 0\;\lra\;\cdots \right)\; ,
 \]
 mapping quasiisomorphically to $T$ via $\eps_0$ in degree $0$ since the image of $\al_0$ is isomorphic to $T$ via $T\lraiso\Img\al_0$, $x\lramaps\sum_{i\in \Z/g}\sigma^i x$. So we have
 \[
 \ba{l}
 \RHom_{T\wr G}(T,T) \= \\
 \left(\cdots\;\lla\; T \;\llaa{\Tr_{L|K}}\; T \;\mla{x^\sa - x\lamaps x}\; T \;\llaa{\Tr_{L|K}}\; T \;
                              \mla{x^\sa - x\lamaps x}\hspace*{-2mm}\ub{T}_\mb{\scr degree $0$}\hspace*{-2mm}\lla\; 0\;\lla\;\cdots \right)\; , \\
 \ea
 \]
 whence for instance $\Ext_{T\wr G}^1(T,T) \iso \{x\in T : \Tr_{L|K}(x) = 0\}/\{ y^\sa - y : y\in T\}$, or 
 $\Ext_{T\wr G}^2(T,T) \iso S/\Tr_{L/K}(T)$. Suppose given an element of, say, $\Ext_{T\wr G}^1(T,T)$, represented by an element $e\in T$ with trace 
 $0$. To apply Yoneda multiplication, we need to represent it as an element of $\liu{\KK^-(T\wr G)}(P_0,P_0[1])$. Let us consider the necessary 
 construction of a morphism of complexes. 
 We attempt to construct a periodic resolution of $e$ with period $2$. If in odd degrees, the morphism of complexes is given by 
 $1\lramaps \sum_{i\in\Z/g} \sa^i x_i$, and in even degrees $\geq 2$ by $1\lramaps \sum_{i\in\Z/g} \sa^i y_i$, where $x_i,\, y_i\,\in\, T$, 
 then these coefficients are subject to the following conditions.
 \[
 \barcl
 \sum_i x_i                       & = & e \vspace*{1mm}\\
 (\sum_i \sa^i x_i)(\sum_j \sa^j) & = & (\sa - 1)(\sum_i \sa^i y_i) \vspace*{1mm}\\
 (\sum_j \sa^j)(\sum_i \sa^i x_i) & = & (\sum_i \sa^i y_i)(\sa - 1)\; , \\
 \ea
 \]
 i.e.\ 
 \[
 \ba{rcll}
 \sum_i x_i                     & = & e \\
 \sum_{i+j\con_g k} x_i^{\sa^j} & = & y_{k-1} - y_k     & \mb{for all $k\in\Z/g$}\\
 e                              & = & y_{k-1}^\sa - y_k & \mb{for all $k\in\Z/g$}\; . \\
 \ea
 \]
 As far as we can see, solving this system requires knowledge of the operation of $\sa$ on $T$ in terms of an $S$-linear basis of $T$. 
 Cf.\ e.g.\ (\ref{ExHelp1}) or appendix \ref{Subsecp3}.
 \end{Remark}
\eq

\subsection{A projective resolution}
\label{SubsecProjres}

We maintain the notation and assumptions of section \ref{SubsecComp}. In particular, $G = C_p = \spi{\sa}$ is cyclic of order equal to the residue 
characteristic $p$ of $S$, and $b = -1 + \val_t(t^\sa - t)\geq 1$. To dispose of the equality $\Xi = \LD$, we assume that 
$\val_s(p)\geq b - \ul{b}$ (\ref{ThFT16}).

Let
\[
\ba{rclcrclrl}
\Xi & \lraa{\al}    & \Xi & : & 1 & \lramaps & \w\mu_{p-1,p-\ol{b}} & := & s^{b - \ul{b} - 1}\,\eps_{\ol{b},p-\ol{b}}              \\
\Xi & \lraa{\be}    & \Xi & : & 1 & \lramaps & \mu_{1,\ol{b}}       &  = & s^{\ul{b}}\sum_{i\in [0,p-1]} i\eps_{i,\ol{b}}         \\
\Xi & \lraa{\chi_k} & T   & : & 1 & \lramaps & t^k                  &    & \mb{for $k\in [0,p-1]$} \\ 
\ea
\]
be $\Xi$-linear maps (cf.\ \ref{CorFT14}). Here $\w\mu_{p-1,p-\ol{b}}$ is a `shifted version' of $\mu_{p-1,p-\ol{b}}$ (cf.\ \ref{RemFT15}).

\begin{Proposition}
\label{LemAG6}
The complex of $\Xi$-linear maps
\[
\cdots\;\lraa{\be}\;\Xi\;\lraa{\al}\;\Xi\;\lraa{\be}\;\Xi\;\lraa{\al}\;\Xi\;\lraa{\be}\Xi\;\lraa{\al}\;\cdots\; , 
\]
periodic of period $2$, is acyclic. The image of $\al$ is isomorphic to $T$ as a module over $\Xi$; more precisely, we have a factorization
\[
(\Xi\lraa{\al}\Xi) \= (\Xi\lraa{\chi_0} T\lra\Xi),
\]
with $\chi_0$ surjective and $T\lra\Xi$ injective. 

\rm
Denote by
\[
P \; :=\; (\cdots\;\lra\;\Xi\;\lraa{\al}\;\Xi\;\lraa{\be}\hspace*{-3mm}\ub{\Xi}_\mb{\scr\rm degree $0$}\hspace*{-3mm}\lra\; 0\;\lra\;\cdots) 
\;\in\; \KK^-(\Xi)
\]
the resulting projective resolution of $T$.

{\it Proof.} We claim exactness of $\Xi \lraa{\al} \Xi \lraa{\be} \Xi$. To prove that $\al\be = 0$, we calculate
\[
\barcl
\mu_{1,\ol{b}}\,\w\mu_{p-1,p-\ol{b}}
& = & s^{b-1}\sum_{i\in [0,p-1]} i\eps_{i,\ol{b}}\,\eps_{\ol{b},p-\ol{b}} \\
& \aufgl{\mb{\scr (\ref{LemFT11})}} &  s^{b-1}\sum_{i\in [0,p-1]} i\dell_{\,\ol{b},\ol{i+b}}\,\eps_{i,p} \\
& = & 0\; . \\
\ea
\]
Now, let $B\tm\Lambda$ be the kernel of the $\Lambda$-linear map $\Lambda\lra\Lambda$, $1\lramaps\mu_{1,\ol{b}}$.
To prove that $\Img\al = \Kern\be$, we shall show that the $S$-linear colengths of both submodules in $B$ coincide.  For $j,k\in [0,p-1]$, we have
\[
\barcl
\mu_{1,\ol{b}}\,\eps_{j,k} 
& = & s^{\ul{b}} \sum_{i\in [0,p-1]} i\eps_{i,\ol{b}}\,\eps_{j,k} \\
& \aufgl{\mb{\scr (\ref{LemFT11})}} & s^{\ul{b}}\sum_{i\in [0,p-1]} i \dell_{j,\ol{i+b}}\,\eps_{i,k+\ol{b}}\vspace*{2mm} \\
& = & s^{\ul{b}}\;\ol{j-b}\;\eps_{\ol{j-b},k+\ol{b}}\; , 
\ea
\]
and thus $B$ has the $S$-linear basis $(\eps_{\ol{b},k})_{k\in [0,p-1]}$. By (\ref{CorFT14}), an element $\sum_{k\in [0,p-1]} x_k\eps_{\ol{b},k}$, 
where $x_k\in S$, is in $\Xi$ if and only if 
\[
\val_s(x_k)\;\geq\; 1 + \ul{b(p-1)-k-1}
\]
for $k\in [0,p-1]$. By (\ref{LemFT12_5}), we obtain the colength of $\Kern\beta$ in $B$ to be equal to 
$\sum_{k\in [0,p-1]} (1+\ul{b(p-1)-k-1}) = b(p-1)$. On the other hand we obtain
\[
\barcl
\mu_{l,j}\al 
& = & \w\mu_{p-1,p-\ol{b}}\,\mu_{l,j} \vspace*{1mm}\\
& = & s^{b - \ul{b} + \ul{bl - j - 1}} \sum_{i\in [0,p-1]} \smatze{i}{l} \eps_{\ol{b},p-\ol{b}}\,\eps_{i,j} \\
& \aufgl{\mb{\scr (\ref{LemFT11})}} & s^{b - \ul{b} + \ul{bl - j - 1}} \sum_{i\in [0,p-1]} \smatze{i}{l} \dell_{i,0}\,\eps_{\ol{b},j+p-\ol{b}} \\
& = & \dell_{l,0} s^{b - \ul{b} - 1} \eps_{\ol{b},j + p - \ol{b}} \\
& = & 
\left\{
\ba{ll}
s^{b - \ul{b} - 1} \eps_{\ol{b},j + p - \ol{b}}  & \mb{for $l = 0$ and $j\in [0,\ol{b}-1]$} \\
s^{b - \ul{b}}\,\eps_{\ol{b},j - \ol{b}}         & \mb{for $l = 0$ and $j\in [\ol{b},p-1]$} \\
0 & \mb{for $l\in [1,p-1]\; $,} \\
\ea
\right. \\
\ea 
\]
where $l,j\in [0,p-1]$, whence the colength of $\Img\al$ in $B$ equals $(b - \ul{b} - 1)p + (p - \ol{b}) = b(p-1)$, too. Thus $\Img\al = \Kern\be$.

Moreover, since
\[
\barcl
\mu_{l,j}\chi_0 
& = & s^{1 + \ul{bl - j - 1}}\sum_{i\in [0,p-1]}\smatze{i}{l}t^0\eps_{i,j} \\
& = & \dell_{l,0} t^j\; , \\
\ea
\]
where $l,j\in [0,p-1]$, the $\Xi$-linear isomorphism
\[
\barcl
T   & \lraiso  & \Xi\al \\
t^j & \lramaps & s^{b - \ul{b} - 1}\eps_{\ul{b},j+p-\ul{b}} \\
\ea
\]
yields the commutativity
\[
(\;\Xi\;\lraa{\al}\;\Xi\al\;) \= (\;\Xi\;\lraa{\chi_0}\; T\;\lraiso\;\Xi\al\;) \; .
\]

We claim exactness of $\Xi \lraa{\be} \Xi \lraa{\al} \Xi$. To prove that $\be\al = 0$, we calculate 
\[
\barcl
\w\mu_{p-1,p-\ol{b}}\,\mu_{1,\ol{b}}
& = & s^{b-1}\sum_{i\in [0,p-1]} i\eps_{\ol{b},p-\ol{b}}\,\eps_{i,\ol{b}} \\
& \aufgl{\mb{\scr (\ref{LemFT11})}} &  s^{b-1}\sum_{i\in [0,p-1]} i \dell_{i,0}\,\eps_{\ol{b},p} \\
& = & 0\; . \\
\ea
\] 
Now, let $A\tm\Lambda$ be the kernel of the $\Lambda$-linear map $\Lambda\lra\Lambda$, $1\lramaps\w\mu_{p-1,p-\ol{b}}$.
To prove that $\Img\be = \Kern\al$, we shall show that the $S$-linear colengths of both submodules in $A$ coincide.  For $j,k\in [0,p-1]$, we have
\[
\eps_{j,k}\al \= s^{b - \ul{b} - 1}\eps_{\ol{b},p-\ol{b}}\,\eps_{j,k} 
\;\aufgl{\mb{\scr (\ref{LemFT11})}}\; s^{b - \ul{b} - 1}\dell_{j,0}\,\eps_{\ol{b},p-\ol{b}+k}\; ,
\]
and therefore $A$ has the $S$-linear basis $(\eps_{j,k})_{j\in [1,p-1],\; k\in [0,p-1]}$. By (\ref{CorFT14}), $\Kern\al = \Xi\cap A$ has the 
$S$-linear basis 
$(\mu_{l,j})_{l\in [1,p-1],\; j\in [0,p-1]}$, whence the colength of $\Kern\al$ in $A$ equals $bp(p-1)/2$ (cf.\ pf.\ of \ref{CorFT14}). 
On the other hand, using the $S$-linear basis 
\[
\left(s^{1 + \ul{bl - j - 1}} \sum_{k\in [0,p-1]} \smatze{k}{l} \eps_{\ol{k+b+1},j}\right)_{\!\! l,j\in [0,p-1]}
\]
of $\Xi$ (\ref{RemFT15}), we obtain
\[
\ba{cl}
  & \left(s^{1 + \ul{bl - j - 1}} \sum_{k\in [0,p-1]} \smatze{k}{l} \eps_{\ol{k+b+1},j}\right)\be \\
= & s^{\ul{b} + 1 + \ul{bl-j-1}}\sum_{i\in [0,p-1]}\sum_{k\in [0,p-1]} i\smatze{k}{l}\eps_{i,\ol{b}}\,\eps_{\ol{k+b+1},j} \\
\aufgl{\mb{\scr (\ref{LemFT11})}} & 
   s^{\ul{b} + 1 + \ul{bl-j-1}}\sum_{i\in [0,p-1]}\sum_{k\in [0,p-1]} i\smatze{k}{l}\dell_{\,\ol{k+b+1},\ol{i+b}}\,\eps_{i,\ol{b} + j}\vspace*{1mm}\\
= & s^{\ul{b} + 1 + \ul{bl-j-1}}\sum_{i\in [l+1,p-1]} i\smatze{i-1}{l}\eps_{i,\ol{b} + j} \\
= & s^{\ul{b} + 1 + \ul{bl-j-1} + \ul{\ol{b} + j}}\sum_{i\in [l+1,p-1]} i\smatze{i-1}{l}\eps_{i,\ol{b + j}}\; . \\
\ea
\]
This yields the colength of $\Img\be$ in $A$ to be equal to 
\[
\ba{cl}
  & \sum_{j\in [0,p-1]}\sum_{l\in [0,p-2]} (\ul{\ol{b} + j} + \ul{b} + 1 + \ul{bl-j-1}) \\
\auf{\mb{\scr (\ref{LemFT12_5})}}{=} & \ol{b}(p-1) + (\ul{b} + 1)p(p-1) + \sum_{l\in [0,p-2]} (bl - p) \\
= & bp(p-1)/2\; ,\\
\ea
\]
too. Thus $\Img\be = \Kern\al$.\qed
\end{Proposition}

\bq
\begin{Example}[cont.]
\label{ExHelp5}\rm
For $\Xi \iso \Z_{(3)}[\pi_2]\wr C_3$, we have $b = 1$ and $\ul{b} = 0$, thus obtaining the projective resolution
\[
P \= \left(\cdots\;\lra\;\Xi\;\mramh{\rsmateckdd{0}{0}{0}{0}{0}{1}{3\cdot 2}{0}{0}\cdot (-)}
                       \;\Xi\;\mramh{\rsmateckdd{0}{0}{0}{3\cdot 1}{0}{0}{0}{0}{0}\cdot (-)}
                       \;\Xi\;\mramh{\rsmateckdd{0}{0}{0}{0}{0}{1}{3\cdot 2}{0}{0}\cdot (-)}
       \hspace*{-2mm}\ub{\Xi}_\mb{\scr degree $0$}\hspace*{-2mm}\lra\; 0\;\lra\;\cdots\right)
\]
of $T$, with quasiisomorphism given by $\Xi\lraa{\chi_0} T$ in degree $0$, sending 
$\xi = \rsmateckdd{a_{0,0}}{a_{0,1}}{a_{0,2}}{3 a_{1,2}}{a_{1,0}}{a_{1,1}}{3 a_{2,1}}{3 a_{2,2}}{a_{2,0}}\in\Xi$ to 
$t^0\xi = a_{0,0} t^0 + a_{0,1}t^1 + a_{0,2}t^2$. Representing elements of $T$ as row vectors with entries in $S$ with respect to the basis 
$(t^0,t^1,t^2)$, the map $\chi_0$ is given by $\smatecked{1}{0}{0}\cdot(-)$.
\end{Example}
\eq

\subsection{The Yoneda ring}
\label{SubsecYonring}

\begin{Remark}
\label{LemAG7}
An $S$-linear basis of $\liu{\Xi}{(\Xi,T)}$ is given by $(\chi_i)_{i\in [0,p-1]}$. The restriction map 
\[
\liu{\Lambda}{(\Lambda,T)} \;\lra\; \liu{\Xi}{(\Xi,T)}
\]
is an isomorphism.
\end{Remark}

\begin{Proposition}
\label{PropAG8}
We have
\[
\ba{rcll}
\Ext^0_{T\wr G}(T,T)      & =    & S\spi{\chi_0} & \\
                          & \iso & S \; ,        & \vspace*{1mm}\\
\Ext^{2i}_{T\wr G}(T,T)   & =    & S\spi{\chi_0}/S\spi{s^{b - \ul{b}}\cdot\chi_0} & \vspace*{1mm}\\
                          & \iso & S/s^{b - \ul{b}} & \mb{for $\; i\geq 1\; $, and}  \vspace*{1mm}\\
\Ext^{2i+1}_{T\wr G}(T,T) & = & \fracd{S\spi{\;\;\;\;\;\;\chi_0,\dots,\;\;\;\;\;\;\chi_{\ol{b}-1},\;\;\;\chi_{\ol{b}+1},\dots,\;\;\;\chi_{p-1}}}
             {S\spi{s^{\ul{b} + 1}\chi_0,\dots,s^{\ul{b} + 1}\chi_{\ol{b}-1},s^{\ul{b}}\chi_{\ol{b}+1},\dots,s^{\ul{b}}\chi_{p-1}}} & \vspace{2mm}\\
    & \iso & \left({\dis\Ds_{k\in [0,\ol{b}-1]}} S/s^{\ul{b}+1}\right) \ds \left({\dis\Ds_{k\in [\ol{b}+1,p-1]}} S/s^{\ul{b}}\right) & 
            \mb{for $\; i\geq 0\;$.} \\
\ea
\]
The element represented by $\chi_0$ in $\Ext^{2i}_{T\wr G}(T,T)$ shall be written $\chi_0^{(2i)}$, the element represented by $\chi_j$ in 
$\Ext^{2i+1}_{T\wr G}(T,T)$ shall be written $\chi_j^{(2i+1)}$, where $i\geq 0$, $j\in [0,p-1]\ohne\{ 1\}$.

\rm
{\it Proof.} We calculate $\liu{\Xi}{(\Xi,T)}\lraa{\al^\ast}\liu{\Xi}{(\Xi,T)}$, $\chi\lramaps \al\chi$. Given $k\in [0,p-1]$, we get
\[
\barcl
1_\Xi(\chi_k \al^\ast)
& = & (1_\Xi\al)\chi_k \vspace*{1mm}\\
& = & (s^{b - \ul{b} - 1}\eps_{\ol{b},p-\ol{b}}) \chi_k \vspace*{1mm}\\
& = & s^{b - \ul{b} - 1} (1_\Xi\,\chi_k) \eps_{\ol{b},p-\ol{b}}\vspace*{1mm}\\
& = & \dell_{\,\ol{b},k} s^{b - \ul{b}}  \vspace*{1mm}\\
& = & \dell_{\,\ol{b},k} s^{b - \ul{b}}  (1_\Xi\,\chi_0)\; , \\
\ea
\]
i.e.\
\[
\chi_k \al^\ast \= \dell_{\,\ol{b},k} s^{b - \ul{b}} \cdot \chi_0\; .
\]
We calculate $\liu{\Xi}{(\Xi,T)}\lraa{\be^\ast}\liu{\Xi}{(\Xi,T)}$, $\chi\lramaps \be\chi$. Given $k\in [0,p-1]$, we get
\[
\barcl
1_\Xi(\chi_k \be^\ast)
& = & (1_\Xi\be)\chi_k \vspace*{1mm}\\
& = & (s^{\ul{b}}\sum_{i\in [0,p-1]} i\eps_{i,\ol{b}}) \chi_k \vspace*{1mm}\\
& = & s^{\ul{b}}\sum_{i\in [0,p-1]} i \dell_{i,k} t^{\ol{k+b}} s^{\ul{k+\ol{b}}} \vspace*{1mm}\\
& = & k s^{\ul{k+\ol{b}} + \ul{b}} (1_\Xi\,\chi_{\ol{k+b}})\; , \\
\ea
\]
i.e.\
\[
\chi_k \be^\ast \= k s^{\ul{k+\ol{b}} + \ul{b}} \cdot \chi_{\ol{k+b}}\; .
\]
The shape of the $\Ext$-groups now follows by (\ref{LemAG6}).\qed
\end{Proposition}

\begin{Lemma}
\label{LemAG9}
We have
\[
s^{b + \ul{j-2b}}\left((\ol{b-j})^{-1}\eps_{\ol{2b-j},\ol{j-2b}} + (\ol{j-b})^{-1}\eps_{\ol{b},\ol{j-2b}}\right)\;\in\;\Xi
\]
for $j\in [0,p-1]\ohne\{ \ol{b}\}$.

\rm
{\it Proof.} First of all, we have $b + \ul{j-2b} \geq b + \ul{-2b} = b - 1 - \ul{2b - 1}$ and\vspace{2mm} 
$p(b-1 -\ul{2b - 1}) = (p-2)(b-1) + \ol{2b - 1} - 1 \geq 0$, both if $b = 1$ or if $b > 1$. So the element in question is in $\Lambda$.

By (\ref{CorFT14}), we have to prove that for $l\in [0,p-1]$,
\[
\val_s\!\left((-1)^{\ol{2b-j}}\smatze{l\ru{-2}}{\ol{2b-j}} (\ol{b-j})^{-1} + (-1)^{\ol{b}}\smatze{l\ru{-2}}{\ol{b}} (\ol{j-b})^{-1}\right) 
\;\geq\; 1 + \ul{bl - \ol{j-2b} - 1} - (b + \ul{j-2b})\; .
\]

If $l\in [0,p-2]$, then
\[ 
1 + \ul{bl - \ol{j-2b} - 1} \;\leq\; 1 + \ul{b(p-2) - \ol{j-2b} - 1} \= b + \ul{j-2b}\; .
\]
If $l = p-1$ then
\[
1 + \ul{b(p-1) - \ol{j-2b} - 1} - (b + \ul{j-2b}) \= - \ul{j - b}\; .
\]
Now $(-1)^h \smatze{p-1}{h} \con_p 1$ for $h\in [0,p-1]$ and $(\ol{b-j})^{-1} + (\ol{j-b})^{-1} \con_p 0$ together with 
$\val_s(p)\,\geq\, b-\ul{b}\,\geq\,  - \ul{j - b}$, both if $b = 1$ or if $b > 1$, yield the result.\qed
\end{Lemma}

Let
\[
\ba{rclcrcl}
\Xi & \lraa{\mu_j}  & \Xi & : & 1 & \lramaps & \mu_{0,j} \=  \sum_{i\in [0,p-1]} \eps_{i,j}              \\ 
\Xi & \lraa{\nu_j}  & \Xi & : & 1 & \lramaps & s^{b + \ul{j-2b}}
                                               \left((\ol{b-j})^{-1}\eps_{\ol{2b-j},\ol{j-2b}} + (\ol{j-b})^{-1}\eps_{\ol{b},\ol{j-2b}}\right)  \\
\ea
\]
be $\Xi$-linear maps for $j\in [0,p-1]\ohne\{\ol{b}\}$. By (\ref{LemAG9}), the map $\nu_j$ is welldefined.

\begin{Lemma}
\label{LemAG10}
For $j \in [0,p-1]\ohne\{ \ol{b}\}$, we obtain
\[
\ba{lcl}
\nu_j\be    & = & \al\mu_j \\
\mu_j\al    & = & \be\nu_j \\
\mu_j\chi_0 & = & \chi_j. \\
\ea
\]
That is, for $i\geq 0$, we obtain a representative
\begin{center}
\begin{picture}(1350,350)
\put(-100, 300){$\cdots$}
\put(   0, 310){\vector(1,0){80}}
\put( 100, 300){$\Xi$}
\put( 150, 310){\vector(1,0){130}}
\put( 210, 325){$\scm\be$}
\put( 300, 300){$\Xi$}
\put( 350, 310){\vector(1,0){130}}
\put( 410, 325){$\scm\al$}
\put( 500, 300){$\Xi$}
\put( 550, 310){\vector(1,0){130}}
\put( 610, 325){$\scm\be$}
\put( 700, 300){$\Xi$}
\put( 750, 310){\vector(1,0){130}}
\put( 810, 325){$\scm\al$}
\put( 829, 300){$\ob{\Xi}^\mb{\scr\rm degree $2i+1$}$}
\put( 950, 310){\vector(1,0){130}}
\put(1010, 325){$\scm\be$}
\put(1100, 300){$\Xi$}
\put(1150, 310){\vector(1,0){130}}
\put(1210, 325){$\scm\al$}
\put(1330, 300){$\cdots$}

\put( 115, 280){\vector(0,-1){130}}
\put( 130, 210){$\scm\mu_j$}
\put( 315, 280){\vector(0,-1){130}}
\put( 330, 210){$\scm\nu_j$}
\put( 515, 280){\vector(0,-1){130}}
\put( 530, 210){$\scm\mu_j$}
\put( 715, 280){\vector(0,-1){130}}
\put( 730, 210){$\scm\nu_j$}
\put( 915, 280){\vector(0,-1){130}}
\put( 870, 210){$\scm\mu_j$}
\put(1115, 280){\vector(0,-1){130}}

\put(-100, 100){$\cdots$}
\put(   0, 110){\vector(1,0){80}}
\put( 100, 100){$\Xi$}
\put( 150, 110){\vector(1,0){130}}
\put( 210, 125){$\scm\al$}
\put( 300, 100){$\Xi$}
\put( 350, 110){\vector(1,0){130}}
\put( 410, 125){$\scm\be$}
\put( 500, 100){$\Xi$}
\put( 550, 110){\vector(1,0){130}}
\put( 610, 125){$\scm\al$}
\put( 700, 100){$\Xi$}
\put( 750, 110){\vector(1,0){130}}
\put( 810, 125){$\scm\be$}
\put( 900, 100){$\Xi$}
\put( 950, 110){\line(1,0){27}}
\put( 997, 110){\vector(1,0){83}}
\put(1105, 100){$0$}
\put(1150, 110){\vector(1,0){130}}
\put(1330, 100){$\cdots$}

\put( 930, 280){\vector(1,-3){80}}
\put( 975, 190){$\scm\chi_j$}
\put( 940,  90){\vector(1,-1){50}}
\put( 920,  55){$\scm \chi_0$}
\put(1000,   0){$T\;$.}
\end{picture}
\end{center}
in $\liu{\KK^-(\Xi)}{(P,P[2i+1])}$ of $\chi_j^{(2i+1)}\in\Ext^{2i+1}_{T\wr G}(T,T)$.

Moreover, for $i\geq 0$, we obtain a representative
\begin{center}
\begin{picture}(1350,350)
\put(-100, 300){$\cdots$}
\put(   0, 310){\vector(1,0){80}}
\put( 100, 300){$\Xi$}
\put( 150, 310){\vector(1,0){130}}
\put( 210, 325){$\scm\al$}
\put( 300, 300){$\Xi$}
\put( 350, 310){\vector(1,0){130}}
\put( 410, 325){$\scm\be$}
\put( 500, 300){$\Xi$}
\put( 550, 310){\vector(1,0){130}}
\put( 610, 325){$\scm\al$}
\put( 700, 300){$\Xi$}
\put( 750, 310){\vector(1,0){130}}
\put( 810, 325){$\scm\be$}
\put( 855, 300){$\ob{\Xi}^\mb{\scr\rm degree $2i$}$}
\put( 950, 310){\vector(1,0){130}}
\put(1010, 325){$\scm\al$}
\put(1100, 300){$\Xi$}
\put(1150, 310){\vector(1,0){130}}
\put(1210, 325){$\scm\be$}
\put(1330, 300){$\cdots$}

\put( 115, 280){\vector(0,-1){130}}
\put( 130, 210){$\scm 1_\Xi$}
\put( 315, 280){\vector(0,-1){130}}
\put( 330, 210){$\scm 1_\Xi$}
\put( 515, 280){\vector(0,-1){130}}
\put( 530, 210){$\scm 1_\Xi$}
\put( 715, 280){\vector(0,-1){130}}
\put( 730, 210){$\scm 1_\Xi$}
\put( 915, 280){\vector(0,-1){130}}
\put( 870, 210){$\scm 1_\Xi$}
\put(1115, 280){\vector(0,-1){130}}

\put(-100, 100){$\cdots$}
\put(   0, 110){\vector(1,0){80}}
\put( 100, 100){$\Xi$}
\put( 150, 110){\vector(1,0){130}}
\put( 210, 125){$\scm\al$}
\put( 300, 100){$\Xi$}
\put( 350, 110){\vector(1,0){130}}
\put( 410, 125){$\scm\be$}
\put( 500, 100){$\Xi$}
\put( 550, 110){\vector(1,0){130}}
\put( 610, 125){$\scm\al$}
\put( 700, 100){$\Xi$}
\put( 750, 110){\vector(1,0){130}}
\put( 810, 125){$\scm\be$}
\put( 900, 100){$\Xi$}
\put( 950, 110){\line(1,0){27}}
\put( 997, 110){\vector(1,0){83}}
\put(1105, 100){$0$}
\put(1150, 110){\vector(1,0){130}}
\put(1330, 100){$\cdots$}

\put( 930, 280){\vector(1,-3){80}}
\put( 975, 190){$\scm\chi_0$}
\put( 940,  90){\vector(1,-1){50}}
\put( 920,  55){$\scm \chi_0$}
\put(1000,   0){$T\;$.}
\end{picture}
\end{center}
in $\liu{\KK^-(\Xi)}{(P,P[2i])}$ of $\chi_0^{(2i)}\in\Ext^{2i}_{T\wr G}(T,T)$. 

\rm
{\it Proof.} We claim that $\nu_j\be = \al\mu_j$. On the one hand, we obtain
\[
\barcl
1_\Xi(\nu_j \be)
& = & s^{\ul{b} + b + \ul{j - 2b}}\sum_{i\in [0,p-1]} \left( i(\ol{b-j})^{-1}\eps_{i,\ol{b}}\,\eps_{\ol{2b-j},\ol{j-2b}}  
                                                                      + i(\ol{j-b})^{-1}\eps_{i,\ol{b}}\,\eps_{\ol{b},  \ol{j-2b}}\right) \\
& \aufgl{\mb{\scr (\ref{LemFT11})}} & 
      s^{\ul{b} + b + \ul{j - 2b}}\sum_{i\in [0,p-1]} i (\ol{b-j})^{-1}\dell_{\,\ol{2b-j},\ol{i+b}}\,\eps_{i,\ol{b} + \ol{j - 2b}} \\
& = & s^{\ul{b} + b + \ul{j - 2b}}\,\eps_{\ol{b-j},\ol{b} + \ol{j - 2b}} \\
& = & s^{b - \ul{b} - 1} \eps_{\ol{b-j},j+p-\ol{b}}\; . \\
\ea
\]
On the other hand, we obtain
\[
\barcl
1_\Xi(\al\mu_j)
& = & s^{b - \ul{b} - 1} \sum_{i\in [0,p-1]} \eps_{i,j}\,\eps_{\ol{b},p-\ol{b}} \\
& \aufgl{\mb{\scr (\ref{LemFT11})}} & s^{b - \ul{b} - 1} \sum_{i\in [0,p-1]} \dell_{\,\ol{b},\ol{i+j}}\,\eps_{i,j+p-\ol{b}}  \\
& = & s^{b - \ul{b} - 1} \eps_{\ol{b-j},j+p-\ol{b}}\; . \\
\ea
\]
We claim that $\mu_j\al = \be\nu_j$. On the one hand, we obtain
\[
\barcl
1_\Xi(\be\nu_j)
& = & s^{\ul{b} + b + \ul{j - 2b}}\sum_{i\in [0,p-1]} \left(i(\ol{b-j})^{-1}\eps_{\ol{2b-j},\ol{j-2b}}\,\eps_{i,\ol{b}} 
                                                                    + i(\ol{j-b})^{-1}\eps_{\ol{b},   \ol{j-2b}}\,\eps_{i,\ol{b}}\right) \\
& \aufgl{\mb{\scr (\ref{LemFT11})}} & 
      s^{\ul{b} + b + \ul{j - 2b}}\sum_{i\in [0,p-1]} i(\ol{j-b})^{-1}\dell_{i,\ol{j-b}}\,\eps_{\ol{b},\ol{b}+\ol{j - 2b}} \\
& = & s^{\ul{b} + b + \ul{j - 2b}} \,\eps_{\ol{b},\ol{b}+\ol{j - 2b}} \\
& = & s^{b - \ul{b} - 1} \eps_{\ol{b},j-\ol{b}+p}\; . \\
\ea
\]
On the other hand, we obtain
\[
\barcl
1_\Xi(\mu_j\al)
& = & s^{b - \ul{b} - 1} \sum_{i\in [0,p-1]} \eps_{\ol{b},p-\ol{b}}\eps_{i,j} \\
& \aufgl{\mb{\scr (\ref{LemFT11})}} & s^{b - \ul{b} - 1} \sum_{i\in [0,p-1]} \dell_{i,0}\eps_{\ol{b},j-\ol{b}+p}  \\
& = & s^{b - \ul{b} - 1} \eps_{\ol{b},j-\ol{b}+p} \; . \\
\ea
\]
We claim that $\mu_j\chi_0 = \chi_j$. In fact,
\[
\barcl
1_\Xi(\mu_j\chi_0)
& = & \sum_{i\in [0,p-1]} t^0\eps_{i,j} \\
& = & \sum_{i\in [0,p-1]} \dell_{i,0} t^j \\
& = & t^j  \\
& = & 1_\Xi\,\chi_j\; . \\
\ea
\]
\qedfl{-45}
\end{Lemma}

\bq
\begin{Example}[cont.]
\label{ExHelp5_5}\rm
If $\Xi\iso\Z_{(3)}[\pi_2]\wr C_3$, then $b = \ol{b} = 1$; $\nu_0 = \smateckdd{0}{\; 0\;}{0}{0}{0}{1/2}{3}{0}{0}\cdot (-)$, \linebreak
$\mu_0 = \smateckdd{1}{0}{0}{0}{1}{0}{0}{0}{1}\cdot (-)$; $\nu_2 = \smateckdd{3/2}{0}{\; 0\;}{0}{3}{0}{0}{0}{0}\cdot (-)$, 
$\mu_2 = \smateckdd{0}{0}{1}{3}{0}{0}{0}{3}{0}\cdot (-)$.
\end{Example}
\eq

\begin{Theorem}
\label{ThAG11}
We have isomorphisms of graded $S$-algebras
\[
\fbox{$\;\;\ru{27}
\barcl
\fracd{S\Big[h_0^{(1)},\dots,h_{\ol{b}-1}^{(1)};h_{\ol{b}+1}^{(1)},\dots,h_{p-1}^{(1)};h_0^{(2)}\Big]}
{\left(\ba{l}\ru{5}
s^{\ul{b}+1}h_0^{(1)},\dots,s^{\ul{b}+1} h_{\ol{b}-1}^{(1)};s^{\ul{b}} h_{\ol{b}+1}^{(1)},\dots,s^{\ul{b}} h_{p-1}^{(1)}; s^{b - \ul{b}}h_0^{(2)}; 
\vspace*{3mm}\\ 
h_j^{(1)}\cdot h_k^{(1)} - \dell_{\,\ol{j+k},\ol{2b}}\, s^{b + \ul{j+k-2b}}\, (\ol{b-j})^{-1} h_0^{(2)} \vspace*{1mm}\\
\hspace*{32mm}\mb{\rm for $j,k\in [0,p-1]\ohne\{ \ol{b}\}$}\\
\ea\right)}
& \lraiso & \Ext^\ast_{T\wr G}(T,T) \; \lraiso \; \HH^\ast(G,T;S) \vspace*{2mm} \\
h^{(1)}_j & \lramaps & \chi_j^{(1)}\; ,\;\;\;\; j\in [0,p-1]\ohne\{ \ol{b}\} \vspace*{1mm} \\
h^{(2)}_0 & \lramaps & \chi_0^{(2)}\; , \\
\ea
\;\;$}
\]
as quotient of the {\rm graded commutative} polynomial ring 
$S[h_0^{(1)},\dots,h_{\ol{b}-1}^{(1)};h_{\ol{b}+1}^{(1)},\dots,h_{p-1}^{(1)};h_0^{(2)}]$ with grading determined by 
$\deg h^{(1)}_j = 1$ for $j\in [0,p-1]\ohne\{\ol{b}\}$ and $\deg h^{(2)}_0 = 2$.

\rm
{\it Proof.} The isomorphism $\Ext^\ast_{T\wr G}(T,T) \lraiso \HH^\ast(G,T;S)$ of graded $S$-algebras with respect to the Yoneda product resp.\ to the cup 
product is a consequence of (\ref{PropCup4}).

We shall exhibit the ring structure on the graded $S$-module $\Ext^\ast_{T\wr G}(T,T)$ (cf.\ \ref{PropAG8}).
By (\ref{LemAG10}), we obtain
\[
\ba{rclclcl}
\chi_0^{(2i)} & \cdot & \chi_0^{(2j)}   & = & \chi_0^{(2i+2j)}   & \in & \Ext_{T\wr G}^{2i+2j}(T,T) \\
\chi_0^{(2i)} & \cdot & \chi_k^{(2j+1)} & = & \chi_k^{(2i+2j+1)} & \in & \Ext_{T\wr G}^{2i+2j+1}(T,T) \\
\ea
\]
for $i,\, j\,\geq\, 0$ and $k\in [0,p-1]\ohne\{\ol{b}\}$. It remains to calculate $\chi_j^{(1)}\cdot\chi_k^{(1)}\,\in\,\Ext^2_{T\wr G}(T,T)$ for 
$j,\,k\,\in\, [0,p-1]\ohne\{ \ol{b}\}$. Using (\ref{LemAG10}) to represent $\chi_j$ in $\liu{\KK^-(\Xi)}{(P,P[1])}$, this product is represented by 
the $2$-cocycle $\nu_j\chi_k\in\liu{\Xi}{(\Xi,T)}$. We obtain
\[
\barcl
1_\Xi(\nu_j\chi_k)
& = & t^k s^{b + \ul{j-2b}}\,\left((\ol{b-j})^{-1}\eps_{\ol{2b-j},\ol{j-2b}} + (\ol{j-b})^{-1}\eps_{\ol{b},\ol{j-2b}}\right)\vspace*{1mm}\\
& = & t^k s^{b + \ul{j-2b}}\,(\ol{b-j})^{-1}\eps_{\ol{2b-j},\ol{j-2b}}\vspace*{1mm}\\
& = & s^{b + \ul{j-2b}}\,(\ol{b-j})^{-1}\dell_{k,\ol{2b-j}}\,t^{\ol{k+j-2b}}s^{\ul{k+\ol{j-2b}}}\vspace*{1mm}\\
& = & s^{b + \ul{j+k-2b}}\,(\ol{b-j})^{-1}\dell_{\,\ol{j+k},\ol{2b}}\, t^0\; ,\vspace*{1mm}\\
\ea
\]
i.e.\ $\nu_j\chi_k = \dell_{\,\ol{j+k},\ol{2b}}\, s^{b + \ul{j+k-2b}}\, (\ol{b-j})^{-1}\chi_0$. Hence
\[
\chi_j^{(1)}\cdot\chi_k^{(1)}\= \dell_{\,\ol{j+k},\ol{2b}}\, s^{b + \ul{j+k-2b}}\, (\ol{b-j})^{-1} \chi_0^{(2)}\; .
\]
\qedfl{-70}
\end{Theorem}

\bq
 \begin{Remark}
 \label{RemAG12}\rm
 Since $\ol{j+k} = \ol{2b}$ implies $(\ol{b-j})^{-1} \con_p -(\ol{b-k})^{-1}$, and since 
 \mb{$\val_s(p)\geq b - \ul{b}$,} we obtain the graded commutativity of $\Ext^\ast_{T\wr G}(T,T)$ without reverting to the graded commutativity of 
 the cup product.

 For instance, if $b = p + 1$, then $\chi_0^{(1)}\cdot\chi_2^{(1)} = s^{p-1}\chi_0^{(2)}\neq 0$. In particular, $\Ext^\ast_{T\wr G}(T,T)$ is not 
 commutative in this case.
 \end{Remark}
\eq

\begin{Corollary}
\label{CorAG13}
If $b = 1$, then we have isomorphisms of graded $S$-algebras
\[
S[h^{(1)},h^{(2)}]/(sh^{(1)},\, sh^{(2)},\, (h^{(1)})^2) \;\lraiso\; \Ext^\ast_{T\wr G}(T,T) \;\lraiso\; \HH^\ast(G,T;S)\; ,
\]
as quotient of the {\rm commutative} polynomial ring $S[h^{(1)},h^{(2)}]$ with grading determined by $\deg h^{(1)} = 1$, $\deg h^{(2)} = 2$. 

\rm
{\it Proof.} In fact, we have $\ul{b} = 0$. There are no nonzero products of homogeneous elements of odd degree, so we may use the commutative polynomial ring.\qed
\end{Corollary}

\bq
 \begin{Example}[cont.]
 \label{ExHelp6}\rm
 We have $\HH^\ast(C_3,\Z_{(3)}[\pi_2];\Z_{(3)}) \iso \Z_{(3)}[h^{(1)},h^{(2)}]/(3h^{(1)},3h^{(2)},(h^{(1)})^2)$. Note that 
 $\HH^0(C_3,\Z_{(3)}[\pi_2];\Z_{(3)}) \iso \Z_{(3)}$, and that $\HH^i(C_3,\Z_{(3)}[\pi_2];\Z_{(3)})\iso\Fu{3}$ for $i\geq 1$.
 \end{Example}
\eq
\section{Applications}
\label{SecAppl}

\bq
We give some applications, refraining, however, from a repetition of (\ref{ThAG11}) in different instances.
\eq

\subsection{Lubin-Tate extensions}
\label{SubsecLT}

\bq
The results apply to certain of the extensions of local fields described by {\sc Lubin} and {\sc Tate} {\bf\cite{LT65}}. An introduction to their 
theory is also given in {\bf\cite[\rm p.\ 146 ff.]{Se62}}.
\eq

Let $p \geq 3$ be a prime. Let $B$ be a local field with discrete valuation ring $R$, whose maximal ideal is generated by $\pi$. Assume that 
$R/\pi\iso\Fu{p}$. We choose the Lubin-Tate series $f(X) = X^p + \pi X \in R[[X]]$, and obtain the unique commutative formal group
\[
\barcl
F(X,Y)
& = & X+Y \\
& - & (\pi - \pi^p)^{-1}((X+Y)^p - (X^p + Y^p)) \\
& - & p(\pi - \pi^p)^{-1}(\pi - \pi^{2p - 1})^{-1}
\left(
\ba{rl}
  & \pi^{p-1}(X+Y)^{p-1}(X^p + Y^p) \\
- & (X+Y)^{p-1}((X+Y)^p - (X^p + Y^p)) \\ 
- & \pi^{p-1}(X^{2p-1} + Y^{2p-1}) \\
\ea
\right)
\\
& + & \Ol(\mb{degree $3p-2$}) \;\in\; R[[X,Y]] \\
\ea
\]
such that $F(f(X),f(Y)) = f(F(X,Y))$. There is an injective ring morphism
\[
\barcl
R & \lra     & \End F \\
a & \lramaps & [a](X)\; , \\
\ea
\]
where $[a](X)\in R[[X]]$ is uniquely determined by $[a](X)\con_{X^2} aX$ and the endomorphism property $F([a](X),[a](Y)) = [a](F(X,Y))$. So for 
instance, $[\pi](X) = f(X) = X^p + \pi X$. We write $P_n(X) := [\pi^n](X)\in R[X]$ for $n\geq 0$, so that $P_0(X) = X$, 
$P_1(X) = f(X) = X^p + \pi X$ and $P_n(X) = P_{n-1}(X)^p + \pi P_{n-1}(X)$. Moreover, $P_n(0) = 0$, $\deg P_n(X) = p^n$, $P_n(X)\con_\pi X^{p^n}$
and $P_n'(X) \con_p \pi^n$.

Let $\b B$ be an algebraic closure of $B$, and let $\b\mfk = \{ x\in\b B : \Nrm_{B(x)|B}(x)\in R\pi\}\tm \b B$
be the maximal ideal of its valuation ring, which becomes an abelian group $(\b\mfk,\ast)$ via $x\ast y := F(x,y)$. Moreover, $\b\mfk$ 
becomes an $R$-module via
\[
\barcl
R & \lra     & \End (\b\mfk,\ast) \\
a & \lramaps & (\; [a]\; :\; x\lramaps [a]\cdot x := [a](x)\;)\; . \\
\ea
\]
For $n\geq 1$, we let
\[
\mufat_n \; :=\; \mb{ann}_{[\pi^n]} \b\mfk\; \; =\; \{x\in\b K\; :\; P_n(x) = 0\}\; .
\]
By separability of $P_n(X)$, we have $\#\mufat_n = p^n$ for each $n\geq 1$, whence $\mufat_n \iso R/\pi^n$ as $R$-modules. Let 
$\tht_n$ be an $R$-linear generator of $\mufat_n$, chosen in such a way that $[\pi](\tht_n) = \tht_{n-1}$. We have 
$\mu_{\tht_n,B}(X) = P_n(X)/P_{n-1}(X) = P_{n-1}(X)^{p-1} + \pi$, whence $B(\mufat_n) = B(\tht_n)$ is galois over $B$ with 
\[
\barcl
(R/\pi^n)^\ast & \lraiso  & \Gal(B(\mufat_n)|B) \\
u              & \lramaps & (\;\spi{u}\; :\; \tht_n\lramaps [u](\tht_n)\;)\; . \\
\ea
\]
Now $R[\tht_n]|R$ is purely ramified, and as different we obtain
\[
\barcl
\Dfk_{R[\tht_n]|R} 
& = & \left(\mu_{\tht_n,B}'(\tht_n)\right) \vspace*{1mm}\\
& = & \left(P_n'(\tht_n)/P_{n-1}(\tht_n)\right) \vspace*{1mm}\\
& = & \left(P_{n-1}'(\tht_n)(p \tht_1^{p-1} + \pi)/\tht_1\right) \vspace*{1mm} \\
& = & \left(P_{n-2}'(\tht_n)(p \tht_2^{p-1} + \pi)(p \tht_1^{p-1} + \pi)/\tht_1\right)\vspace*{1mm} \\
& = & \cdots \\
& = & \left(\tht_1^{-1}\prod_{i\in [1,n]} (p \tht_i^{p-1} + \pi)\right) \\
& = & \left( \tht_1^{-1}\pi^n\right)\; , \\
\ea
\]
whence $\Dfk_{R[\tht_n]|R[\tht_{n-1}]} = (\pi) = (\tht_n^{p^{n-1}(p-1)})$. 

\begin{Example}
\label{ExLT1}\rm
Let $n\geq 2$. We may apply (\ref{ThFT16}, \ref{ThAG11}) to 
\[
\ba{rclcrcl}
S      & = & R[\tht_{n-1}]   & \Icm & s      & = & \tht_{n-1} \\
T      & = & R[\tht_n]       &      & t      & = & \tht_n \\
b      & = & p^{n-1} - 1     &      &        &   & \\
\ul{b} & = & p^{n-2} - 1     &      & \ol{b} & = & p - 1\; .\\
\ea
\]
The value of $b$ results from the different by the formula $\val_t(\Dfk_{T|S}) = (p-1)(b+1)$ (\ref{RemFT13_1}). 
Moreover, $\val_s(p)\geq\val_s(\pi) = p^{n-2}(p-1) = b - \ul{b}\,$.
\end{Example}

We have $\Fu{p}^\ast\hra (R/\pi^n)^\ast$ by sending $j\lramaps j^{p^{n-1}}$. For $n\geq 1$, we let 
\[
\pi_n \; :=\; \prod_{j\in \Fu{p}^\ast} \tht_n^{\spi{j^{p^{n-1}}}} 
      \; =\; \prod_{j\in [1,p-1]} [j^{p^{n-1}}](\tht_n) \; .
\]
Then $R[\pi_n]$ is purely ramified over $R$, of degree $p^{n-1}$ and with maximal ideal generated by $\pi_n$. In particular, $\pi_1 = \pi$. 
Moreover, $\pi_n = \Nrm_{B(\tht_n)|B(\pi_n)}(\tht_n)$.

As different, we obtain
\[
\barcl
\Dfk_{R[\pi_n]|R[\pi_{n-1}]} 
& = & \Dfk_{R[\tht_n]|R[\pi_n]}^{-1}\Dfk_{R[\tht_n]|R[\tht_{n-1}]}\Dfk_{R[\tht_{n-1}]|R[\pi_{n-1}]} \\
& = & (\tht_n^{p-2})^{-1} (\tht_n^{p^{n-1}(p-1)}) (\tht_n^{p(p-2)}) \\
& = & (\pi_n^{p^{n-1} + p - 2}) \\
\ea
\]
{\bf\cite[\rm III.\S3, prop.\ 13]{Se62}}.

\begin{Example}
\label{ExA4}\rm
Let $n\geq 2$, $p\geq 3$. We may apply (\ref{ThFT16}, \ref{ThAG11}, \ref{PropND7}) to 
\[
\ba{rclcrcl}
S      & = & R[\pi_{n-1}]        & \Icm & s      & = & \pi_{n-1} \\
T      & = & R[\pi_n]            &      & t      & = & \pi_n \\
U      & = & R[\tht_n]           &      & u      & = & \tht_n \\
b      & = & (p^{n-1} - 1)/(p-1) &      &        &   & \\
\ul{b} & = & (p^{n-2} - 1)/(p-1) &      & \ol{b} & = & 1 \; ,\\
\ea
\]
The value of $b$ results from the different by the formula $\val_t(\Dfk_{T|S}) = (p-1)(b+1)$ (\ref{RemFT13_1}). 
Moreover, $\val_s(p)\geq\val_s(\pi) = p^{n-2} = b - \ul{b}\,$. 

In particular, (\ref{PropCup4}, \ref{PropND7}) yield isomorphisms of graded $S$-algebras
\[
\HH^\ast(C_p\ti C_{p-1},U;S)\;\iso\; \Ext_{U\wr (C_p\ti C_{p-1})}^\ast(U,U)\;\iso\; \Ext_{T\wr C_p}^\ast(T,T)\;\iso\;\HH^\ast(C_p,T;S)\; .
\] 
\end{Example}

\subsection{Cyclotomic number field extensions}
\label{SubsecCNF}

\bq
Passing to completions without changing cohomology, we may consider cyclotomic number field extensions as particular Lubin-Tate extensions. For 
sake of illustration, we recall the cyclotomic framework; in it, there is no need for completion, since the formal group law is given by the 
polynomial $F(X,Y) = X + Y + XY$. Strictly speaking, since in section \ref{SubsecCNF} we choose a different Lubin-Tate series as in 
section \ref{SubsecLT}, viz.\ $(X+1)^p - 1$ instead of $X^p + pX$, we are not directly specializing to this cyclotomic case. So keeping the 
notation of section \ref{SubsecLT} is a slight abuse.
\eq

Let $p\geq 3$ be a prime. For $n \geq 1$, we let $\zeta_{p^n}$ be a primitive $p^n$th root of unity over $\Q$. We make choices in such a manner that 
$\zeta_{p^n}^p = \zeta_{p^{n-1}}$ for $n\geq 2$ and denote $\tht_n := \zeta_{p^n} - 1$.  Let 
\[
\pi_n \; =\; \prod_{j\in [1,\, p-1]} (\zeta_{p^n}^{j^{p^{n-1}}} - 1)\; .
\]
Then $\Q(\tht_n) = \Q(\zeta_{p^n})$, $\Q(\pi_n) = \Fix_{C_{p-1}} \Q(\tht_n)$ and $\pi_n = \Nrm_{\Q(\tht_n)|\Q(\pi_n)}(\tht_n)$.

We have $\Nrm_{\Q(\tht_n)|\Q(\tht_{n-1})}(\tht_n) = \tht_{n-1}$ and $\Nrm_{\Q(\pi_n)|\Q(\pi_{n-1})}(\pi_n) = \pi_{n-1}$. Note that $\pi_1 = p$. 

The integral closure of $\Z_{(p)}$ in $\Q(\tht_n)$ is given by the discrete valuation ring $\Z_{(p)}[\tht_n]$, with maximal ideal generated by 
$\tht_n$, purely ramified over $\Z_{(p)}$; the integral closure of $\Z_{(p)}$ in $\Q(\pi_n)$ is given by the discrete valuation ring 
$\Z_{(p)}[\pi_n]$, with maximal ideal generated by $\pi_n$, purely ramified over $\Z_{(p)}$. 

\begin{Example}
\label{ExCyc1}\rm
Let $n\geq 2$. We may apply (\ref{ThFT16}, \ref{ThAG11}) to 
\[
\ba{rclcrcl}
S      & = & \Z_{(p)}[\tht_{n-1}]   & \Icm & s      & = & \tht_{n-1} \\
T      & = & \Z_{(p)}[\tht_n]       &      & t      & = & \tht_n \\
b      & = & p^{n-1} - 1            &      &        &   & \\
\ul{b} & = & p^{n-2} - 1            &      & \ol{b} & = & p - 1\; .\\
\ea
\]
We remark that $\val_s(p) = p^{n-2}(p-1) = b - \ul{b}$.
\end{Example}

\begin{Example}
\label{ExCyc2}\rm
Let $n\geq 2$, $p\geq 3$. We may apply (\ref{ThFT16}, \ref{ThAG11}, \ref{PropND7}) to 
\[
\ba{rclcrcl}
S      & = & \Z_{(p)}[\pi_{n-1}]   & \Icm & s      & = & \pi_{n-1} \\
T      & = & \Z_{(p)}[\pi_n]       &      & t      & = & \pi_n \\
U      & = & \Z_{(p)}[\tht_n]      &      & u      & = & \tht_n \\
b      & = & (p^{n-1} - 1)/(p-1)   &      &        &   & \\
\ul{b} & = & (p^{n-2} - 1)/(p-1)   &      & \ol{b} & = & 1 \; ,\\
\ea
\]
where $b$ is e.g.\ calculated using {\bf\cite[\rm VI.\S 1, prop.\ 3]{Se62}}. We remark that $\val_s(p) = p^{n-2} = b - \ul{b}$. In particular, 
(\ref{PropCup4}, \ref{PropND7}) yield
\[
\HH^\ast(C_p\ti C_{p-1},U;S)\;\iso\;\HH^\ast(C_p,T;S)\; .
\]
Hence, for instance,
\[
\HH^\ast((\Z/p^2)^\ast,\Z_{(p)}[\zeta_{p^2}];\Z_{(p)})\;\iso\; \Z_{(p)}[h^{(1)},h^{(2)}]/(ph^{(1)},\, ph^{(2)},\, (h^{(1)})^2)
\]
(cf.\ \ref{CorAG13}). 
\end{Example}

\bq
\begin{Remark}
\label{RemCyc3}\rm
In the same manner, (\ref{ThFT16}, \ref{ThAG11}, \ref{PropND7}) may be applied to certain cyclotomic function field extensions as defined by 
{\sc Carlitz} and {\sc Hayes} (cf.\ {\bf\cite{Ca38}}, {\bf\cite{Ha74}}; see also {\bf\cite[\rm sec.\ 6.1]{KW03}}). Up to completion, these also 
form a particular case of Lubin-Tate-extensions; again, there is no need for completion, the formal group law being given by $F(X,Y) = X + Y$.
\end{Remark}
\eq
\appendix

\section{The case $C_{p^2}$: a conjecture and an experiment}
\label{AppExp}

\bq
So far, we have essentially only treated the case of an extension with galois group $C_p$. The galois group $C_{p^2}$ seems to yield a somewhat more
involved twisted group ring, which we would like to illustrate in the case of $\Z_{(p)}[\pi_3]\wr C_{p^2}$. The calculations were carried out using
{\sc Magma} {\bf\cite{Magma}}.
\eq

\begin{footnotesize}
\subsection{The conjectural situation}

Suppose given a prime $p\geq 3$. We maintain the notation of section \ref{SubsecCNF} concerning $\pi_n$. Let 
\[
\ba{rclcrcl}
S & = & \Z_{(p)}\; ,        & \Icm & s & = & \pi_1 \= p\; , \\
T & = & \Z_{(p)}[\pi_2]\; , &      & t & = & \pi_2\; , \\ 
U & = & \Z_{(p)}[\pi_3]\; , &      & u & = & \pi_3\; , \\
\ea
\]
and let $G = C_{p^2}$ be generated by the restriction of $\zeta_{p^3}\lramapsa{\sa}\zeta_{p^3}^{p+1}$ from $\Q(\zeta_{p^3})$ to $\Q(\pi_3)$. 
(The role of $U$ in this appendix differs from the role of $U$ in section \ref{SecNebe}, where it has been a `cohomologically inessential' 
extension of $T$.)

For some peculiar reason, $t$ will not play a role at all. Instead, we consider a {\it Sen element} (cf.\ \linebreak {\bf\cite[\rm Lem.\ 1]{Sen69}})
\[
v \; := \; \prod_{i\in [0,p-1]} u^{\sa^i} \; .
\]

The $S$-linear colength of the Wedderburn embedding
\[
U\wr C_{p^2} \;\hraa{\omega}\; \Gamma \; :=\; \End_S U
\]
is $p^2(p^2 + (p^2 - p - 2)/2)$ {\bf\cite[\rm (2.17)]{K01}}. We fix the $S$-linear basis 
\[
(u^i v^j)_{i\in [0,p-1],\; j\in [0,p-1]} \= (u^0 v^0, u^0 v^1, \dots, u^0 v^{p-1}, u^1 v^0, \dots, u^1 v^{p-1}, \; \dots\; , 
u^{p-1} v^0,\dots, u^{p-1} v^{p-1}) 
\]
of $U$ with respect to which we represent elements of $\Gamma$ as matrices, i.e.\ by means of which we identify $\Gamma = S^{p^2\ti p^2}$. 

\begin{Remark}
\label{RemA1}
We have
\[
\barcl
\val_u(u^\sa - u) & = & 1 + 1 \\
\val_u(v^\sa - v) & = & 1 + 2p \; . \\
\ea
\]

\rm
{\it Proof.} The second congruence is equivalent to $\val_u(u^{\sa^p} - u) = 1 + (p + 1)$, so that both assertions follow from 
{\bf\cite[\rm VI.\S 1, prop.\ 3]{Se62}}.\qed
\end{Remark}

As usual, let $\Xi$ denote the image of the Wedderburn embedding 
\[ 
\barcl
U\wr C_{p^2} & \lraa{\omega} & \Gamma    \\
u            & \lramaps      & (\dot u :   x\lramaps xu) \\
\sa          & \lramaps      & (\dot\sa :  x\lramaps x^\sa)\; , \\
\ea
\]
and let 
\[
\Lambda \; :=\; \Big\{ f\in\Gamma \; :\; (U u^i) f \tm U u^i  \mb{\rm\ for all } i\in [0,p-1]\Big\} \;\tm\; \Gamma \; .
\]

By (\ref{RemA1}) we obtain the intermediate ring
\[
\ba{rcrrl}
\Xi 
& \tm & \sLD  & := & \Lambda((\dot u,\dot v), (p-1, p-1), (2,2p + 1))_{\dot u\Lambda} \\
&     &  & =  & \{f\in\Lambda\; :\; D_{\dot u}^i\circ D_{\dot v}^j(f) \con_{\dot u^{2i + (2p + 1)j}\Lambda} 0\mb{\rm\ for all } i,j\in [0,p-1] \} \\
& \tm & \Lambda\; , &   & 
\ea
\]
cf.\ (\ref{LemFT0_5}). Presumably, this is the smallest intermediate ring between $\Xi$ and $\Lambda$ that can be defined by derivations.

\begin{Conjecture}
\label{ConjA2}\Absit
\begin{itemize}
\item[{\rm (i)}] Given $\ta\in C_{p^2}$, we conjecture that
\[
2(v - v^\ta) + (u^\ta)^p (u - u^\ta) u^{p-1} + ((u^\ta)^{2p-1} - (u^\ta)^{2p})(u - u^\ta) \;\con_{u^{2p+3}}\; 0\; .
\]
\item[{\rm (ii)}] Moreover, we conjecture that 
\[
\sLDE \; :=\; \{ f\in\sLD\; :\; D_{\dot u}^i\circ E_{\dot u,\dot v}^j(f) \con_{\dot u^{2i + (2p + 3)j}\Lambda} 0\mb{\rm\ for all } i,j\in [0,p-1]\} 
\]
contains $\Xi$, where
\[
E_{\dot u,\dot v}(f) \; :=\; 2 D_{\dot v}(f) + \dot u^p D_{\dot u}(f) \dot u^{p-1} + (\dot u^{2p-1} - \dot u^{2p}) D_{\dot u}(f) \; .
\]
\end{itemize}
\end{Conjecture}

\begin{Remark}
\label{RemA4}\rm\Absit
\begin{itemize}
\item[(i)] Conjecture (\ref{ConjA2}.i) holds for $p\in \{ 3,5,7\}$.
\item[(ii)] If $p = 3$, we obtain the colengths
\[
\Xi \= \sLDE\;\auf{9\cdot (3 - 1)}{\tm}\; \sLD\;\auf{9\cdot (3\cdot (3+1)/2 - 1)}{\tm}\; \Lambda \;\auf{9\cdot (9-1)/2}{\tm}\; \Gamma\; .
\]
\item[(iii)] If $p = 5$, we obtain the colengths 
\[
\Xi\;\auf{25\cdot (5 - 1)(5-3)/2}{\tm}\;\sLDE\;\auf{25\cdot (5 - 1)}{\tm}\;\sLD\;\auf{25\cdot (5\cdot (5+1)/2 - 1)}{\tm}\;\Lambda
\;\auf{25\cdot (25-1)/2}{\tm}\;\Gamma\; . \Icm
\]
\item[(iv)] According to our wishful thinking, (\ref{ConjA2}.i) should be part of a series of congruences in $U$ that completely describes $\Xi$ 
by adding the according congruences in $\Lambda$ to the provisional definition of $\sLDE$ given in (\ref{ConjA2}.ii).
\item[(v)] Since $E_{\dot u,\dot v}^j$ is not a derivation, we do not know whether (\ref{ConjA2}.i) implies (\ref{ConjA2}.ii).
\end{itemize}
\end{Remark}

\subsection{Simplifying $\dot u$, $\dot v$}
\label{Subsecp3}

\subsubsection{The case $p = 3$}

Suppose $p = 3$, i.e.\ $S = \Z_{(3)}$, $u = \pi_3 = (\zeta_{27} - 1)(\zeta_{27}^{-1} - 1)$, 
$v = u\; u^\sa u^{\sa^2} = (\zeta_{27} - 1)(\zeta_{27}^{-1} - 1)(\zeta_{27}^4 - 1)\cdot$\linebreak 
$\cdot(\zeta_{27}^{-4} - 1)(\zeta_{27}^{16} - 1)(\zeta_{27}^{-16} - 1)$ and $U = \Z_{(3)}[\pi_3]$. 

With respect to the basis $(u^0 v^0, u^0 v^1, u^0 v^2, u^1 v^0, u^1 v^1, u^1 v^2, u^2 v^0, u^2 v^1, u^2 v^2)$, the multiplication by $u$ on $U$ 
is given by
\[
\ba{l}
\dot u \= 
\etwasenger
\left[
\ba{ccc|ccc|ccc}
\scm 0 &\scm 0 &\scm 0 &\scm 1 &\scm 0 &\scm 0 &\scm 0 &\scm 0 &\scm 0\\
\scm 0 &\scm 0 &\scm 0 &\scm 0 &\scm 1 &\scm 0 &\scm 0 &\scm 0 &\scm 0\\
\scm 0 &\scm 0 &\scm 0 &\scm 0 &\scm 0 &\scm 1 &\scm 0 &\scm 0 &\scm 0\\\hline\rule{0mm}{2mm}
\scm 0 &\scm 0 &\scm 0 &\scm 0 &\scm 0 &\scm 0 &\scm 1 &\scm 0 &\scm 0\\
\scm 0 &\scm 0 &\scm 0 &\scm 0 &\scm 0 &\scm 0 &\scm 0 &\scm 1 &\scm 0\\
\scm 0 &\scm 0 &\scm 0 &\scm 0 &\scm 0 &\scm 0 &\scm 0 &\scm 0 &\scm 1\\\hline\rule{0mm}{2mm}
\scm -\frac{2433}{7217} &\scm \frac{3119}{1031} &\scm -\frac{2151}{7217} &\scm -\frac{6249}{7217} &\scm -\frac{12225}{1031} &\scm \frac{1168}{1031} &\scm \frac{19680}{7217} &\scm \frac{25540}{7217} &\scm -\frac{306}{1031}\\
\scm -\frac{4050}{7217} &\scm \frac{3777}{1031} &\scm \frac{593}{7217} &\scm \frac{18804}{7217} &\scm -\frac{20163}{1031} &\scm \frac{987}{1031} &\scm -\frac{6885}{7217} &\scm \frac{64005}{7217} &\scm -\frac{278}{1031}\\
\scm -\frac{111}{7217} &\scm -\frac{1026}{1031} &\scm \frac{50145}{7217} &\scm \frac{18189}{7217} &\scm -\frac{14604}{1031} &\scm -\frac{12306}{1031} &\scm -\frac{6684}{7217} &\scm \frac{29259}{7217} &\scm \frac{6603}{1031} \\
\ea
\right]
\weiter\vspace*{2mm}\\
\;\con_{\dot u^6\Lambda}\;
\enger
\left[
\ba{rrr|rrr|rrr}
\scm 0 &\scm\;\; 0 &\scm\;\; 0 &\scm  1 &\scm  0 &\scm  0 &\scm  0 &\scm  0 &\scm  0 \\
\scm 0 &\scm     0 &\scm     0 &\scm  0 &\scm  1 &\scm  0 &\scm  0 &\scm  0 &\scm  0\\
\scm 0 &\scm     0 &\scm     0 &\scm  0 &\scm  0 &\scm  1 &\scm  0 &\scm  0 &\scm  0\\\hline\rule{0mm}{2mm}
\scm 0 &\scm     0 &\scm     0 &\scm  0 &\scm  0 &\scm  0 &\scm  1 &\scm  0 &\scm  0\\
\scm 0 &\scm     0 &\scm     0 &\scm  0 &\scm  0 &\scm  0 &\scm  0 &\scm  1 &\scm  0\\
\scm 0 &\scm     0 &\scm     0 &\scm  0 &\scm  0 &\scm  0 &\scm  0 &\scm  0 &\scm  1\\\hline\rule{0mm}{2mm}
\scm 0 &\scm     1 &\scm     0 &\scm  0 &\scm  0 &\scm -1 &\scm  0 &\scm -1 &\scm  0\\
\scm 0 &\scm     0 &\scm     1 &\scm -3 &\scm  0 &\scm  0 &\scm  0 &\scm  0 &\scm -1\\
\scm 3 &\scm     0 &\scm     0 &\scm  0 &\scm -3 &\scm  0 &\scm -3 &\scm  0 &\scm  0\\
\ea
\right]
\; =: \; \dddot{u}
\;\con_{\dot u^3\Lambda}\;
\enger
\left[
\ba{rrr|rrr|rrr}
\scm 0 &\scm 0 &\scm 0 &\scm 1 &\scm 0 &\scm 0 &\scm 0 &\scm 0 &\scm 0\\
\scm 0 &\scm 0 &\scm 0 &\scm 0 &\scm 1 &\scm 0 &\scm 0 &\scm 0 &\scm 0\\
\scm 0 &\scm 0 &\scm 0 &\scm 0 &\scm 0 &\scm 1 &\scm 0 &\scm 0 &\scm 0\\\hline\rule{0mm}{2mm}
\scm 0 &\scm 0 &\scm 0 &\scm 0 &\scm 0 &\scm 0 &\scm 1 &\scm 0 &\scm 0\\
\scm 0 &\scm 0 &\scm 0 &\scm 0 &\scm 0 &\scm 0 &\scm 0 &\scm 1 &\scm 0\\
\scm 0 &\scm 0 &\scm 0 &\scm 0 &\scm 0 &\scm 0 &\scm 0 &\scm 0 &\scm 1\\\hline\rule{0mm}{2mm}
\scm 0 &\scm 1 &\scm 0 &\scm 0 &\scm 0 &\scm 0 &\scm 0 &\scm 0 &\scm 0\\
\scm 0 &\scm 0 &\scm 1 &\scm 0 &\scm 0 &\scm 0 &\scm 0 &\scm 0 &\scm 0\\
\scm 3 &\scm 0 &\scm 0 &\scm 0 &\scm 0 &\scm 0 &\scm 0 &\scm 0 &\scm 0\\
\ea
\right]
\weiter
=: \ddot u\; , \\
\ea
\]
the multiplication by the Sen-element $v$ by
\[
\ba{l}
\dot v \=
\etwasenger
\left[
\ba{ccc|ccc|ccc}
\scm 0 &\scm 1 &\scm 0 &\scm 0 &\scm 0 &\scm 0 &\scm 0 &\scm 0 &\scm 0\\
\scm 0 &\scm 0 &\scm 1 &\scm 0 &\scm 0 &\scm 0 &\scm 0 &\scm 0 &\scm 0\\
\scm \frac{2256}{1031} &\scm -\frac{19548}{1031} &\scm \frac{31788}{1031} &\scm \frac{2484}{1031} &\scm -\frac{3456}{1031} &\scm -\frac{31626}{1031} &\scm -\frac{495}{1031} &\scm -\frac{7083}{1031} &\scm \frac{10224}{1031}\\\hline\rule{0mm}{2mm}
\scm 0 &\scm 0 &\scm 0 &\scm 0 &\scm 1 &\scm 0 &\scm 0 &\scm 0 &\scm 0\\
\scm 0 &\scm 0 &\scm 0 &\scm 0 &\scm 0 &\scm 1 &\scm 0 &\scm 0 &\scm 0\\
\scm \frac{27891}{7217} &\scm -\frac{37620}{1031} &\scm \frac{494226}{7217} &\scm \frac{69981}{7217} &\scm -\frac{19980}{1031} &\scm -\frac{97587}{1031} &\scm -\frac{11043}{7217} &\scm -\frac{186021}{7217} &\scm \frac{35910}{1031}\\\hline\rule{0mm}{2mm}
\scm 0 &\scm 0 &\scm 0 &\scm 0 &\scm 0 &\scm 0 &\scm 0 &\scm 1 &\scm 0\\
\scm 0 &\scm 0 &\scm 0 &\scm 0 &\scm 0 &\scm 0 &\scm 0 &\scm 0 &\scm 1\\
\scm \frac{104247}{7217} &\scm -\frac{137862}{1031} &\scm \frac{1734570}{7217} &\scm \frac{186300}{7217} &\scm -\frac{7866}{1031} &\scm -\frac{385245}{1031} &\scm -\frac{15474}{7217} &\scm -\frac{809595}{7217} &\scm \frac{140031}{1031}\\
\ea
\right]
\weiter\vspace*{2mm}\\
\;\con_{\dot u^{9 + 3\cdot (3+1)/2}\Lambda}\;
\enger
\left[
\ba{rrr|rrr|rrr}
\scm 0 &\scm 1 &\scm 0 &\scm 0 &\scm 0 &\scm 0 &\scm 0 &\scm 0 &\scm 0\\
\scm 0 &\scm 0 &\scm 1 &\scm 0 &\scm 0 &\scm 0 &\scm 0 &\scm 0 &\scm 0\\
\scm 3 &\scm 0 &\scm 0 &\scm 0 &\scm 0 &\scm 0 &\scm 9 &\scm 0 &\scm 0\\\hline\rule{0mm}{2mm}
\scm 0 &\scm 0 &\scm 0 &\scm 0 &\scm 1 &\scm 0 &\scm 0 &\scm 0 &\scm 0\\
\scm 0 &\scm 0 &\scm 0 &\scm 0 &\scm 0 &\scm 1 &\scm 0 &\scm 0 &\scm 0\\
\scm 0 &\scm 9 &\scm 0 &\scm 3 &\scm 0 &\scm 0 &\scm 0 &\scm 0 &\scm 0\\\hline\rule{0mm}{2mm}
\scm 0 &\scm 0 &\scm 0 &\scm 0 &\scm 0 &\scm 0 &\scm 0 &\scm 1 &\scm 0\\
\scm 0 &\scm 0 &\scm 0 &\scm 0 &\scm 0 &\scm 0 &\scm 0 &\scm 0 &\scm 1\\
\scm 0 &\scm 0 &\scm 0 &\scm 0 &\scm 9 &\scm 0 &\scm 3 &\scm 0 &\scm 0\\
\ea
\right]
\weiter
\; =: \; \dddot v 
\;\con_{\dot u^{9+3}\Lambda}\; 
\enger
\left[
\ba{rrr|rrr|rrr}
\scm 0 &\scm 1 &\scm 0 &\scm 0 &\scm 0 &\scm 0 &\scm 0 &\scm 0 &\scm 0\\
\scm 0 &\scm 0 &\scm 1 &\scm 0 &\scm 0 &\scm 0 &\scm 0 &\scm 0 &\scm 0\\
\scm 3 &\scm 0 &\scm 0 &\scm 0 &\scm 0 &\scm 0 &\scm 0 &\scm 0 &\scm 0\\\hline\rule{0mm}{2mm}
\scm 0 &\scm 0 &\scm 0 &\scm 0 &\scm 1 &\scm 0 &\scm 0 &\scm 0 &\scm 0\\
\scm 0 &\scm 0 &\scm 0 &\scm 0 &\scm 0 &\scm 1 &\scm 0 &\scm 0 &\scm 0\\
\scm 0 &\scm 0 &\scm 0 &\scm 3 &\scm 0 &\scm 0 &\scm 0 &\scm 0 &\scm 0\\\hline\rule{0mm}{2mm}
\scm 0 &\scm 0 &\scm 0 &\scm 0 &\scm 0 &\scm 0 &\scm 0 &\scm 1 &\scm 0\\
\scm 0 &\scm 0 &\scm 0 &\scm 0 &\scm 0 &\scm 0 &\scm 0 &\scm 0 &\scm 1\\
\scm 0 &\scm 0 &\scm 0 &\scm 0 &\scm 0 &\scm 0 &\scm 3 &\scm 0 &\scm 0\\
\ea
\right]
\weiter
\; =:\; \ddot v\; , \\ 
\ea
\]
and the operation of $\sa$ by
\[
\dot\sa \= 
\etwasenger
\left[
\ba{ccc|ccc|ccc}
\scm 1 &\scm 0 &\scm 0 &\scm 0 &\scm 0 &\scm 0 &\scm 0 &\scm 0 &\scm 0\\
\scm-\frac{24876}{7217} &\scm\frac{25976}{1031} &\scm -\frac{31737}{7217} &\scm\frac{55068}{7217} &\scm -\frac{25392}{1031} &\scm\frac{5540}{1031} &\scm -\frac{16308}{7217} &\scm \frac{36132}{7217} &\scm -\frac{1176}{1031}\\
\scm -\frac{62766}{1031} &\scm\frac{413856}{1031} &\scm -\frac{149176}{1031} &\scm\frac{117030}{1031} &\scm -\frac{647484}{1031} &\scm \frac{228306}{1031} &\scm -\frac{34326}{1031} &\scm\frac{147615}{1031} &\scm -\frac{51159}{1031}\\\hline\rule{0mm}{2mm}
\scm 0 &\scm 0 &\scm 0 &\scm 4 &\scm 0 &\scm 0 &\scm -1 &\scm 0 &\scm 0\\
\scm -\frac{35412}{7217} &\scm\frac{29460}{1031} &\scm -\frac{75624}{7217} &\scm \frac{164160}{7217} &\scm -\frac{57244}{1031} &\scm \frac{17712}{1031} &\scm -\frac{51540}{7217} &\scm \frac{97432}{7217} &\scm -\frac{4033}{1031}\\
\scm -\frac{845298}{7217} &\scm \frac{815757}{1031} &\scm -\frac{3392322}{7217} &\scm \frac{2508327}{7217} &\scm -\frac{2001978}{1031} &\scm\frac{810182}{1031} &\scm -\frac{796137}{7217} &\scm \frac{3442200}{7217} &\scm -\frac{184868}{1031}\\\hline\rule{0mm}{2mm}
\scm -\frac{210}{1031} &\scm -\frac{2776}{1031} &\scm -\frac{206}{1031} &\scm \frac{13095}{1031} &\scm \frac{563}{1031} &\scm \frac{679}{1031} &\scm -\frac{2419}{1031} &\scm -\frac{347}{1031} &\scm -\frac{162}{1031}\\
\scm -\frac{9945}{1031} &\scm \frac{52791}{1031} &\scm -\frac{35236}{1031} &\scm \frac{72018}{1031} &\scm -\frac{171069}{1031} &\scm \frac{62570}{1031} &\scm -\frac{21546}{1031} &\scm \frac{43586}{1031} &\scm -\frac{14447}{1031}\\
\scm -\frac{2471781}{7217} &\scm \frac{2428305}{1031} &\scm -\frac{12112977}{7217} &\scm \frac{8656818}{7217} &\scm -\frac{7025124}{1031} &\scm \frac{2951589}{1031} &\scm -\frac{2779557}{7217} &\scm \frac{12296484}{7217} &\scm -\frac{676060}{1031}\\
\ea    
\right]\; .
\weiter
\] 
We observe that we may replace $\dot u$, $\dot v$ by $\ddot u$, $\ddot v$ resp.\ by $\dddot u$, $\dddot v$ to obtain
\[
\sLD\= \Lambda((\ddot u,\ddot v), (2,2), (2,7))_{\dot u\Lambda} 
   \= \{f\in\Lambda\; :\; D_{\ddot u}^i\circ D_{\ddot v}^j(f) \con_{\dot u^{2i + 7j}\Lambda} 0\mb{\rm\ for all } i,j\in [0,2] \} 
\]
and 
\[
\Xi\=\sLDE 
\= \{ f\in\sLD\; :\; D_{\;\dddot u}^i\circ E_{\;\dddot u\, ,\, \dddot v}^j(f) \con_{\dot u^{2i + 9j}\Lambda} 0\mb{\rm\ for all } i,j\in [0,2]\}\; .
\]

\subsubsection{The case $p = 5$}

Suppose $p = 5$, i.e.\ $S = \Z_{(5)}$, 
$u = \pi_3 = (\zeta_{125} - 1)(\zeta_{125}^{-1} - 1)(\zeta_{125}^{57} - 1)(\zeta_{125}^{-57} - 1)$,
$\zeta_{125}^\sa = \zeta_{125}^6$, \linebreak 
$v = u\; u^{\sa} u^{\sa^2} u^{\sa^3} u^{\sa^4}$ and $U = \Z_{(5)}[\pi_3]$. With respect to the basis 
\[
\barcl
( & u^0 v^0, u^0 v^1, u^0 v^2, u^0 v^3, u^0 v^4, & \\
  & u^1 v^0, u^1 v^1, u^1 v^2, u^1 v^3, u^1 v^4, & \\
  & u^2 v^0, u^2 v^1, u^2 v^2, u^2 v^3, u^2 v^4, & \\
  & u^3 v^0, u^3 v^1, u^3 v^2, u^3 v^3, u^3 v^4, & \\
  & u^4 v^0, u^4 v^1, u^4 v^2, u^4 v^3, u^4 v^4, & )\; , \\
\ea
\]
the matrix describing the multiplication by $u$ on $U$ reduces to
\[
\barcl
\dot u & \con_{\dot u^{31}\Lambda} &
{\enger
\left[
\ba{rrrrr|rrrrr|rrrrr|rrrrr|rrrrr}
 \scm 0 &\scm 0 &\scm 0 &\scm 0 &\scm 0 
&\scm\;\;\; 1 &\scm\;\;\; 0 &\scm\;\;\; 0 &\scm\;\;\; 0 &\scm\;\;\; 0 
&\scm 0 &\scm 0 &\scm 0 &\scm 0 &\scm 0 
&\scm 0 &\scm 0 &\scm 0 &\scm 0 &\scm 0 
&\scm 0 &\scm 0 &\scm 0 &\scm 0 &\scm 0 \\ 
 \scm 0 &\scm 0 &\scm 0 &\scm 0 &\scm 0 
&\scm 0 &\scm 1 &\scm 0 &\scm 0 &\scm 0 
&\scm 0 &\scm 0 &\scm 0 &\scm 0 &\scm 0 
&\scm 0 &\scm 0 &\scm 0 &\scm 0 &\scm 0 
&\scm 0 &\scm 0 &\scm 0 &\scm 0 &\scm 0 \\ 
 \scm 0 &\scm 0 &\scm 0 &\scm 0 &\scm 0 
&\scm 0 &\scm 0 &\scm 1 &\scm 0 &\scm 0 
&\scm 0 &\scm 0 &\scm 0 &\scm 0 &\scm 0 
&\scm 0 &\scm 0 &\scm 0 &\scm 0 &\scm 0 
&\scm 0 &\scm 0 &\scm 0 &\scm 0 &\scm 0 \\ 
 \scm 0 &\scm 0 &\scm 0 &\scm 0 &\scm 0 
&\scm 0 &\scm 0 &\scm 0 &\scm 1 &\scm 0 
&\scm 0 &\scm 0 &\scm 0 &\scm 0 &\scm 0 
&\scm 0 &\scm 0 &\scm 0 &\scm 0 &\scm 0 
&\scm 0 &\scm 0 &\scm 0 &\scm 0 &\scm 0 \\ 
 \scm 0 &\scm 0 &\scm 0 &\scm 0 &\scm 0 
&\scm 0 &\scm 0 &\scm 0 &\scm 0 &\scm 1 
&\scm 0 &\scm 0 &\scm 0 &\scm 0 &\scm 0 
&\scm 0 &\scm 0 &\scm 0 &\scm 0 &\scm 0 
&\scm 0 &\scm 0 &\scm 0 &\scm 0 &\scm 0 \\\hline\rule{0mm}{2mm} 
 \scm 0 &\scm 0 &\scm 0 &\scm 0 &\scm 0 
&\scm 0 &\scm 0 &\scm 0 &\scm 0 &\scm 0 
&\scm 1 &\scm 0 &\scm 0 &\scm 0 &\scm 0 
&\scm 0 &\scm 0 &\scm 0 &\scm 0 &\scm 0 
&\scm 0 &\scm 0 &\scm 0 &\scm 0 &\scm 0 \\ 
 \scm 0 &\scm 0 &\scm 0 &\scm 0 &\scm 0 
&\scm 0 &\scm 0 &\scm 0 &\scm 0 &\scm 0 
&\scm 0 &\scm 1 &\scm 0 &\scm 0 &\scm 0 
&\scm 0 &\scm 0 &\scm 0 &\scm 0 &\scm 0 
&\scm 0 &\scm 0 &\scm 0 &\scm 0 &\scm 0 \\ 
 \scm 0 &\scm 0 &\scm 0 &\scm 0 &\scm 0 
&\scm 0 &\scm 0 &\scm 0 &\scm 0 &\scm 0 
&\scm 0 &\scm 0 &\scm 1 &\scm 0 &\scm 0 
&\scm 0 &\scm 0 &\scm 0 &\scm 0 &\scm 0 
&\scm 0 &\scm 0 &\scm 0 &\scm 0 &\scm 0 \\ 
 \scm 0 &\scm 0 &\scm 0 &\scm 0 &\scm 0 
&\scm 0 &\scm 0 &\scm 0 &\scm 0 &\scm 0 
&\scm 0 &\scm 0 &\scm 0 &\scm 1 &\scm 0 
&\scm 0 &\scm 0 &\scm 0 &\scm 0 &\scm 0 
&\scm 0 &\scm 0 &\scm 0 &\scm 0 &\scm 0 \\ 
 \scm 0 &\scm 0 &\scm 0 &\scm 0 &\scm 0 
&\scm 0 &\scm 0 &\scm 0 &\scm 0 &\scm 0 
&\scm 0 &\scm 0 &\scm 0 &\scm 0 &\scm 1 
&\scm 0 &\scm 0 &\scm 0 &\scm 0 &\scm 0 
&\scm 0 &\scm 0 &\scm 0 &\scm 0 &\scm 0 \\\hline\rule{0mm}{2mm} 
 \scm 0 &\scm 0 &\scm 0 &\scm 0 &\scm 0 
&\scm 0 &\scm 0 &\scm 0 &\scm 0 &\scm 0 
&\scm 0 &\scm 0 &\scm 0 &\scm 0 &\scm 0 
&\scm 1 &\scm 0 &\scm 0 &\scm 0 &\scm 0 
&\scm 0 &\scm 0 &\scm 0 &\scm 0 &\scm 0 \\ 
 \scm 0 &\scm 0 &\scm 0 &\scm 0 &\scm 0 
&\scm 0 &\scm 0 &\scm 0 &\scm 0 &\scm 0 
&\scm 0 &\scm 0 &\scm 0 &\scm 0 &\scm 0 
&\scm 0 &\scm 1 &\scm 0 &\scm 0 &\scm 0 
&\scm 0 &\scm 0 &\scm 0 &\scm 0 &\scm 0 \\ 
 \scm 0 &\scm 0 &\scm 0 &\scm 0 &\scm 0 
&\scm 0 &\scm 0 &\scm 0 &\scm 0 &\scm 0 
&\scm 0 &\scm 0 &\scm 0 &\scm 0 &\scm 0 
&\scm 0 &\scm 0 &\scm 1 &\scm 0 &\scm 0 
&\scm 0 &\scm 0 &\scm 0 &\scm 0 &\scm 0 \\ 
 \scm 0 &\scm 0 &\scm 0 &\scm 0 &\scm 0 
&\scm 0 &\scm 0 &\scm 0 &\scm 0 &\scm 0 
&\scm 0 &\scm 0 &\scm 0 &\scm 0 &\scm 0 
&\scm 0 &\scm 0 &\scm 0 &\scm 1 &\scm 0 
&\scm 0 &\scm 0 &\scm 0 &\scm 0 &\scm 0 \\ 
 \scm 0 &\scm 0 &\scm 0 &\scm 0 &\scm 0 
&\scm 0 &\scm 0 &\scm 0 &\scm 0 &\scm 0 
&\scm 0 &\scm 0 &\scm 0 &\scm 0 &\scm 0 
&\scm 0 &\scm 0 &\scm 0 &\scm 0 &\scm 1 
&\scm 0 &\scm 0 &\scm 0 &\scm 0 &\scm 0 \\\hline\rule{0mm}{2mm} 
 \scm 0 &\scm 0 &\scm 0 &\scm 0 &\scm 0 
&\scm 0 &\scm 0 &\scm 0 &\scm 0 &\scm 0 
&\scm 0 &\scm 0 &\scm 0 &\scm 0 &\scm 0 
&\scm 0 &\scm 0 &\scm 0 &\scm 0 &\scm 0 
&\scm 1 &\scm 0 &\scm 0 &\scm 0 &\scm 0 \\ 
 \scm 0 &\scm 0 &\scm 0 &\scm 0 &\scm 0 
&\scm 0 &\scm 0 &\scm 0 &\scm 0 &\scm 0 
&\scm 0 &\scm 0 &\scm 0 &\scm 0 &\scm 0 
&\scm 0 &\scm 0 &\scm 0 &\scm 0 &\scm 0 
&\scm 0 &\scm 1 &\scm 0 &\scm 0 &\scm 0 \\ 
 \scm 0 &\scm 0 &\scm 0 &\scm 0 &\scm 0 
&\scm 0 &\scm 0 &\scm 0 &\scm 0 &\scm 0 
&\scm 0 &\scm 0 &\scm 0 &\scm 0 &\scm 0 
&\scm 0 &\scm 0 &\scm 0 &\scm 0 &\scm 0 
&\scm 0 &\scm 0 &\scm 1 &\scm 0 &\scm 0 \\ 
 \scm 0 &\scm 0 &\scm 0 &\scm 0 &\scm 0 
&\scm 0 &\scm 0 &\scm 0 &\scm 0 &\scm 0 
&\scm 0 &\scm 0 &\scm 0 &\scm 0 &\scm 0 
&\scm 0 &\scm 0 &\scm 0 &\scm 0 &\scm 0 
&\scm 0 &\scm 0 &\scm 0 &\scm 1 &\scm 0 \\ 
 \scm 0 &\scm 0 &\scm 0 &\scm 0 &\scm 0 
&\scm 0 &\scm 0 &\scm 0 &\scm 0 &\scm 0 
&\scm 0 &\scm 0 &\scm 0 &\scm 0 &\scm 0 
&\scm 0 &\scm 0 &\scm 0 &\scm 0 &\scm 0 
&\scm 0 &\scm 0 &\scm 0 &\scm 0 &\scm 1 \\\hline\rule{0mm}{2mm} 
 \scm  -5 &\scm  -9 &\scm   0 &\scm   3 &\scm   3 
&\scm  10 &\scm   0 &\scm   0 &\scm   3 &\scm   1 
&\scm -10 &\scm  10 &\scm   0 &\scm   1 &\scm   2 
&\scm   0 &\scm -10 &\scm  -1 &\scm   1 &\scm   1 
&\scm  -5 &\scm  -6 &\scm   1 &\scm   0 &\scm   2 \\ 
 \scm  15 &\scm  -5 &\scm  -9 &\scm   0 &\scm   3 
&\scm   5 &\scm  10 &\scm   0 &\scm   0 &\scm   3 
&\scm  10 &\scm -10 &\scm  10 &\scm   0 &\scm   1 
&\scm   5 &\scm   0 &\scm  10 &\scm  -1 &\scm   1 
&\scm  10 &\scm  -5 &\scm  -6 &\scm   1 &\scm   0 \\ 
 \scm  15 &\scm  15 &\scm  -5 &\scm  -9 &\scm   0 
&\scm  15 &\scm   5 &\scm  10 &\scm   0 &\scm   0 
&\scm   5 &\scm  10 &\scm -10 &\scm  10 &\scm   0 
&\scm   5 &\scm   5 &\scm   0 &\scm  10 &\scm  -1 
&\scm   0 &\scm  10 &\scm  -5 &\scm  -6 &\scm   1 \\ 
 \scm   0 &\scm  15 &\scm  15 &\scm  -5 &\scm  -9 
&\scm   0 &\scm  15 &\scm   5 &\scm  10 &\scm   0 
&\scm   0 &\scm   5 &\scm  10 &\scm -10 &\scm  10 
&\scm  -5 &\scm   5 &\scm   5 &\scm   0 &\scm  10 
&\scm   5 &\scm   0 &\scm  10 &\scm  -5 &\scm  -6 \\ 
 \scm -45 &\scm   0 &\scm  15 &\scm  15 &\scm  -5 
&\scm   0 &\scm   0 &\scm  15 &\scm   5 &\scm  10 
&\scm  50 &\scm   0 &\scm   5 &\scm  10 &\scm -10 
&\scm  50 &\scm  -5 &\scm   5 &\scm   5 &\scm   0 
&\scm  -5 &\scm   5 &\scm   0 &\scm  10 &\scm  -5 \\
\ea
\right]}
\;\; =: \;\; \dddot{u} \vs\\
& \con_{\dot u^5\Lambda} &
{\enger
\left[
\ba{rrrrr|rrrrr|rrrrr|rrrrr|rrrrr}
 \scm 0 &\scm 0 &\scm 0 &\scm 0 &\scm 0 
&\scm 1 &\scm 0 &\scm 0 &\scm 0 &\scm 0 
&\scm 0 &\scm 0 &\scm 0 &\scm 0 &\scm 0 
&\scm 0 &\scm 0 &\scm 0 &\scm 0 &\scm 0 
&\scm 0 &\scm 0 &\scm 0 &\scm 0 &\scm 0 \\ 
 \scm 0 &\scm 0 &\scm 0 &\scm 0 &\scm 0 
&\scm 0 &\scm 1 &\scm 0 &\scm 0 &\scm 0 
&\scm 0 &\scm 0 &\scm 0 &\scm 0 &\scm 0 
&\scm 0 &\scm 0 &\scm 0 &\scm 0 &\scm 0 
&\scm 0 &\scm 0 &\scm 0 &\scm 0 &\scm 0 \\ 
 \scm 0 &\scm 0 &\scm 0 &\scm 0 &\scm 0 
&\scm 0 &\scm 0 &\scm 1 &\scm 0 &\scm 0 
&\scm 0 &\scm 0 &\scm 0 &\scm 0 &\scm 0 
&\scm 0 &\scm 0 &\scm 0 &\scm 0 &\scm 0 
&\scm 0 &\scm 0 &\scm 0 &\scm 0 &\scm 0 \\ 
 \scm 0 &\scm 0 &\scm 0 &\scm 0 &\scm 0 
&\scm 0 &\scm 0 &\scm 0 &\scm 1 &\scm 0 
&\scm 0 &\scm 0 &\scm 0 &\scm 0 &\scm 0 
&\scm 0 &\scm 0 &\scm 0 &\scm 0 &\scm 0 
&\scm 0 &\scm 0 &\scm 0 &\scm 0 &\scm 0 \\ 
 \scm 0 &\scm 0 &\scm 0 &\scm 0 &\scm 0 
&\scm 0 &\scm 0 &\scm 0 &\scm 0 &\scm 1 
&\scm 0 &\scm 0 &\scm 0 &\scm 0 &\scm 0 
&\scm 0 &\scm 0 &\scm 0 &\scm 0 &\scm 0 
&\scm 0 &\scm 0 &\scm 0 &\scm 0 &\scm 0 \\\hline\rule{0mm}{2mm} 
 \scm 0 &\scm 0 &\scm 0 &\scm 0 &\scm 0 
&\scm 0 &\scm 0 &\scm 0 &\scm 0 &\scm 0 
&\scm 1 &\scm 0 &\scm 0 &\scm 0 &\scm 0 
&\scm 0 &\scm 0 &\scm 0 &\scm 0 &\scm 0 
&\scm 0 &\scm 0 &\scm 0 &\scm 0 &\scm 0 \\ 
 \scm 0 &\scm 0 &\scm 0 &\scm 0 &\scm 0 
&\scm 0 &\scm 0 &\scm 0 &\scm 0 &\scm 0 
&\scm 0 &\scm 1 &\scm 0 &\scm 0 &\scm 0 
&\scm 0 &\scm 0 &\scm 0 &\scm 0 &\scm 0 
&\scm 0 &\scm 0 &\scm 0 &\scm 0 &\scm 0 \\ 
 \scm 0 &\scm 0 &\scm 0 &\scm 0 &\scm 0 
&\scm 0 &\scm 0 &\scm 0 &\scm 0 &\scm 0 
&\scm 0 &\scm 0 &\scm 1 &\scm 0 &\scm 0 
&\scm 0 &\scm 0 &\scm 0 &\scm 0 &\scm 0 
&\scm 0 &\scm 0 &\scm 0 &\scm 0 &\scm 0 \\ 
 \scm 0 &\scm 0 &\scm 0 &\scm 0 &\scm 0 
&\scm 0 &\scm 0 &\scm 0 &\scm 0 &\scm 0 
&\scm 0 &\scm 0 &\scm 0 &\scm 1 &\scm 0 
&\scm 0 &\scm 0 &\scm 0 &\scm 0 &\scm 0 
&\scm 0 &\scm 0 &\scm 0 &\scm 0 &\scm 0 \\ 
 \scm 0 &\scm 0 &\scm 0 &\scm 0 &\scm 0 
&\scm 0 &\scm 0 &\scm 0 &\scm 0 &\scm 0 
&\scm 0 &\scm 0 &\scm 0 &\scm 0 &\scm 1 
&\scm 0 &\scm 0 &\scm 0 &\scm 0 &\scm 0 
&\scm 0 &\scm 0 &\scm 0 &\scm 0 &\scm 0 \\\hline\rule{0mm}{2mm} 
 \scm 0 &\scm 0 &\scm 0 &\scm 0 &\scm 0 
&\scm 0 &\scm 0 &\scm 0 &\scm 0 &\scm 0 
&\scm 0 &\scm 0 &\scm 0 &\scm 0 &\scm 0 
&\scm 1 &\scm 0 &\scm 0 &\scm 0 &\scm 0 
&\scm 0 &\scm 0 &\scm 0 &\scm 0 &\scm 0 \\ 
 \scm 0 &\scm 0 &\scm 0 &\scm 0 &\scm 0 
&\scm 0 &\scm 0 &\scm 0 &\scm 0 &\scm 0 
&\scm 0 &\scm 0 &\scm 0 &\scm 0 &\scm 0 
&\scm 0 &\scm 1 &\scm 0 &\scm 0 &\scm 0 
&\scm 0 &\scm 0 &\scm 0 &\scm 0 &\scm 0 \\ 
 \scm 0 &\scm 0 &\scm 0 &\scm 0 &\scm 0 
&\scm 0 &\scm 0 &\scm 0 &\scm 0 &\scm 0 
&\scm 0 &\scm 0 &\scm 0 &\scm 0 &\scm 0 
&\scm 0 &\scm 0 &\scm 1 &\scm 0 &\scm 0 
&\scm 0 &\scm 0 &\scm 0 &\scm 0 &\scm 0 \\ 
 \scm 0 &\scm 0 &\scm 0 &\scm 0 &\scm 0 
&\scm 0 &\scm 0 &\scm 0 &\scm 0 &\scm 0 
&\scm 0 &\scm 0 &\scm 0 &\scm 0 &\scm 0 
&\scm 0 &\scm 0 &\scm 0 &\scm 1 &\scm 0 
&\scm 0 &\scm 0 &\scm 0 &\scm 0 &\scm 0 \\ 
 \scm 0 &\scm 0 &\scm 0 &\scm 0 &\scm 0 
&\scm 0 &\scm 0 &\scm 0 &\scm 0 &\scm 0 
&\scm 0 &\scm 0 &\scm 0 &\scm 0 &\scm 0 
&\scm 0 &\scm 0 &\scm 0 &\scm 0 &\scm 1 
&\scm 0 &\scm 0 &\scm 0 &\scm 0 &\scm 0 \\\hline\rule{0mm}{2mm} 
 \scm 0 &\scm 0 &\scm 0 &\scm 0 &\scm 0 
&\scm 0 &\scm 0 &\scm 0 &\scm 0 &\scm 0 
&\scm 0 &\scm 0 &\scm 0 &\scm 0 &\scm 0 
&\scm 0 &\scm 0 &\scm 0 &\scm 0 &\scm 0 
&\scm 1 &\scm 0 &\scm 0 &\scm 0 &\scm 0 \\ 
 \scm 0 &\scm 0 &\scm 0 &\scm 0 &\scm 0 
&\scm 0 &\scm 0 &\scm 0 &\scm 0 &\scm 0 
&\scm 0 &\scm 0 &\scm 0 &\scm 0 &\scm 0 
&\scm 0 &\scm 0 &\scm 0 &\scm 0 &\scm 0 
&\scm 0 &\scm 1 &\scm 0 &\scm 0 &\scm 0 \\ 
 \scm 0 &\scm 0 &\scm 0 &\scm 0 &\scm 0 
&\scm 0 &\scm 0 &\scm 0 &\scm 0 &\scm 0 
&\scm 0 &\scm 0 &\scm 0 &\scm 0 &\scm 0 
&\scm 0 &\scm 0 &\scm 0 &\scm 0 &\scm 0 
&\scm 0 &\scm 0 &\scm 1 &\scm 0 &\scm 0 \\ 
 \scm 0 &\scm 0 &\scm 0 &\scm 0 &\scm 0 
&\scm 0 &\scm 0 &\scm 0 &\scm 0 &\scm 0 
&\scm 0 &\scm 0 &\scm 0 &\scm 0 &\scm 0 
&\scm 0 &\scm 0 &\scm 0 &\scm 0 &\scm 0 
&\scm 0 &\scm 0 &\scm 0 &\scm 1 &\scm 0 \\ 
 \scm 0 &\scm 0 &\scm 0 &\scm 0 &\scm 0 
&\scm 0 &\scm 0 &\scm 0 &\scm 0 &\scm 0 
&\scm 0 &\scm 0 &\scm 0 &\scm 0 &\scm 0 
&\scm 0 &\scm 0 &\scm 0 &\scm 0 &\scm 0 
&\scm 0 &\scm 0 &\scm 0 &\scm 0 &\scm 1 \\\hline\rule{0mm}{2mm} 
 \scm 0 &\scm 1 &\scm 0 &\scm 0 &\scm 0
&\scm 0 &\scm 0 &\scm 0 &\scm 0 &\scm 0
&\scm 0 &\scm 0 &\scm 0 &\scm 0 &\scm 0
&\scm 0 &\scm 0 &\scm 0 &\scm 0 &\scm 0
&\scm 0 &\scm 0 &\scm 0 &\scm 0 &\scm 0 \\
 \scm 0 &\scm 0 &\scm 1 &\scm 0 &\scm 0
&\scm 0 &\scm 0 &\scm 0 &\scm 0 &\scm 0
&\scm 0 &\scm 0 &\scm 0 &\scm 0 &\scm 0
&\scm 0 &\scm 0 &\scm 0 &\scm 0 &\scm 0
&\scm 0 &\scm 0 &\scm 0 &\scm 0 &\scm 0 \\
 \scm 0 &\scm 0 &\scm 0 &\scm 1 &\scm 0
&\scm 0 &\scm 0 &\scm 0 &\scm 0 &\scm 0
&\scm 0 &\scm 0 &\scm 0 &\scm 0 &\scm 0
&\scm 0 &\scm 0 &\scm 0 &\scm 0 &\scm 0
&\scm 0 &\scm 0 &\scm 0 &\scm 0 &\scm 0 \\
 \scm 0 &\scm 0 &\scm 0 &\scm 0 &\scm 1
&\scm 0 &\scm 0 &\scm 0 &\scm 0 &\scm 0
&\scm 0 &\scm 0 &\scm 0 &\scm 0 &\scm 0
&\scm 0 &\scm 0 &\scm 0 &\scm 0 &\scm 0
&\scm 0 &\scm 0 &\scm 0 &\scm 0 &\scm 0 \\
 \scm 5 &\scm 0 &\scm 0 &\scm 0 &\scm 0
&\scm 0 &\scm 0 &\scm 0 &\scm 0 &\scm 0
&\scm 0 &\scm 0 &\scm 0 &\scm 0 &\scm 0
&\scm 0 &\scm 0 &\scm 0 &\scm 0 &\scm 0
&\scm 0 &\scm 0 &\scm 0 &\scm 0 &\scm 0 \\
\ea
\right]}
\;\; =: \;\; \ddot{u}\; , \\
\ea
\]
and the matrix describing the multiplication by the Sen-element $v$ reduces to
\[
\dot v  \;\;\con_{\dot u^{25 + 5\cdot (5+1)/2}\Lambda}\;\;
{\enger
\left[
\ba{rrrrr|rrrrr|rrrrr|rrrrr|rrrrr}
 \scm 0 &\scm 1 &\scm 0 &\scm 0 &\scm 0 
&\scm 0 &\scm 0 &\scm 0 &\scm 0 &\scm 0 
&\scm 0 &\scm 0 &\scm 0 &\scm 0 &\scm 0 
&\scm 0 &\scm 0 &\scm 0 &\scm 0 &\scm 0 
&\scm 0 &\scm 0 &\scm 0 &\scm 0 &\scm 0 \\ 
 \scm 0 &\scm 0 &\scm 1 &\scm 0 &\scm 0 
&\scm 0 &\scm 0 &\scm 0 &\scm 0 &\scm 0 
&\scm 0 &\scm 0 &\scm 0 &\scm 0 &\scm 0 
&\scm 0 &\scm 0 &\scm 0 &\scm 0 &\scm 0 
&\scm 0 &\scm 0 &\scm 0 &\scm 0 &\scm 0 \\ 
 \scm 0 &\scm 0 &\scm 0 &\scm 1 &\scm 0 
&\scm 0 &\scm 0 &\scm 0 &\scm 0 &\scm 0 
&\scm 0 &\scm 0 &\scm 0 &\scm 0 &\scm 0 
&\scm 0 &\scm 0 &\scm 0 &\scm 0 &\scm 0 
&\scm 0 &\scm 0 &\scm 0 &\scm 0 &\scm 0 \\ 
 \scm 0 &\scm 0 &\scm 0 &\scm 0 &\scm 1 
&\scm 0 &\scm 0 &\scm 0 &\scm 0 &\scm 0 
&\scm 0 &\scm 0 &\scm 0 &\scm 0 &\scm 0 
&\scm 0 &\scm 0 &\scm 0 &\scm 0 &\scm 0 
&\scm 0 &\scm 0 &\scm 0 &\scm 0 &\scm 0 \\ 
 \scm  5 &\scm\;\;  75 &\scm 0 &\scm\;\;\;\; 0 &\scm\;\;\;\; 0 
&\scm  0 &\scm\;\;\; 0 &\scm 0 &\scm\;\;\;\; 0 &\scm\;\;\;\; 0 
&\scm  0 &\scm\;\;\; 0 &\scm 0 &\scm\;\;\;\; 0 &\scm\;\;\;\; 0 
&\scm  0 &\scm -25 &\scm\;\;\; 0 &\scm\;\;\;\; 0 &\scm\;\;\;\; 0 
&\scm 25 &\scm  25 &\scm\;\;\;\; 0 &\scm\;\;\;\; 0 &\scm\;\;\;\; 0 \\\hline\rule{0mm}{2mm} 
 \scm 0 &\scm 0 &\scm 0 &\scm 0 &\scm 0 
&\scm 0 &\scm 1 &\scm 0 &\scm 0 &\scm 0 
&\scm 0 &\scm 0 &\scm 0 &\scm 0 &\scm 0 
&\scm 0 &\scm 0 &\scm 0 &\scm 0 &\scm 0 
&\scm 0 &\scm 0 &\scm 0 &\scm 0 &\scm 0 \\ 
 \scm 0 &\scm 0 &\scm 0 &\scm 0 &\scm 0 
&\scm 0 &\scm 0 &\scm 1 &\scm 0 &\scm 0 
&\scm 0 &\scm 0 &\scm 0 &\scm 0 &\scm 0 
&\scm 0 &\scm 0 &\scm 0 &\scm 0 &\scm 0 
&\scm 0 &\scm 0 &\scm 0 &\scm 0 &\scm 0 \\ 
 \scm 0 &\scm 0 &\scm 0 &\scm 0 &\scm 0 
&\scm 0 &\scm 0 &\scm 0 &\scm 1 &\scm 0 
&\scm 0 &\scm 0 &\scm 0 &\scm 0 &\scm 0 
&\scm 0 &\scm 0 &\scm 0 &\scm 0 &\scm 0 
&\scm 0 &\scm 0 &\scm 0 &\scm 0 &\scm 0 \\ 
 \scm 0 &\scm 0 &\scm 0 &\scm 0 &\scm 0 
&\scm 0 &\scm 0 &\scm 0 &\scm 0 &\scm 1 
&\scm 0 &\scm 0 &\scm 0 &\scm 0 &\scm 0 
&\scm 0 &\scm 0 &\scm 0 &\scm 0 &\scm 0 
&\scm 0 &\scm 0 &\scm 0 &\scm 0 &\scm 0 \\ 
 \scm 0 &\scm  25 &\scm 25 &\scm 0 &\scm 0 
&\scm 5 &\scm  75 &\scm  0 &\scm 0 &\scm 0 
&\scm 0 &\scm   0 &\scm  0 &\scm 0 &\scm 0 
&\scm 0 &\scm   0 &\scm  0 &\scm 0 &\scm 0 
&\scm 0 &\scm -50 &\scm  0 &\scm 0 &\scm 0 \\\hline\rule{0mm}{2mm}
 \scm 0 &\scm 0 &\scm 0 &\scm 0 &\scm 0 
&\scm 0 &\scm 0 &\scm 0 &\scm 0 &\scm 0 
&\scm 0 &\scm 1 &\scm 0 &\scm 0 &\scm 0 
&\scm 0 &\scm 0 &\scm 0 &\scm 0 &\scm 0 
&\scm 0 &\scm 0 &\scm 0 &\scm 0 &\scm 0 \\ 
 \scm 0 &\scm 0 &\scm 0 &\scm 0 &\scm 0 
&\scm 0 &\scm 0 &\scm 0 &\scm 0 &\scm 0 
&\scm 0 &\scm 0 &\scm 1 &\scm 0 &\scm 0 
&\scm 0 &\scm 0 &\scm 0 &\scm 0 &\scm 0 
&\scm 0 &\scm 0 &\scm 0 &\scm 0 &\scm 0 \\ 
 \scm 0 &\scm 0 &\scm 0 &\scm 0 &\scm 0 
&\scm 0 &\scm 0 &\scm 0 &\scm 0 &\scm 0 
&\scm 0 &\scm 0 &\scm 0 &\scm 1 &\scm 0 
&\scm 0 &\scm 0 &\scm 0 &\scm 0 &\scm 0 
&\scm 0 &\scm 0 &\scm 0 &\scm 0 &\scm 0 \\ 
 \scm 0 &\scm 0 &\scm 0 &\scm 0 &\scm 0 
&\scm 0 &\scm 0 &\scm 0 &\scm 0 &\scm 0 
&\scm 0 &\scm 0 &\scm 0 &\scm 0 &\scm 1 
&\scm 0 &\scm 0 &\scm 0 &\scm 0 &\scm 0 
&\scm 0 &\scm 0 &\scm 0 &\scm 0 &\scm 0 \\ 
 \scm 0 &\scm  0 &\scm -50 &\scm 0 &\scm 0 
&\scm 0 &\scm 25 &\scm  25 &\scm 0 &\scm 0 
&\scm 5 &\scm 75 &\scm   0 &\scm 0 &\scm 0 
&\scm 0 &\scm  0 &\scm   0 &\scm 0 &\scm 0 
&\scm 0 &\scm  0 &\scm   0 &\scm 0 &\scm 0 \\\hline\rule{0mm}{2mm} 
 \scm 0 &\scm 0 &\scm 0 &\scm 0 &\scm 0 
&\scm 0 &\scm 0 &\scm 0 &\scm 0 &\scm 0 
&\scm 0 &\scm 0 &\scm 0 &\scm 0 &\scm 0 
&\scm 0 &\scm 1 &\scm 0 &\scm 0 &\scm 0 
&\scm 0 &\scm 0 &\scm 0 &\scm 0 &\scm 0 \\ 
 \scm 0 &\scm 0 &\scm 0 &\scm 0 &\scm 0 
&\scm 0 &\scm 0 &\scm 0 &\scm 0 &\scm 0 
&\scm 0 &\scm 0 &\scm 0 &\scm 0 &\scm 0 
&\scm 0 &\scm 0 &\scm 1 &\scm 0 &\scm 0 
&\scm 0 &\scm 0 &\scm 0 &\scm 0 &\scm 0 \\ 
 \scm 0 &\scm 0 &\scm 0 &\scm 0 &\scm 0 
&\scm 0 &\scm 0 &\scm 0 &\scm 0 &\scm 0 
&\scm 0 &\scm 0 &\scm 0 &\scm 0 &\scm 0 
&\scm 0 &\scm 0 &\scm 0 &\scm 1 &\scm 0 
&\scm 0 &\scm 0 &\scm 0 &\scm 0 &\scm 0 \\ 
 \scm 0 &\scm 0 &\scm 0 &\scm 0 &\scm 0 
&\scm 0 &\scm 0 &\scm 0 &\scm 0 &\scm 0 
&\scm 0 &\scm 0 &\scm 0 &\scm 0 &\scm 0 
&\scm 0 &\scm 0 &\scm 0 &\scm 0 &\scm 1 
&\scm 0 &\scm 0 &\scm 0 &\scm 0 &\scm 0 \\ 
 \scm 0 &\scm  0 &\scm   0 &\scm 0 &\scm 0 
&\scm 0 &\scm  0 &\scm -50 &\scm 0 &\scm 0 
&\scm 0 &\scm 25 &\scm  25 &\scm 0 &\scm 0 
&\scm 5 &\scm 75 &\scm   0 &\scm 0 &\scm 0 
&\scm 0 &\scm  0 &\scm   0 &\scm 0 &\scm 0 \\\hline\rule{0mm}{2mm} 
 \scm 0 &\scm 0 &\scm 0 &\scm 0 &\scm 0 
&\scm 0 &\scm 0 &\scm 0 &\scm 0 &\scm 0 
&\scm 0 &\scm 0 &\scm 0 &\scm 0 &\scm 0 
&\scm 0 &\scm 0 &\scm 0 &\scm 0 &\scm 0 
&\scm 0 &\scm 1 &\scm 0 &\scm 0 &\scm 0 \\ 
 \scm 0 &\scm 0 &\scm 0 &\scm 0 &\scm 0 
&\scm 0 &\scm 0 &\scm 0 &\scm 0 &\scm 0 
&\scm 0 &\scm 0 &\scm 0 &\scm 0 &\scm 0 
&\scm 0 &\scm 0 &\scm 0 &\scm 0 &\scm 0 
&\scm 0 &\scm 0 &\scm 1 &\scm 0 &\scm 0 \\ 
 \scm 0 &\scm 0 &\scm 0 &\scm 0 &\scm 0 
&\scm 0 &\scm 0 &\scm 0 &\scm 0 &\scm 0 
&\scm 0 &\scm 0 &\scm 0 &\scm 0 &\scm 0 
&\scm 0 &\scm 0 &\scm 0 &\scm 0 &\scm 0 
&\scm 0 &\scm 0 &\scm 0 &\scm 1 &\scm 0 \\ 
 \scm 0 &\scm 0 &\scm 0 &\scm 0 &\scm 0 
&\scm 0 &\scm 0 &\scm 0 &\scm 0 &\scm 0 
&\scm 0 &\scm 0 &\scm 0 &\scm 0 &\scm 0 
&\scm 0 &\scm 0 &\scm 0 &\scm 0 &\scm 0 
&\scm 0 &\scm 0 &\scm 0 &\scm 0 &\scm 1 \\ 
 \scm 0 &\scm  0 &\scm   0 &\scm 0 &\scm 0 
&\scm 0 &\scm  0 &\scm   0 &\scm 0 &\scm 0 
&\scm 0 &\scm  0 &\scm -50 &\scm 0 &\scm 0 
&\scm 0 &\scm 25 &\scm  25 &\scm 0 &\scm 0 
&\scm 5 &\scm 75 &\scm   0 &\scm 0 &\scm 0 \\
\ea
\right]}
\;\; =: \;\; \dddot{v} 
\]
\[
\con_{\dot u^{25+5}\Lambda} \;\;
{\enger
\left[
\ba{rrrrr|rrrrr|rrrrr|rrrrr|rrrrr}
 \scm 0 &\scm 1 &\scm 0 &\scm 0 &\scm 0 
&\scm 0 &\scm 0 &\scm 0 &\scm 0 &\scm 0 
&\scm 0 &\scm 0 &\scm 0 &\scm 0 &\scm 0 
&\scm 0 &\scm 0 &\scm 0 &\scm 0 &\scm 0 
&\scm 0 &\scm 0 &\scm 0 &\scm 0 &\scm 0 \\ 
 \scm 0 &\scm 0 &\scm 1 &\scm 0 &\scm 0 
&\scm 0 &\scm 0 &\scm 0 &\scm 0 &\scm 0 
&\scm 0 &\scm 0 &\scm 0 &\scm 0 &\scm 0 
&\scm 0 &\scm 0 &\scm 0 &\scm 0 &\scm 0 
&\scm 0 &\scm 0 &\scm 0 &\scm 0 &\scm 0 \\ 
 \scm 0 &\scm 0 &\scm 0 &\scm 1 &\scm 0 
&\scm 0 &\scm 0 &\scm 0 &\scm 0 &\scm 0 
&\scm 0 &\scm 0 &\scm 0 &\scm 0 &\scm 0 
&\scm 0 &\scm 0 &\scm 0 &\scm 0 &\scm 0 
&\scm 0 &\scm 0 &\scm 0 &\scm 0 &\scm 0 \\ 
 \scm 0 &\scm 0 &\scm 0 &\scm 0 &\scm 1 
&\scm 0 &\scm 0 &\scm 0 &\scm 0 &\scm 0 
&\scm 0 &\scm 0 &\scm 0 &\scm 0 &\scm 0 
&\scm 0 &\scm 0 &\scm 0 &\scm 0 &\scm 0 
&\scm 0 &\scm 0 &\scm 0 &\scm 0 &\scm 0 \\ 
 \scm  5 &\scm   0 &\scm 0 &\scm 0 &\scm 0 
&\scm  0 &\scm   0 &\scm 0 &\scm 0 &\scm 0 
&\scm  0 &\scm   0 &\scm 0 &\scm 0 &\scm 0 
&\scm  0 &\scm   0 &\scm 0 &\scm 0 &\scm 0 
&\scm  0 &\scm   0 &\scm 0 &\scm 0 &\scm 0 \\\hline\rule{0mm}{2mm} 
 \scm 0 &\scm 0 &\scm 0 &\scm 0 &\scm 0 
&\scm 0 &\scm 1 &\scm 0 &\scm 0 &\scm 0 
&\scm 0 &\scm 0 &\scm 0 &\scm 0 &\scm 0 
&\scm 0 &\scm 0 &\scm 0 &\scm 0 &\scm 0 
&\scm 0 &\scm 0 &\scm 0 &\scm 0 &\scm 0 \\ 
 \scm 0 &\scm 0 &\scm 0 &\scm 0 &\scm 0 
&\scm 0 &\scm 0 &\scm 1 &\scm 0 &\scm 0 
&\scm 0 &\scm 0 &\scm 0 &\scm 0 &\scm 0 
&\scm 0 &\scm 0 &\scm 0 &\scm 0 &\scm 0 
&\scm 0 &\scm 0 &\scm 0 &\scm 0 &\scm 0 \\ 
 \scm 0 &\scm 0 &\scm 0 &\scm 0 &\scm 0 
&\scm 0 &\scm 0 &\scm 0 &\scm 1 &\scm 0 
&\scm 0 &\scm 0 &\scm 0 &\scm 0 &\scm 0 
&\scm 0 &\scm 0 &\scm 0 &\scm 0 &\scm 0 
&\scm 0 &\scm 0 &\scm 0 &\scm 0 &\scm 0 \\ 
 \scm 0 &\scm 0 &\scm 0 &\scm 0 &\scm 0 
&\scm 0 &\scm 0 &\scm 0 &\scm 0 &\scm 1 
&\scm 0 &\scm 0 &\scm 0 &\scm 0 &\scm 0 
&\scm 0 &\scm 0 &\scm 0 &\scm 0 &\scm 0 
&\scm 0 &\scm 0 &\scm 0 &\scm 0 &\scm 0 \\ 
 \scm 0 &\scm   0 &\scm  0 &\scm 0 &\scm 0 
&\scm 5 &\scm   0 &\scm  0 &\scm 0 &\scm 0 
&\scm 0 &\scm   0 &\scm  0 &\scm 0 &\scm 0 
&\scm 0 &\scm   0 &\scm  0 &\scm 0 &\scm 0 
&\scm 0 &\scm   0 &\scm  0 &\scm 0 &\scm 0 \\\hline\rule{0mm}{2mm}
 \scm 0 &\scm 0 &\scm 0 &\scm 0 &\scm 0 
&\scm 0 &\scm 0 &\scm 0 &\scm 0 &\scm 0 
&\scm 0 &\scm 1 &\scm 0 &\scm 0 &\scm 0 
&\scm 0 &\scm 0 &\scm 0 &\scm 0 &\scm 0 
&\scm 0 &\scm 0 &\scm 0 &\scm 0 &\scm 0 \\ 
 \scm 0 &\scm 0 &\scm 0 &\scm 0 &\scm 0 
&\scm 0 &\scm 0 &\scm 0 &\scm 0 &\scm 0 
&\scm 0 &\scm 0 &\scm 1 &\scm 0 &\scm 0 
&\scm 0 &\scm 0 &\scm 0 &\scm 0 &\scm 0 
&\scm 0 &\scm 0 &\scm 0 &\scm 0 &\scm 0 \\ 
 \scm 0 &\scm 0 &\scm 0 &\scm 0 &\scm 0 
&\scm 0 &\scm 0 &\scm 0 &\scm 0 &\scm 0 
&\scm 0 &\scm 0 &\scm 0 &\scm 1 &\scm 0 
&\scm 0 &\scm 0 &\scm 0 &\scm 0 &\scm 0 
&\scm 0 &\scm 0 &\scm 0 &\scm 0 &\scm 0 \\ 
 \scm 0 &\scm 0 &\scm 0 &\scm 0 &\scm 0 
&\scm 0 &\scm 0 &\scm 0 &\scm 0 &\scm 0 
&\scm 0 &\scm 0 &\scm 0 &\scm 0 &\scm 1 
&\scm 0 &\scm 0 &\scm 0 &\scm 0 &\scm 0 
&\scm 0 &\scm 0 &\scm 0 &\scm 0 &\scm 0 \\ 
 \scm 0 &\scm  0 &\scm   0 &\scm 0 &\scm 0 
&\scm 0 &\scm  0 &\scm   0 &\scm 0 &\scm 0 
&\scm 5 &\scm  0 &\scm   0 &\scm 0 &\scm 0 
&\scm 0 &\scm  0 &\scm   0 &\scm 0 &\scm 0 
&\scm 0 &\scm  0 &\scm   0 &\scm 0 &\scm 0 \\\hline\rule{0mm}{2mm} 
 \scm 0 &\scm 0 &\scm 0 &\scm 0 &\scm 0 
&\scm 0 &\scm 0 &\scm 0 &\scm 0 &\scm 0 
&\scm 0 &\scm 0 &\scm 0 &\scm 0 &\scm 0 
&\scm 0 &\scm 1 &\scm 0 &\scm 0 &\scm 0 
&\scm 0 &\scm 0 &\scm 0 &\scm 0 &\scm 0 \\ 
 \scm 0 &\scm 0 &\scm 0 &\scm 0 &\scm 0 
&\scm 0 &\scm 0 &\scm 0 &\scm 0 &\scm 0 
&\scm 0 &\scm 0 &\scm 0 &\scm 0 &\scm 0 
&\scm 0 &\scm 0 &\scm 1 &\scm 0 &\scm 0 
&\scm 0 &\scm 0 &\scm 0 &\scm 0 &\scm 0 \\ 
 \scm 0 &\scm 0 &\scm 0 &\scm 0 &\scm 0 
&\scm 0 &\scm 0 &\scm 0 &\scm 0 &\scm 0 
&\scm 0 &\scm 0 &\scm 0 &\scm 0 &\scm 0 
&\scm 0 &\scm 0 &\scm 0 &\scm 1 &\scm 0 
&\scm 0 &\scm 0 &\scm 0 &\scm 0 &\scm 0 \\ 
 \scm 0 &\scm 0 &\scm 0 &\scm 0 &\scm 0 
&\scm 0 &\scm 0 &\scm 0 &\scm 0 &\scm 0 
&\scm 0 &\scm 0 &\scm 0 &\scm 0 &\scm 0 
&\scm 0 &\scm 0 &\scm 0 &\scm 0 &\scm 1 
&\scm 0 &\scm 0 &\scm 0 &\scm 0 &\scm 0 \\ 
 \scm 0 &\scm  0 &\scm   0 &\scm 0 &\scm 0 
&\scm 0 &\scm  0 &\scm   0 &\scm 0 &\scm 0 
&\scm 0 &\scm  0 &\scm   0 &\scm 0 &\scm 0 
&\scm 5 &\scm  0 &\scm   0 &\scm 0 &\scm 0 
&\scm 0 &\scm  0 &\scm   0 &\scm 0 &\scm 0 \\\hline\rule{0mm}{2mm} 
 \scm 0 &\scm 0 &\scm 0 &\scm 0 &\scm 0 
&\scm 0 &\scm 0 &\scm 0 &\scm 0 &\scm 0 
&\scm 0 &\scm 0 &\scm 0 &\scm 0 &\scm 0 
&\scm 0 &\scm 0 &\scm 0 &\scm 0 &\scm 0 
&\scm 0 &\scm 1 &\scm 0 &\scm 0 &\scm 0 \\ 
 \scm 0 &\scm 0 &\scm 0 &\scm 0 &\scm 0 
&\scm 0 &\scm 0 &\scm 0 &\scm 0 &\scm 0 
&\scm 0 &\scm 0 &\scm 0 &\scm 0 &\scm 0 
&\scm 0 &\scm 0 &\scm 0 &\scm 0 &\scm 0 
&\scm 0 &\scm 0 &\scm 1 &\scm 0 &\scm 0 \\ 
 \scm 0 &\scm 0 &\scm 0 &\scm 0 &\scm 0 
&\scm 0 &\scm 0 &\scm 0 &\scm 0 &\scm 0 
&\scm 0 &\scm 0 &\scm 0 &\scm 0 &\scm 0 
&\scm 0 &\scm 0 &\scm 0 &\scm 0 &\scm 0 
&\scm 0 &\scm 0 &\scm 0 &\scm 1 &\scm 0 \\ 
 \scm 0 &\scm 0 &\scm 0 &\scm 0 &\scm 0 
&\scm 0 &\scm 0 &\scm 0 &\scm 0 &\scm 0 
&\scm 0 &\scm 0 &\scm 0 &\scm 0 &\scm 0 
&\scm 0 &\scm 0 &\scm 0 &\scm 0 &\scm 0 
&\scm 0 &\scm 0 &\scm 0 &\scm 0 &\scm 1 \\ 
 \scm 0 &\scm  0 &\scm   0 &\scm 0 &\scm 0 
&\scm 0 &\scm  0 &\scm   0 &\scm 0 &\scm 0 
&\scm 0 &\scm  0 &\scm   0 &\scm 0 &\scm 0 
&\scm 0 &\scm  0 &\scm   0 &\scm 0 &\scm 0 
&\scm 5 &\scm  0 &\scm   0 &\scm 0 &\scm 0 \\
\ea
\right]}
\;\; =: \;\; \ddot{v}\; .
\]

We observe that we may replace $\dot u$, $\dot v$ by $\ddot u$, $\ddot v$ resp.\ by $\dddot u$, $\dddot v$ to obtain
\[
\sLD\= \Lambda((\ddot u,\ddot v), (4,4), (2,11))_{\dot u\Lambda} 
   \= \{f\in\Lambda\; :\; D_{\ddot u}^i\circ D_{\ddot v}^j(f) \con_{\dot u^{2i + 11j}\Lambda} 0\mb{\rm\ for all } i,j\in [0,4] \} 
\]
and 
\[
\sLDE \= 
\{ f\in\sLD\; :\; D_{\;\dddot u}^i\circ E_{\;\dddot u\, ,\, \dddot v}^j(f) \con_{\dot u^{2i + 13j}\Lambda} 0\mb{\rm\ for all } i,j\in [0,4]\}\; .
\]

\subsection{A spectral sequence}

\bq
Alternatively, there is a Lyndon-Hochschild-Serre-Grothendieck spectral sequence that might perhaps help in calculating cohomology in the case 
$C_{p^2}$ using the result in the case $C_p$ instead of using the Wedderburn embedding of $U\wr C_{p^2}$ (cf.\ preceding sections). 
Due to varying ground rings, we have to apply a (hardly visible) modification to the usual Lyndon-Hochschild-Serre spectral sequence. 
\eq

Let $S\tm T\tm U$ be an iterated finite extension of discrete valuation rings, $U|S$ galois with $H = \Gal(U|S)$, $T|S$ galois with 
$G = \Gal(T|S)$. Let $N$ be the kernel of the restriction map $H\lra G$, so $N = \Gal(U|T)$ and $G \iso H/N$. In this section, modules are 
not necessarily finitely generated.

The Grothendieck spectral sequence of the composition
\[
\Modr U\wr H \;\mra{\liu{U\wr N}{(U,-)}}\; \Modr T\wr G \;\mra{\liu{T\wr G}{(T,-)}}\; \Modr S
\]
is given by
$$
E_2^{m,n} \; :=\; \Ext^m_{T\wr G}(T,\Ext^n_{U\wr N}(U,X))\; ,
\leqno (\ast)
$$
where $X\in\modr U\wr H$, $m,\, n\,\geq\, 0$. For $X\in\Modr U\wr H$, the $T\wr G$-module structure on the image $\liu{U\wr N}{(U,X)}$ is induced 
by the $\liu{U\wr N}{(U,X)} \iso \liu{TN}{(T,X)}$ and the left $T\wr G$-module structure on $T$.

To prove that it converges to $\Ext^{m+n}_{U\wr H}(U,X)$, it suffices to show that an 
injective $U\wr H$-module $I$ is mapped to an injective module. In fact, for $Y\in\Modr T\wr G$ we calculate
\[
\barcl
\liu{T\wr G}{(Y,\liu{U\wr N}{(U,I)})}
& \iso & \liu{T\wr G}{(Y,\liu{TN}{(T,I)})} \\
& \iso & \liu{TN}{(Y\ts_{T\wr G} T,I)} \\
& \iso & \liu{TN}{(Y\ts_{T\wr G} T,\liu{U\wr H}(U\wr H, I))} \\
& \iso & \liu{U\wr H}{(Y\ts_{T\wr G} T\ts_{TN} U\wr H, I)}\; , \\
\ea
\]
so that the assertion follows by injectivity of $I$ and by projectivity as a left $T\wr G$-module of 
\[
\ba{rcl}
T\ts_{TN} U\wr H
& \liu{T\wr G\,}{\iso} & T\ts_{TN} (T\wr H)^{(\# N)} \\
& \liu{T\wr G\,}{\iso} & (T\wr G)^{(\# N)}\; . \\
\ea
\]

Using adjunction, we may rewrite $(\ast)$ in the familiar shape
$$
E_2^{m,n} \; =\; \HH^m(G,\HH^n(N,X;T);S) \;\imp\; \HH^{m+n}(H,X;S)\; ,
\leqno (\ast\ast)
$$
applicable to $X\in\Modr U\wr H$ --- so e.g.\ to $X = U$. Concerning the cup product, cf.\ {\bf\cite[\rm sec.\ 3.9]{Be91}}. 

Now, if $H = C_{p^2}$, $N = C_p$, $G = C_p$ and $X = U$, and our remaining conditions are satisfied (pure ramification, $\val(p)$ big enough), then
(\ref{PropAG8}) already calculates $E_2^{m,0}$ for $m\geq 0$. The first step to take when pursuing this spectral sequence approach, using $(\ast)$ 
rather than $(\ast\ast)$, would be to calculate $\Ext^n_{U\wr N}(U,U)$ {\it as a $T\wr G$-module} for $n\geq 1$. We do not know whether this 
approach is actually viable.
\end{footnotesize}         
\parskip0.0ex
\begin{footnotesize}

\parskip1.2ex

\vspace*{1cm}

\begin{flushright}
Matthias K\"unzer\\
German University of Cairo\\
New Cairo City\\
Main Entrance El Tagamoa El Khames\\
Egypt\\
Matthias.Kuenzer@guc.edu.eg\\
\vspace*{1cm}
Harald Weber\\
Universit\"at Stuttgart\\
Institut f\"ur Algebra und Zahlentheorie\\
Pfaffenwaldring 57\\
D-70569 Stuttgart\\
Harald.Weber@iaz.uni-stuttgart.de\\
\end{flushright}
\end{footnotesize}

\end{document}